%% file: root_sa_main.tex
\documentclass[11pt,twoside]{article}
\usepackage{fullpage}
\usepackage{epsf}
\usepackage{fancyhdr}
\usepackage{graphics}
\usepackage{graphicx}
\usepackage{psfrag}
\usepackage{microtype}
\usepackage{subfigure}
\usepackage{algorithmic} \usepackage{color,xcolor}
\usepackage{subfiles}
\usepackage{color} \usepackage{tabularx} \usepackage{amsthm}
\usepackage{amsfonts} \usepackage{amsmath} \usepackage{amssymb,bbm}
\usepackage{stackengine} \usepackage{footnote}
\makesavenoteenv{tabular} \makesavenoteenv{table}
\usepackage{algorithm}

\input{final_macros}

\usepackage{url}

\usepackage[colorlinks=True,linkcolor=magenta,citecolor=blue,urlcolor=blue,backref=true]{hyperref}

\usepackage{nicefrac} \usepackage{comment}

\usepackage{chngpage}

 \usepackage{tabularx}%

\usepackage{enumitem}
\usepackage{booktabs}
\usepackage{caption}

\usepackage{bm,bbm}
\usepackage{mathtools}

\newtheorem{assumption}{Assumption}

\usepackage[mathscr]{euscript} 

\theoremstyle{definition} 
\makeatletter \long\def\@makecaption#1#2{
        \vskip 0.8ex \setbox\@tempboxa\hbox{\small {\bf #1:} #2}
        \parindent 1.5em 
        \dimen0=\hsize \advance\dimen0 by -3em \ifdim \wd\@tempboxa
        >\dimen0 \hbox to \hsize{ \parindent 0em \hfil
                        \parbox{\dimen0}{\def\baselinestretch{0.96}\small
                          {\bf #1.} #2
                                }\hfil} \else \hbox to \hsize{\hfil
          \box\@tempboxa \hfil} \fi } \makeatother

\def\ROOTSGD{\texttt{ROOT-SGD}}

  \newcommand{\stationary}{\xi} \newcommand{\avgcost}{\mu}
  \newcommand{\mixingtime}{t_{\mathrm{mix}}}
  \newcommand{\totalvariation}{d_{\mathrm{TV}}}

\newcommand{\myrhotil}{\ensuremath{\plainrhotil_\numobs}}

 \newcommand{\rewardMax}{r_{\max}}
\newcommand{\costMax}{c_{\max}}

\newcommand{\DudGauss}{\ensuremath{\dudley_2}}
\newcommand{\DudExp}{\ensuremath{\dudley_1}}

\newcommand{\UpsSet}{\ensuremath{\mathcal{C}}}

\newcommand{\myorder}{\ensuremath{\mathcal{O}}}

\newcommand{\MySet}{\ensuremath{C}}

\newcommand{\delrad}{\ensuremath{s}}
\newcommand{\delradstar}{\ensuremath{\delrad^*_\numobs}}
\newcommand{\delradastar}{\ensuremath{\delrad^*_{\MySet, \numobs}}}

\newcommand{\QOp}{\mathbf{h}} 
 \newcommand{\Qstar}{\theta^\star}
\newcommand{\rewardQ}{r} \newcommand{\Q}{\theta}
\newcommand{\infNorm}[1]{\| #1 \|_\infty} 
 \newcommand{\AQ}{\contraction
  \TransMat^{\pi_\Q}} \newcommand{\Astar}{\contraction
  \TransMat^{\pi_\star}} \newcommand{\statesQ}{\mathcal{S}}
\newcommand{\actionsQ}{\mathcal{A}}

\newcommand{\plainrhotil}{\ensuremath{\rho}}

\newcommand{\localgausssemi}{\ensuremath{\localgauss_\MySet}}
\newcommand{\localvarsemi}{\ensuremath{\localvar_{\MySet } }}

\newcommand{\Martin}{\ensuremath{M}}
\newcommand{\MartinStar}{\ensuremath{\Martin^*}}

\newcommand{\noisegausssemi}{\ensuremath{\noisegauss_\MySet}}
\newcommand{\noisevarsemi}{\ensuremath{\noisevar_\MySet}}

\newcommand{\Qmat}{\ensuremath{Q}}
\newcommand{\QmatTil}{\ensuremath{\widetilde{\Qmat}}}

\newcommand{\mysqrt}[1]{\ensuremath{\big( #1 \big)^{1/2}}}
\newcommand{\mybiggersqrt}[1]{\ensuremath{\Big( #1 \Big)^{1/2}}}
\newcommand{\mynorm}[1]{\ensuremath{ \| #1 \|}}
\newcommand{\mybignorm}[1]{\ensuremath{ \big\| #1 \big\|}}
\newcommand{\mybiggernorm}[1]{\ensuremath{ \Big\| #1 \Big\|}}

\newcommand{\myinfnorm}[1]{\ensuremath{ \| #1 \|_{\infty}}}

\newcommand{\sqrtlog}[1]{\ensuremath{\sqrt{\log{#1}}}}

\newcommand{\MuCon}{\ensuremath{\mu}}

\newcommand{\myV}{v} \newcommand{\myZ}{z}
\newcommand{\lecamRad}{\Delta}
\newcommand{\admissiblePar}{\kappa}

\newcommand{\myspan}{\ensuremath{\mathrm{span}}}
\newcommand{\avgcoststar}{\ensuremath{\avgcost_*}}

\newcommand{\HOT}{\ensuremath{\mathcal{H}}}
\newcommand{\FullHot}{\ensuremath{\HOT_\numobs(\delta, \stepsize)}}

\newcommand{\SuperHot}[1]{\ensuremath{\HOT^{#1}_\numobs( \delta, \stepsize)}}

\newcommand{\NewHot}{\SuperHot{\dagger}}

\newcommand{\ThreeHot}{\SuperHot{\circ}}

\newcommand{\myunder}[2]{\ensuremath
  \underbrace{#1}_{\mbox{\footnotesize{#2}}}}

\newcommand{\maxvar}{\ensuremath{\fullmaxvar{\dualBall}}}
\newcommand{\maxstd}{\ensuremath{\fullmaxstd{\dualBall}}}

\newcommand{\fullmaxvar}[1]{\ensuremath{\sigma^2_{#1}}}
\newcommand{\fullmaxstd}[1]{\ensuremath{\sigma_{#1}}}

\newcommand{\userho}{\ensuremath{\rho}}

 \usepackage{fullpage}

\usepackage{cleveref}

\usepackage[firstinits=true,
  maxnames = 99,
  backref=true,
  backend=bibtex,
  style=alphabetic,
  sortlocale=de_DE,
  natbib=true,
  doi=false,
  isbn=false,
  url=false]{biblatex}
\newenvironment{carlist}
 {\begin{list}{$\bullet$}
 {\setlength{\topsep}{0in} \setlength{\partopsep}{0in}
  \setlength{\parsep}{0in} \setlength{\itemsep}{\parskip}
  \setlength{\leftmargin}{0.15in} \setlength{\rightmargin}{0.08in}
  \setlength{\listparindent}{0in} \setlength{\labelwidth}{0.08in}
  \setlength{\labelsep}{0.1in} \setlength{\itemindent}{0in}}}
 {\end{list}}

\newcommand{\bcar}{\begin{carlist}}
\newcommand{\ecar}{\end{carlist}}

\newcommand{\myorderlog}{\ensuremath{\widetilde{\myorder}}}

\newcommand{\newvar}{\ensuremath{\localvar}}

\newcommand{\weightmax}{\ensuremath{w_{\mbox{\scriptsize{max}}}}}
\newcommand{\weightmin}{\ensuremath{w_{\mbox{\scriptsize{min}}}}}

\newcommand{\PLAINHOT}{\ensuremath{\mathcal{H}_\numobs}}

\begin{document}

\begin{center}
{\bf{\LARGE{Optimal variance-reduced stochastic approximation in
      Banach spaces}}}

\vspace*{.2in} {\large{
 \begin{tabular}{ccc}
  Wenlong Mou$^{\star, \diamond}$ & Koulik Khamaru$^{\star, \dagger}$
  & Martin J. Wainwright$^{\diamond, \dagger, \circ}$\\
 \end{tabular}
 \begin{tabular}
 {cc} Peter L. Bartlett$^{\diamond, \dagger}$ & Michael
 I. Jordan$^{\diamond, \dagger}$
 \end{tabular}

}}

\vspace*{.2in}

 \begin{tabular}{c}
 Department of Electrical Engineering and Computer
 Sciences$^\diamond$\\ Department of Statistics$^\dagger$ \\ UC
 Berkeley\\
 \end{tabular}

\vspace*{.2in}

 \begin{tabular}{c}
    Department of Electrical Engineering and Computer Sciences$^\circ$
    \\ Department of Mathematics$^\circ$ \\ Massachusetts Institute of
    Technology
 \end{tabular}

   \vspace*{.2in}

   \today

\vspace*{.2in}

\begin{abstract}
We study the problem of estimating the fixed point of a contractive
operator defined on a separable Banach space. Focusing on a stochastic
query model that provides noisy evaluations of the operator, we
analyze a variance-reduced stochastic approximation scheme, and
establish non-asymptotic bounds for both the operator defect and the
estimation error, measured in an arbitrary semi-norm.  In contrast to
worst-case guarantees, our bounds are instance-dependent, and achieve
the local asymptotic minimax risk non-asymptotically. For linear
operators, contractivity can be relaxed to multi-step contractivity,
so that the theory can be applied to problems like average reward
policy evaluation problem in reinforcement learning.  We illustrate
the theory via applications to stochastic shortest path problems,
two-player zero-sum Markov games, as well as average-reward policy
evaluation.
\let\thefootnote\relax\footnotetext{$\star$ Wenlong Mou and Koulik
  Khamaru contributed equally to this work.}
\end{abstract}

\end{center}


\section{Introduction}
\label{new--SecIntro}

In this paper, we consider a class of stochastic fixed-point problems
defined in Banach spaces.  In particular, let $\vecspace$ be a
separable Banach space with its associated norm $\norm{\cdot}$, and
suppose that $\hpop: \vecspace \rightarrow \vecspace$ is an operator
on the Banach space.  Of interest to us are solutions $\thetastar$ to
the fixed-point equation
\begin{align}
\label{eq:fixed-point}  
\thetastar = \hpop(\thetastar).
\end{align}
When the operator $\hpop$ is contractive, the Banach fixed point
theorem (e.g.,~\cite{Dugundji}) ensures the existence and uniqueness
of the fixed point.  The bulk of our analysis focuses on this
contractive case, but we also allow for weaker multi-stage contraction
in certain settings.

Fixed points of this type lie at the core of many mathematical areas,
including differential and integral
equations~\cite{teschl2012ordinary,kirsch2011introduction}, game
theory~\cite{stokey1989recursive}, optimization and variational
inequalities~\cite{nesterov2003introductory,
  rockafellar2009variational}, as well as dynamic programming and
reinforcement learning~\cite{bertsekas2019reinforcement,Puterman05}.
In these settings, the contraction property not only plays an
instrumental role in existence and uniqueness proofs, but also leads
to efficient methods for computing fixed points. Our focus will be on
the extension of such methods to problems in which the operator
$\hpop$ can be observed only via a stochastic oracle that, when given
a query point $\theta$, returns a noisy version of the operator
evaluation $\hpop(\theta)$. Such random observation models necessitate
the use of stochastic approximation schemes. A fundamental question
associated with such schemes is their \emph{statistical complexity}:
how many noisy operator evaluations are required to estimate the fixed
point $\thetastar$ to a pre-specified accuracy? In this paper, we
undertake a fine-grained yet relatively general analysis of this
question.  Notably, our analysis captures the way in which statistical
complexity depends on the geometry of the Banach space, as well as the
structure of the fixed point $\thetastar$ itself.

An important sub-class of Banach spaces are Hilbert spaces, with the
Euclidean case ($\vecspace = \real^d$ with the usual inner product)
being one special example. The behavior of stochastic approximation
for many Hilbert spaces is relatively well understood. In this case,
the space $\vecspace$ is endowed with an inner product
$\inprod{\cdot}{\cdot}_\vecspace$ that induces the norm $\norm{x} =
\sqrt{\inprod{x}{x}_\vecspace}$. For example, for the Euclidean space
$(\real^\usedim, \| \, \cdot \, \|_2)$, if we set \mbox{$\hpop(x)
  \mydefn x - \beta^{-1} \nabla f (x)$} for a \mbox{$\beta$-smooth}
and strongly convex function $f$, then solving the fixed-point
equation~\eqref{eq:fixed-point} is equivalent to minimizing the
function $f$. A rich theory has been developed around this stochastic
optimization
problem~\cite{bottou2018optimization,nemirovski2009robust}, giving
rise to the concepts of
averaging~\cite{polyak1992acceleration,ruppert1988efficient},
acceleration~\cite{ghadimi2012optimal,ghadimi2013optimal}, and
variance
reduction~\cite{johnson2013accelerating,nguyen2021inexact,li2020root},
along with associated characterizations of
optimality~\cite{moulines2011non,duchi2016asymptotic,mou2020linear}.

In contrast, relatively less is known in the general setting of Banach
spaces.  One of the simplest examples is $\real^d$ equipped with a
non-Euclidean norm, such as the $\ell_\infty$-norm.  To be clear,
non-Euclidean set-ups of this type have been studied in the literature
on stochastic optimization and stochastic variational inequalities,
with the method of mirror descent being a representative
example~\cite{nemirovski1983problem,juditsky2011solving,kotsalis2020simple}.
Our study, however, deviates from this line of research. The
difference stems from the formulation of the problem itself: the
operator $\hpop$ in equation~\eqref{eq:fixed-point} is a mapping from
$\vecspace$ to itself, whereas the operators studied in variational
inequalities map a Banach space to its dual. This difference leads to
a different path of analysis, as taken here.

At least initially, it might seem that non-Euclidean geometry should
pose little difficulty for stochastic approximation: all norms are
equivalent in the finite-dimensional case, and as is known from
standard theory (e.g.,~\cite{borkar2009stochastic}), asymptotic
convergence depends ultimately on the limiting ODE defined by the
scheme.  From a non-asymptotic point of view, however, the picture
becomes more nuanced: a natural desideratum is that the bounds depend
on the \emph{geometric complexity} of $\vecspace$, as opposed to its
(possibly much larger) ambient dimension.  The difference between the
two can be significant. As one concrete example, when solving
fixed-point equations that arise in tabular Markov decision processes,
the ambient dimension is the size of state-action space, whereas one
can obtain $\ell_\infty$-norm bounds that have only logarithmic
dependency on the dimension (see,
e.g.,~\cite{wainwright2019stochastic,khamaru2020temporal}). Our first
goal, therefore, is to develop a unified and geometry-aware theory for
a certain class of stochastic approximation procedures in Banach
spaces.

Our second goal is to establish bounds that are instance-dependent,
and so move us beyond a classical worst-case analysis.  Any method for
stochastic approximation corresponds to a particular type of recursive
statistical estimator, so that the the classical statistical theory of
local asymptotic minimax can be brought into
play~\cite{hajek1972local,van2000asymptotic}.  This theory provides a
framework for deriving lower bounds on the error of any estimator that
depend explicitly on (a local neighborhood of) the instance under
consideration.  As for the form that such bounds should take in our
setting, recall that a sum of $\mathrm{i.i.d.}$ random variables in
Banach spaces is known (under mild regularity conditions) to satisfy a
central limit theorem (see Ledoux and
Talagrand~\cite{ledoux2013probability}, Section 10). These two lines
of asymptotic analysis, in conjunction, indicate that the ``right''
complexity for estimation in a Banach space $\vecspace$ should involve
the expected norm of a Gaussian random element with covariance
structure specified by the noise in the stochastic oracle.  Given this
fundamental limit, it is natural to seek an estimator whose
\emph{non-asymptotic} risk matches this quantity, with possible
higher-order terms which, again, depends only the geometric complexity
of the norm $\norm{\cdot}$ (and not the ambient dimension).

In order to address these goals, we analyze an extension of the
ROOT-SGD algorithm, a stochastic approximation (SA) algorithm
introduced in past work involving a subset of the current
authors~\cite{li2020root}. We adapt the scheme to solve general
fixed-point problems and establish instance-dependent non-asymptotic
guarantees in general Banach spaces. More specifically:
\bcar
\item We establish sharp non-asymptotic bounds on the \emph{operator
defect} $\norm{\hpop (\theta_\numobs) - \theta_\numobs}$ of the
  iterate $\theta_\numobs$ after $\numobs$ rounds.  The leading-order
  term, defined in terms of a Gaussian complexity induced by the noisy
  evaluations of the operator $\hpop$, matches the the optimal
  Gaussian limit.  To the best of our knowledge, this is the first
  non-asymptotic bound for SA procedures with general non-Euclidean
  norm that depends directly on the geometric complexity of the
  underlying space.
\item Under a local linearization assumption on the operator $\hpop$,
  we establish a sharp instance-dependent upper bound on the
  \emph{estimation error} $\anorm{\theta_\numobs - \thetastar}$,
  measured by any semi-norm $\anorm{\cdot}$ that is dominated by
  $\norm{\cdot}$. The leading-order term of this bound is a Gaussian
  complexity involving the dual ball of the semi-norm $\anorm{\cdot}$,
  and its interaction with locally linear approximations of the
  operator around $\thetastar$.
\item When the operator $\hpop$ is affine, we establish an improved
  result that matches the leading-order term in the nonlinear case,
  and with an even lower sample complexity.  We also generalize this
  result to settings in which $\hpop$ itself is not necessarily
  contractive, but its $\compostep$-step composition is
  contractive.
\item Finally, we illustrate some specific consequences of our theory
  for different examples, including stochastic shortest path problems,
  Markov games, and average-reward policy evaluation.
\ecar  

\subsection{Related work}
In this section, we survey existing literature on stochastic
approximation and its variance-reduced analogues.

\paragraph{Stochastic approximation and asymptotic guarantees:}

The study of stochastic approximation methods dates back to the
seminal work of Robbins and Monro~\cite{robbins1951stochastic}, as
well as Kiefer and Wolfowitz~\cite{kiefer1952stochastic}, who
established asymptotic convergence for various classes of
one-dimensional problems. Subsequent work by
Ljung~\cite{ljung1977analysis,ljung1977positive} and Kushner and
Clark~\cite{kushnerclarks} provided general criteria for convergence
to a stable limit, in particular by using an ordinary differential
equation (ODE) to track the trajectory of SA procedures.  The ODE
method has been substantially refined in a long line of subsequent
work~\cite{benaim1996dynamical,kushner2003stochastic,borkar2009stochastic,benveniste2012adaptive}.
In addition to pointwise convergence, there is a rich body of work
characterizing the asymptotic distribution of SA
trajectories~\cite{khas1966stochastic,
  kushner1984approximation,kushner1984invariant}.  We refer the reader
to the
monographs~\cite{kushner2003stochastic,borkar2009stochastic,benveniste2012adaptive}
for more background and details on these results.

The idea of improving SA schemes by averaging the iterates was
proposed in independent work by Polyak and
Juditsky~\cite{polyak1990new,polyak1992acceleration} as well as
Ruppert~\cite{ruppert1988efficient}.  Averaging the iterates allows
for the use of more aggressive stepsize choices, and Gaussian limiting
behavior is achieved over a broad range.  The form of this limiting
distribution is known to optimal in the sense of local asymptotic
minimax~\cite{hajek1972local,van2000asymptotic,duchi2016asymptotic}. The
idea of iterate averaging underlies many important aspects of
large-scale statistical learning, leading to improved algorithms in
different
settings~\cite{bach2013non,duchi2016asymptotic,tripuraneni2018averaging}
and laying the foundations of online statistical
inference~\cite{chen2020statistical}. The ROOT-SGD
algorithm~\cite{li2020root} that inspired our approach is motivated by
the averaging scheme, but combines variance reduction with averaging
of the gradient sequence (as opposed to the sequence of iterates).


\paragraph{Non-asymptotic guarantees for stochastic approximation:}

Recent years have witnessed significant interest in obtaining
non-asymptotic guarantees of the standard SA scheme (see
equation~\eqref{eqn:Vanilla-SA} in the sequel).  For instance, Qu and
Wierman~\cite{qu2020finite} directly analyzed the iterates of SA
algorithms in the asynchronous setting, whereas Chen et
al.~\cite{chen2020finite} derived non-asymptotic bounds on stochastic
approximation methods using Lyapunov functions.  Using the generalized
Moreau envelope, they constructed a smooth Lyapunov function, and show
that the iterates of a standard SA scheme have a negative drift with
respect to this Lyapunov function.  Such Lyapunov techniques have been
used to derive non-asymptotic guarantees for SA schemes in variety of
settings (e.g.,~\cite{chen2021lyapunov,
  chen2019performance,chen2021finite,zhang2021finite}).
Wainwright~\cite{wainwright2019stochastic} proved non-asymptotic
guarantees for stochastic approximation algorithms under a
cone-contractive assumption.  For general contractive fixed-point
problems in Banach spaces, Gupta et al.,~\cite{gupta2018probabilistic}
developed general criteria for the asymptotic convergence of
mini-batch fixed-point iterations; and recently,
Borkar~\cite{borkar2021concentration} established non-asymptotic
concentration inequalities for the iterates, albeit with potentially
dimension-dependent pre-factors.  It should be noted that the standard
SA scheme~\eqref{eqn:Vanilla-SA}, while guaranteed to converge to the
fixed point, may do so at a sub-optimal rate when measured in a
minimax sense; for example, the
papers~\cite{wainwright2019stochastic,li2021q} demonstrate the
non-optimality of this approach for the $Q$-learning problem in
reinforcement learning.

Non-asymptotic guarantees that are instance-dependent---meaning that
they go beyond worst-case and are adaptive to the difficulty---have
been established for several stochastic approximation procedures.  For
stochastic gradient (SG) methods in the Euclidean setting, such bounds
have been established for Polyak-Ruppert-averaged
SG~\cite{moulines2011non,gadat2017optimal} and variance-reduced SG
algorithms~\cite{frostig2015competing,li2020root}, with the sample
complexity and high-order terms being improved over time. For
reinforcement learning problems, such type of guarantees have been
established in the $\vecnorm{\cdot}{\infty}$ norm for temporal
difference methods~\cite{khamaru2020temporal} and
$Q$-learning~\cite{khamaru2021instance} under a generative model, as
well as Markovian trajectories~\cite{mou2021optimal,li2021accelerated}
under the $\ell_2$-norm. In the context of stochastic optimization,
the paper~\cite{li2020root} provides fine-grained bound for
\textsf{ROOT-SGD} with a \emph{unity} pre-factor on the leading-order
instance-dependent term. The bounds in our paper, on the other hand,
involve constants that need not be optimal in this sense.  It is an
interesting future direction of research to establish similar
non-asymptotic bounds for \rootSA with the sharp unity pre-factor.

\paragraph{Variance-reduced stochastic approximation algorithms:}

In order to obtain optimal SA procedures, different forms of variance
reduction have been analyzed. The idea of variance reduction in
stochastic approximation is classical; in the specific context of
stochastic gradient methods, the
papers~\cite{johnson2013accelerating,defazio2014saga,schmidt2017minimizing}
proposed versions of variance reduction that accelerate convergence by
careful averaging and re-centering of the gradient sequence.  In this
special case of stochastic optimization, the fixed-point
operator~$\hpop$ is obtained from the gradient update operator
(cf.\ the discussion in Section~\ref{new--SecIntro}); under suitable
convexity and smoothness conditions, it is contractive under the
$\ell_2$-norm.  In more recent work, several fully online schemes for
variance-reduced stochastic optimization have been developed and
analyzed, including SARAH~\cite{nguyen2017sarah,nguyen2021inexact},
STORM~\cite{cutkosky2019momentum} and
\textsf{ROOT-SGD}~\cite{li2020root}.  The \textsf{ROOT-SGD} scheme
uses recursive $1/t$-averaging of gradients, and has been been shown
to be optimal for various convex problems in both asymptotic and
non-asymptotic settings; see the paper~\cite{li2020root} and
references therein for more details.

In the context of reinforcement learning (RL) problems, the operator
$\hpop$ often corresponds to some type of Bellman
operator~\cite{bertsekas2012weighted,bertsekas2019reinforcement},
known to be contractive under the \mbox{$\ell_\infty$-norm.}
Unfortunately, the key techniques used to design optimal methods for
RL differ considerably from those used in the stochastic optimization
literature.  Concretely, in order to obtain optimal RL algorithms, it
is often necessary to exploit monotonicity properties of the Bellman
operator, combined with variance reduction
schemes~\cite{sidford2018near,wainwright2019variance,khamaru2020temporal,khamaru2021instance}.
Consequently, the literature is currently lacking a more unified
perspective on how to obtain optimal SA schemes in a general setting.
The main contribution of our paper is to fill this gap by proposing
and analyzing a single variance-reduced stochastic approximation
algorithm for finding the fixed point of \emph{any contractive}
operator. In this way, our analysis does not depend on the exact form
of the contraction norm $\norm{\cdot}$.

\paragraph{Notation:} We use $\vecspace^*$
to denote the dual space of the Banach space $\vecspace$, i.e., the
space of all bounded linear functionals on $\vecspace$. We define the
dual norm \mbox{$\dualnorm{y} \mydefn \sup_{x \in \vecspace \backslash
    \{0\}} \inprod{x}{y}/\norm{x}$.}  We define the unit norm ball
$\ball \mydefn \big \{x \in \vecspace, \norm{x} \leq 1 \big \}$ in
$\vecspace$, as well the dual norm unit ball \mbox{$\ball^* \mydefn
  \big \{ y \in \vecspace^* \mid \dualnorm{y} \leq 1 \big \}$.}

Given a bounded linear operator $A: \vecspace \rightarrow \vecspace$,
the adjoint operator $A^*: \vecspace^* \rightarrow \vecspace^*$ is
characterized by the property
\begin{align*}
  \inprod{A x}{y} = \inprod{x}{A^* y} \qquad \mbox{for all $x
    \in \vecspace$ and $y \in \vecspace^*$.}
\end{align*}
The operator norm of a bounded linear operator $A$ on $\vecspace$ is
given by \mbox{$\matsnorm{A}{\vecspace} \mydefn \sup_{x \in \vecspace
    \backslash \{0\}} \tfrac{\norm{A x}}{\norm{x}}$.}  Similarly, we
can define the operator norm $\matsnorm{\cdot}{\vecspace^*}$ of a
bounded linear operator mapping from $\vecspace^*$ to itself.  For any
bounded linear operator $A$ that maps $\vecspace$ to itself, we have
the equivalence $\matsnorm{A^*}{\vecspace^*} =
\matsnorm{A}{\vecspace}$.


\section{Problem set-up and the \rootSA Algorithm}
\label{SecProblemSetup}

In this section, we begin with a precise description of the class of
problems that we study, along with the assumptions imposed.  We then
describe the \rootSA algorithm analyzed in this paper.

\subsection{Problem formulation}

Consider a separable Banach space $(\vecspace, \norm{\cdot})$, and an
operator $\hpop$ mapping from $\vecspace$ to itself.  Assuming
sufficient regularity to guarantee the existence and uniqueness of the
fixed-point $\thetastar$ of the operator $\hpop$, we study stochastic
approximation procedures for estimating the fixed point, i.e., for
approximately solving the equation $\hpop (\theta) = \theta$. 

In many practical applications, we may not have access to the operator
$\hpop$ itself; instead, at each time $t$, we have access to a
stochastic oracle $\Hstoch_t$ that, when queried at some $\theta
\in \vecspace$, returns a noisy version $\Hstoch_t(\theta)$ of the
operator evaluation $\hpop(\theta)$.  We impose the following
conditions on the stochastic operators $\{\Hstoch_t \}_{t \geq 1}$ and
the population operator $\hpop$:
\paragraph{Assumptions}
\begin{enumerate}[label=(A\arabic*)]
\item \label{new--assume-pop-contractive}
There is a scalar $\contraction \in [0,1)$ such that the operator
  $\hpop: \vecspace \rightarrow \vecspace$ is
  $\contraction$-contractive---viz.
  \begin{align*}
 \norm{\hpop (\theta_1) - \hpop (\theta_2)} \leq \contraction
 \norm{\theta_1 - \theta_2} \quad \mbox{for all $\theta_1, \theta_2
   \in \vecspace$.}
  \end{align*}

\item
\label{new--assume-sample-operator}
For each $t = 1, 2, \ldots$, the stochastic operator
$\Hstoch_t: \vecspace \mapsto \vecspace$ is almost surely (a.s.)
$\Lip$-Lipschitz:
\begin{align*}
  \norm{\Hstoch_t (\theta_1) - \Hstoch_t (\theta_2)} \leq \Lip
  \norm{\theta_1 - \theta_2} \qquad \mbox{a.s. for all $\theta_1,
    \theta_2 \in \vecspace$.}
\end{align*}
\item
\label{new--assume-noise-bound}
For any fixed $\theta \in \vecspace$, the noise variables $
\{\noise_t(\theta) \defn\Hstoch_t(\theta) - \hpop(\theta)\}_{t \geq
  1}$ are zero-mean and i.i.d., and $\norm{\noise_t(\thetastar)} \leq
\boundstar$ almost surely for all $t = 1, 2, \ldots$.
\end{enumerate}
A few remarks are in order. By the Banach fixed point theorem
(e.g.,~\cite{Dugundji}), the contractivity condition in
Assumption~\ref{new--assume-pop-contractive} ensures that $\hpop$ has
a unique fixed point $\thetastar$.  The bulk of our analysis imposes
Assumption~\ref{new--assume-pop-contractive}, with the exception of
Section~\ref{new--sec:multi-step-contraction}, where it is relaxed to
a multi-stage contraction assumption in the special case of linear
operators. Throughout this paper, we assume that $\contraction \geq
\tfrac{3}{4}$ for the ease of presentation. Note that this assumption
can be made without loss of generality, since an operator that is
$\contraction$-contractive for some $\contraction \in[0, 3/4)$ is also
  $3/4$-contractive.

Assumption~\ref{new--assume-sample-operator} requires the stochastic
operator $\Hstoch_t$ to be Lipschitz, with the associated constant
$\Lip$ allowed to be much larger than one---that is, there is no
requirement that $\Hstoch_t$ be contractive or non-expansive.  This
setup should be contrasted with past work on cone-contractive
operators~\cite{wainwright2019stochastic,wainwright2019variance} or
$\ell_\infty$-norm
contractions~\cite{khamaru2020temporal,khamaru2021instance}, in which
the stochastic operator $\Hstoch_t$ itself is required to be
contractive.  In the special case of stochastic optimization in
$\real^\usedim$, this type of sample-level Lipschitz condition is
widely used, especially for variance-reduced procedures
(cf.~\cite{johnson2013accelerating,nguyen2021inexact,li2020root}).

As for Assumption~\ref{new--assume-noise-bound}, it imposes bounds
only on the noise function when evaluated at the fixed point
$\thetastar$ of the operator $\hpop$. In conjunction with
Assumption~\ref{new--assume-sample-operator}, this bound implies that
$\norm{\noise_t(\theta)} \leq \boundstar + (\Lip+\contraction)
\norm{\theta - \thetastar}$, allowing the norm of the noise
$\noise_t(\theta)$ to grow linearly with $\norm{\theta -
  \thetastar}$. It is worth remarking that by using slightly more
involved concentration arguments, it is possible to relax the almost
sure bounds in Assumptions~\ref{new--assume-sample-operator}
and~\ref{new--assume-noise-bound}.
More precisely, it suffices to impose a $p^{th}$-moment
condition on all projections:
\begin{subequations}
\begin{align}
  \sup_{u \in \dualBall} \Exs \left[ \inprod{u}{\Hstoch_1 (\theta_1) -
      \Hstoch_1 (\theta_2)}^p \right] &\leq p! \cdot \Lip^p
  \norm{\theta_1 - \theta_2}^p  \qquad \mbox{for all $\theta_1,
  \theta_2 \in \vecspace$, and} \\
\sup_{u \in \dualBall} \Exs \left[ \inprod{u}{\noise_1(\thetastar)}^p
  \right] & \leq p! \cdot \boundstar^p,
\end{align}
\end{subequations}
for all $p \geq 2$.  Here $\dualBall$ is a skeleton set whose convex
hull generates the dual norm ball.

\subsection{The \rootSA algorithm}
\label{new--sec:rootSgd-algo}

Stochastic approximation algorithms are methods for solving
fixed-point equations based on noisy observations.  In the simplest of
such schemes, one starts with initial point $\theta_0$, and then
performs the recursive update
\begin{align}
  \label{eqn:Vanilla-SA}
  \theta_{t + 1} = \theta_{t} + \SAstep_t
  \left\{\Hstoch_{t}(\theta_t) - \theta_t\right\},
\end{align} 
where $\{\SAstep_t\}_{t \geq 0}$ is a sequence of positive stepsizes,
typically in the interval $(0,1)$.  At any given step $t$, conditioned
on $\theta_t$, the quantity $\Hstoch_t(\theta_t)$ is an unbiased
estimate of $\hpop(\theta_t)$, and the noise in the observation model
is given by $\Hstoch_t(\theta_t) - \hpop(\theta_t)$. Under the
contractivity assumptions~\ref{new--assume-pop-contractive} on the
operator~$\hpop$ and moment bounds on the observation noise
$\{\Hstoch_t(\theta_t) - \hpop(\theta_t)\}_{t\geq 1}$, the sequence
$\{\theta_t\}$ converges almost surely to the unique fixed point
$\thetastar$; moreover, the rate of convergence of $\theta_t$ to
$\thetastar$ is governed by the conditional variance of
$\Hstoch_t(\theta_t)$ around its conditional mean $\hpop(\theta_t)$.
See the standard
texts~\cite{kushner2003stochastic,borkar2009stochastic,benveniste2012adaptive}
for results of this type.

\begin{algorithm}
\caption{~\rootSA: A recursive SA algorithm}
\label{AlgRootSGD}
\begin{algorithmic}[1]
\STATE Given (a) Initialization $\theta_{0} \in \vecspace$, (b)
Burn-in $\burnin \geq 2$, and (c) stepsize $\stepsize > 0$
\FOR{$t = 1, \ldots, T$}
  \IF{$t \leq \burnin$} \STATE $\myV_t = \tfrac{1}{\burnin} \sum
  \limits_{t = 1}^\burnin \left\{\Hstoch_t(\thetainit) - \thetainit\right\}, \quad \text{and}
  \quad \theta_t = \thetainit$.  \ELSE \STATE
  $\myV_t = \big( \Hstoch_{t}(\theta_{t-1}) - \theta_{t - 1} \big) +
  \tfrac{t - 1}{t} \Big \{ \myV_{t - 1} - \big(
  \Hstoch_t(\theta_{t-2}) - \theta_{t - 2} \big) \Big \}$,
  \STATE $\theta_t =  \theta_{t - 1} + \stepsize
  \myV_t$.
\ENDIF \ENDFOR \RETURN $\theta_T$
\end{algorithmic}
\end{algorithm}

The goal of variance reduction is to improve the basic stochastic
approximation scheme~\eqref{eqn:Vanilla-SA} by replacing
$\Hstoch_t(\theta_t) - \theta_t$ with an alternative quantity $\myV_t$
that has lower variance. In this paper, we study a simple version of
such a variance-reduction scheme, as described in
Algorithm~\ref{AlgRootSGD}. Our algorithm is inspired by
Recursive-one-over-t SGD (\textsf{ROOT-SGD}) algorithm proposed and
analyzed in the past work~\cite{li2020root} involving a subset of the
current authors.  The \textsf{ROOT-SGD} algorithm was developed for
stochastic optimization; it exploits a two-time scale framework that
averages the gradient while performing variance reduction. Our \rootSA
algorithm extends this same idea to the more general setting of
stochastic approximation for fixed-point finding in Banach
spaces. While the algorithms are similar in spirit, the analysis in
this paper uses completely different techniques, since it applies to
the Banach-space setting for general operators, as opposed to the
Euclidean setting and gradient operators of convex functions. The key
technical difficulties lie in the absence of inner product structure.


\section{Main results}
\label{new--SecMain}

In this section, we state our main results and discuss some of their
consequences.  At a high level, our main results consist of various
non-asymptotic bounds on the behavior of the \rootSA algorithm in a
number of different (semi)-norms. In
Section~\ref{new--SecRootSgd-upper-bound}, we derive bounds on the
\emph{operator defect} $\norm{\hpop(\theta_t) - \theta_t}$, which
measures how far the $t^{th}$-iterate $\theta_t$ of
Algorithm~\ref{AlgRootSGD} is from being a fixed point of the
population operator~$\hpop$. In other settings, we are interested in
bounds on the \emph{estimation error} $\mynorm{\theta_t -
  \thetastar}$; accordingly, Section~\ref{sec:Estimation-error} is
devoted to such results, along with bounds on various kinds of
semi-norms.  Finally, in~\Cref{new--sec:multi-step-contraction}, we
discuss how to obtain refined results in the special case of linear
operators, for which the contractivity
assumption~\ref{new--assume-pop-contractive} can be relaxed.

Central to our bounds are the second-order properties of the
i.i.d.\ noise sequence $\Seq{\noise_t(\thetastar)}$.  In particular,
we let $\GaussNoise \in \vecspace$ be a zero-mean Gaussian random
element with covariance structure
\begin{align}
\label{EqnGaussCov}
\Exs \left[ \inprod{\GaussNoise}{y} \cdot \inprod{\GaussNoise}{z}
  \right] = \Exs \left[ \inprod{\noise_1(\thetastar)}{y} \cdot
  \inprod{\noise_1(\thetastar)}{z} \right] \quad \text{for all} \;\;
y, z \in \vecspace^*.
\end{align}
Various statistics of this Banach-space-valued random variable,
including its mean $\Exs[\norm{W}]$ and variance in certain
directions, specify the leading instance-dependent terms of our
results.

\subsection{Upper bounds on operator defect}
\label{new--SecRootSgd-upper-bound}

We begin by stating some non-asymptotic upper bounds on the so-called
operator defect $\norm{\hpop(\theta_t) - \theta_t}$, which measures
the error in the iterate $\theta_t$ as a fixed point.  As noted above,
the Gaussian element $W$ with covariance structure~\eqref{EqnGaussCov}
plays a central role.  In addition to the expected norm
$\Exs[\norm{W}]$, our result involves a certain type of maximal
variance over a skeleton set---namely, a subset $\dualBall$ of the
dual ball $\DualSpace$ such that $\ball^* = \mathrm{conv}
(\dualBall)$, so that the norm $\norm{\cdot}$ has the variational
representation $\norm{x} = \sup_{y \in \dualBall} \inprod{x}{y}$.
Given a set of this type, we define the $\dualBall$-maximal variance
\begin{align}
  \maxvar(W) \defn \sup_{u \in \dualBall} \Exs \big[ \inprod{u}{W}^2
    \big].
\end{align}
With this definition, we have the following:
\begin{theorem}
\label{new--thm:bellman-err-bound-general-nonlinear}
Under
Assumptions~\ref{new--assume-pop-contractive}---\ref{new--assume-noise-bound}
and a given failure probability $\delta \in (0,1)$, there is a range
of stepsizes $\stepsize > 0$ such that with burn-in period $\burnin
\defn \tfrac{c}{(1 - \contraction)^2 \stepsize} \log\big(
\tfrac{\numobs}{\delta}\big)$, sample size $\numobs \geq 2 \burnin$, the last iterate $\theta_\numobs$ of
Algorithm~\ref{AlgRootSGD} satisfies
\begin{align}
\label{new--eq:main-prop-bellman-err-bound-general-nonlinear}  
\norm{\hpop (\theta_\numobs) - \theta_\numobs} & \leq \myunder{
  \frac{c}{\sqrt{\numobs}}
  \Big \{\Exs[\norm{\GaussNoise}] + \sqrt{\maxvar(W)
    \log(\tfrac{1}{\delta})} \Big \}}{Instance-dependent} +
\myunder{\frac{c \burnin}{\numobs} \cdot \norm{\thetainit - \hpop
    (\thetainit)}}{Initial error} + \myunder{\FullHot}{Higher-order}
\end{align}
with probability at least $1-\delta$.  Here $\FullHot$ is a
higher-order term defined below (cf. equation~\eqref{EqnDefnFullHot}).
\end{theorem}
\noindent See
Section~\ref{subsec:proof-prop-bellman-err-general-nonlinear} for a
proof of this theorem.

\medskip

\paragraph{Instance-dependent term:}  As we discuss in the sequel,
with a stepsize $\stepsize \asymp 1/\sqrt{\numobs}$, the dominant
quantity in this upper bound is the instance-dependent term defined by
the Gaussian process $W$.  So as to appreciate its significance, we
note that for any $\delta \in (0,1)$, a near-optimal tail bound for
the Gaussian process $W$ is given by
\begin{align*}
\Prob \Big[ \norm{W} \geq \Exs[\norm{W}] + c \, \sqrt{\maxvar(W)
    \log(\tfrac{1}{\delta})} \; \Big] \leq \delta.
\end{align*}
(For instance, see Section 3.1 in Ledoux and
Talagrand~\cite{ledoux2013probability}).  Thus, the instance-dependent
term in Theorem~\ref{new--thm:bellman-err-bound-general-nonlinear}
matches the behavior of the limiting Gaussian random variable $W$, up
to constant factors and high-order terms.

Of course, it is natural to wonder whether this instance-dependent
term is actually optimal for stochastic approximation, or more
generally, for any procedure used to estimate the fixed point based on
stochastic observations.  As we discuss in
Section~\ref{sec:Estimation-error}, for a finite-dimensional space
 $\vecspace$, this term matches the fundamental lower bound provided by
local asymptotic minimax theory, so that---at least in general---it
cannot be improved.


\subsubsection{Stepsizes and higher-order term}

The range of permissible stepsizes and the higher-order term involve
certain Dudley entropy integrals, which we now define. Let $\ball$
denote the unit ball in the space $(\vecspace, \norm{\cdot})$, and let
$\ball^*$ denote the dual norm ball. Recall that $\dualBall$ is a
skeleton set such that $\ball^* = \mathrm{conv} (\dualBall)$.

Given a metric $\userho$ on the dual space, we define (for any $q \geq
1$) the Dudley entropy integral
\begin{align*}
 \dudley_q(\dualBall, \myrhotil) \mydefn \int_0^{\infty} \big[\log
   N(s; \dualBall, \myrhotil) \big]^{1/q} \; ds,
\end{align*}
where $N(s; \dualBall, \userho)$ denotes the cardinality of a minimal
$s$-covering of the skeleton set $\dualBall \subseteq \DualSpace$
under $\userho$.  Of particular interest are the cases $q = 2$ and $q
= 1$, which arise in the cases of sub-Gaussian and sub-exponential
tails, respectively.

In the simplest case (when $\ball^*$ is totally compact, such as the
finite-dimensional case), we can let $\userho$ be the dual norm
$\dualnorm{\cdot}$.  However, to handle the general
infinite-dimensional case and also sharpen our results, we make use of
the following pseudo-metric on the skeleton set
\begin{align}
\label{eq:defn-pseudo-metric-dual}  
    \myrhotil(x, y) \mydefn \sup_{e \in \numobs \obssubset \cap \ball}
    \inprod{x - y}{e} \quad \text{defined for all pairs} \quad x, y
    \in \dualBall,
\end{align}
where $\obssubset$ is the range of the operator.  Note that the
additional restriction $e \in \numobs \obssubset$ makes $\myrhotil$ a
weaker pseudo-metric than the dual norm, and in particular, we have
$\myrhotil(x, y) \leq \dualnorm{x - y}$ for any $x, y
\in \vecspace^*$.  This weakening is especially important in the
infinite-dimensional case, where the skeleton set $\dualBall$ is not
compact under the original norm $\dualnorm{\cdot}$.

With this notation, for a given tolerance probability $\delta \in
(0,1)$, the range of permissible stepsizes is given by
\begin{subequations}
\begin{align}
  \label{eqn:stepsize}   
  \stepsize \in \Big(0, \tfrac{(1 - \contraction)^2}{c \Lip^2
    \DudGauss^2 (\dualBall, \myrhotil)
    \log\big(\tfrac{\numobs}{\delta} \big)} \Big],
\end{align}
and the higher-order term is given by
\begin{align}
\label{EqnDefnFullHot}
\FullHot & \defn c \tfrac{\boundstar}{(1 - \contraction)} \left[
  \tfrac{1}{\numobs} + \tfrac{\stepsize \Lip}{\sqrt{\numobs}} \DudGauss(\dualBall,
  \myrhotil) \log( \tfrac{\numobs}{\delta}) \right] \cdot \Big \{
\DudExp(\dualBall, \myrhotil) + \log(\tfrac{1}{\delta}) \Big \}.
\end{align}
\end{subequations}

\subsubsection{Stepsize choice and restarting}\label{subsubsec:stepsize-choice-and-restarting}

Note that Theorem~\ref{new--thm:bellman-err-bound-general-nonlinear}
holds for a range of stepsizes, and the stepsize plays a role in the
burn-in length $\burnin = \tfrac{c}{(1 - \contraction)^2 \stepsize}
\log\big( \tfrac{\numobs}{\delta}\big)$.  In conjunction, these
requirements induce the following lower bound on the sample size
\begin{subequations}
\begin{align}
 \numobs \geq \tfrac{c}{(1 - \contraction)^4} \Lip^2
 \DudGauss^2(\dualBall, \myrhotil) \log^2(\tfrac{\numobs}{\delta}).
\end{align}
Given a sample size $\numobs$ satisfying this requirement, suppose
that we choose the stepsize
\begin{align}
\label{eq:main-prop-bellman-err-bound-stepsize-choice}  
 \stepsize = \Big \{ \Lip \DudGauss(\dualBall, \myrhotil) \log
 \big(\tfrac{\numobs}{\delta} \big) \sqrt{\numobs} \Big \}^{-1}.
\end{align}
With this choice, when evaluated at iteration $t = \numobs$, it can be
shown that the
bound~\eqref{new--eq:main-prop-bellman-err-bound-general-nonlinear}
simplifies to
\begin{multline}
\label{new--eq:main-prop-bellman-err-bound-general-nonlinear-optimal}
\norm{\hpop(\theta_\numobs) - \theta_\numobs} \leq \frac{c}{\sqrt{\numobs}}
\Big \{\Exs[\norm{\GaussNoise}] + \sqrt{\maxvar(W)
  \log(\tfrac{1}{\delta})} \Big \} + \frac{c \burnin}{\numobs} \cdot
\norm{\thetainit - \hpop(\thetainit)} \\
+ c \frac{\boundstar}{(1 - \contraction) \numobs} \Big \{
\DudExp(\dualBall, \myrhotil) + \log \big(\tfrac{1}{\delta} \big) \Big
\}.
\end{multline}
\end{subequations}
Thus, we see that the instance-dependent term (with its
$1/\sqrt{\numobs}$ decay) dominates the other two terms, which decay
at the faster $1/\numobs$ rate.

The final aspect that can be refined is the dependence of the
bound~\eqref{new--eq:main-prop-bellman-err-bound-general-nonlinear-optimal}
on the initial error.  As stated, this dependence is sub-optimal, but
can be refined via a simple restarting procedure, leading to an improved bound
\begin{align}
    \norm{\hpop(\theta_\numobs) - \theta_\numobs} \leq \frac{c}{\sqrt{\numobs}}
\Big \{\Exs[\norm{\GaussNoise}] + \sqrt{\maxvar(W)
  \log(\tfrac{1}{\delta})} \Big \} 
+ c \frac{\boundstar}{(1 - \contraction) \numobs} \Big \{
\DudExp(\dualBall, \myrhotil) + \log \big(\tfrac{1}{\delta} \big) \Big
\},
\end{align}
as long as the initial operator defect $\norm{\thetainit - \hpop (\thetainit)}$ is controlled by a finite-degree polynomial of the sample size $\numobs$.
See~\Cref{AppRestart} for the details of this procedure. In Corollary~\ref{cor:anorm-suboptimality-bound} and Theorem~\ref{new--thm:main-nonlinear-general} to follow, we assume that such re-starting scheme has been applied, so that the contribution from initial gap $\norm{\thetainit - \hpop (\thetainit)}$ is negligible.

\subsubsection{Semi-norm bounds on the operator defect}
\label{sec:operator-defect-anorm-bound}

There are various practical settings in which it is of interest to
obtain a bound in some semi-norm $\anorm{\cdot}$, as opposed to the
original Banach space norm $\norm{\cdot}$. As a simple example, in the
Euclidean setting, i.e. $\vecspace = \real^d$, one might have an
operator that is contractive in the $\ell_2$-norm, but be interested
in deriving bounds in the $\ell_\infty$-norm.  As a second example, in
various applications, one only cares about the error in some fixed
direction $\dir \in \real^\dims$, so that the semi-norm
$\anorm{\theta} \mydefn \abs{\dir^\top \theta}$ is the relevant
quantity.

In this section, we state a family of bounds applicable to any
semi-norm $\anorm{\cdot}$ of the form
\begin{align}
\label{eqn:semi-norm-defn}
\anorm{\theta} \defn \sup_{v \in \MySet} \inprod{v}{\theta} \qquad
\mbox{where $\MySet \subset \vecspace^\star$ is a symmetric and convex
  subset.}
\end{align}
Note that a wide class of interesting semi-norms can be generated in
this way.

A crude bound can be obtained by relating the semi-norm to the Banach
space norm. In particular, when the \emph{norm domination factor} $\aUB
\defn \sup_{v \in \MySet} \dualnorm{v}$ is finite, then any $\theta
\in \vecspace$ satisfies the upper bound
\begin{align}
\label{eqn:norm-domination}
\anorm{\theta} \leq (\sup_{v \in \MySet}\dualnorm{v}) \cdot
\norm{\theta} \defn \aUB \cdot \norm{\theta},
\end{align}
and a direct application of
Theorem~\ref{new--thm:bellman-err-bound-general-nonlinear} yields
\begin{align*}
\anorm{\hpop(\theta_\numobs) - \theta_\numobs} \leq \aUB \cdot
\norm{\hpop(\theta_\numobs) - \theta_\numobs} \lesssim
\tfrac{\aUB}{\sqrt{\numobs}} \cdot \Big \{ \Exs[\norm{W}] +
\sqrt{\maxvar(W) \log(\tfrac{1}{\delta})} \Big \}.
\end{align*}
This bound is potentially weak for two reasons: (a) the leading term
depends directly on $\aUB$, which can be large and possibly dependent
on the ambient dimension of the problem; and (b) it depends in a
global way on the Gaussian random element $W$, via the skeleton set
$\dualBall$ as opposed to $\MySet$, which can be much smaller.

It is natural to expect that one could prove bounds with a leading
term specified in terms of $\Exs \anorm{\GaussNoise}$ along with the
refined variance functional $\fullmaxvar{\MySet} \defn \sup_{u \in
  \MySet} \Exs[ \inprod{u}{\GaussNoise}^2]$.  This refinement is the
content of the following:
\begin{corollary}
\label{cor:anorm-suboptimality-bound}
Under the conditions of
Theorem~\ref{new--thm:bellman-err-bound-general-nonlinear}, the
iterate $\theta_\numobs$ satisfies the bound
\begin{align}
\label{eqn:anorm-bellman-bound}
\anorm{\hpop(\theta_\numobs) - \theta_\numobs} & \leq
\frac{c}{\sqrt{\numobs}} \left \{ \Exs[\anorm{W}] +
\sqrt{\fullmaxvar{\MySet}(\GaussNoise) \log(\tfrac{1}{\delta})} \right
\} + \NewHot
\end{align}
with probability at least $1 - \delta$, where
\begin{multline*}
\NewHot \defn \tfrac{\aUB }{ (1 - \contraction)} \Big \{ \LipCon
\constMGtype \log (\tfrac{\numobs}{\delta})
\sqrt{\tfrac{\stepsize}{\numobs}} + \tfrac{1}{\numobs \sqrt{
    \stepsize}} \Big \} \; \Big \{ \Exs[\norm{W}] + \sqrt{\maxvar(W)
  \log(\tfrac{\numobs}{\delta})} \Big \} \\
+ \tfrac{c \aUB \Lip \boundstar}{1 - \contraction} \Big \{
\tfrac{\sqrt{\stepsize}}{\numobs} + \tfrac{\stepsize}{\sqrt{\numobs}}
\Big \} \; \DudGauss(\dualBall, \myrhotil) \DudExp(\dualBall,
\myrhotil) \log^2(\tfrac{\numobs}{\delta}).
\end{multline*}
\end{corollary}
\noindent See Section~\ref{sec:proof-of-cor-anorm-suboptimality-bound}
for the proof.

\medskip 
With the stepsize
choice~\eqref{eq:main-prop-bellman-err-bound-stepsize-choice}, the
higher-order term scales as
\begin{align*}
\NewHot & = \Big \{ \Exs[\norm{\GaussNoise}] + \sqrt{\maxvar
  \log(\tfrac{1}{\delta})} \Big \} \cdot \myorderlog \left(
\tfrac{\aUB}{(1 - \contraction)\numobs^{3/4}} \right) + \myorderlog
\left(\tfrac{\aUB}{(1 - \contraction)\numobs}\right),
\end{align*}
where $\myorderlog$ subsumes various constants and logarithmic
factors.  Any dependence on global features of the Banach space
appears only in this higher-order term, which goes to zero at a rate
faster than $1/\sqrt{\numobs}$.


\subsection{Upper bounds on the estimation error}
\label{sec:Estimation-error}
Thus far, our analysis has focused on bounding the operator defect
$\theta_\numobs - \hpop(\theta_\numobs)$ in various (semi)-norms.  In
this section, we turn to problem of deriving upper bounds on the
estimation error $\norm{\theta_\numobs - \thetastar}$, which is the
primary goal in various applications of the SA methodology.  Bounds on
the operator defect imply bounds on this quantity: indeed, some simple
calculation\footnote{ By the triangle inequality, we have
$\norm{\theta_\numobs - \thetastar} \leq \norm{\theta_\numobs -
  \hpop(\theta_\numobs)} + \norm{\hpop(\theta_\numobs) -
  \hpop(\thetastar)}$.  From the contractivity
assumption~\ref{new--assume-pop-contractive}, we have
$\norm{\hpop(\theta_\numobs) - \hpop(\thetastar)} \leq \contraction
\norm{\theta_\numobs - \thetastar}$, and rearranging yields the
claim.}  yields the bound
\begin{align}
  \label{new--eqn:sub-optimality-to-estimation-err}
  \norm{\theta_\numobs - \thetastar} \leq \tfrac{1}{1 - \contraction}
  \cdot \norm{\theta_\numobs - \hpop(\theta_\numobs)}.
\end{align}
Although this bound is useful---and sharp in a worst-case sense--- it
can certainly be improved in general.

In this section, we develop a result (to be stated as
Theorem~\ref{new--thm:main-nonlinear-general}) that gives a sharper
bound on the estimation error $\norm{\theta_\numobs - \thetastar}$
when it is possible to construct linear approximations of the operator
$\hpop$ in a neighborhood of $\thetastar$.  More precisely, we impose
the following local linearity condition.
\paragraph{Assumption: Local linearity}
\begin{enumerate}[label=(A4)] 
\item \label{new--assume:linearization}
For any $\delrad >0$, there exists a set $\linearizationSet_\delrad$
of bounded linear operators on $\vecspace$ such that
\begin{align}
\norm{\theta - \thetastar} \leq \sup_{A \in \linearizationSet_\delrad}
\norm{(I - A)^{-1} \big( \hpop (\theta) - \theta \big)} \qquad
\mbox{for all $\theta \in \ball(\thetastar, \delrad)$.}
\end{align}
\end{enumerate}
As before, let $\GaussNoise$ be a centered Gaussian random variable in
$\vecspace$ with the same covariance structure as $\noise_1(\thetastar
) \defn \Hstoch_1(\thetastar) - \hpop(\thetastar)$---that is
\begin{align*}
\Exs \left[ \inprod{\GaussNoise}{y} \cdot \inprod{\GaussNoise}{z}
  \right] = \Exs \left[ \inprod{\noise_1(\thetastar)}{y} \cdot
  \inprod{\noise_1 (\thetastar)}{z} \right]\qquad \mbox{for all $y, z
  \in \DualSpace$.}
\end{align*}

Our bounds in this section are stated in terms of the solution to a
fixed-point equation involving functionals of the Gaussian noise
$\GaussNoise$.  For any $\delrad > 0$, define
\begin{align}
\label{eqn:complexity-terms-estimation}  
\localgauss(\delrad) \mydefn \Exs \Big[ \sup_{ \substack{y \in \dualBall
      \\ A \in \mathcal{A}_\delrad}} \inprod{\GaussNoise}{(I - A)^{-1}
    y} \Big], \quad \mbox{and} \quad \newvar^2(\delrad) \mydefn
\sup_{ \substack{y \in \dualBall \\ A \in \mathcal{A}_\delrad}} \Exs
\left[ \inprod{y}{(I - A)^{-1} \GaussNoise}^2 \right].
\end{align}
Given a stepsize $\stepsize$ satisfying~\eqref{eqn:stepsize} and a
tolerance probability $\delta \in (0, \tfrac{1}{1 +
  \log(1/(1-\contraction))})$, we define the function
\begin{subequations}
\begin{align}
\label{eqn:higher-order-term}
  \highorder_\numobs(\stepsize, \delta) & \mydefn \tfrac{\log
    (\tfrac{\numobs}{\delta})}{(1 - \contraction)^2} \Big \{ \big[
    \constMGtype \Lip \sqrt{\tfrac{\stepsize}{\numobs}} +
    \tfrac{1}{\numobs \sqrt{\stepsize}} \big] \cdot \Exs[\norm{W}] + \big
            [ \tfrac{\constMGtype \Lip \stepsize}{\sqrt{\numobs}}+
              \tfrac{1}{\numobs} \big] \cdot \boundstar \big[
              \DudExp(\dualBall, \myrhotil) + \log
              (\tfrac{\numobs}{\delta}) \big] \Big \}.
\end{align}
This quantity serves as a higher-order term in our analysis. We
consider the following fixed-point equation in the variable $\delrad$:
\begin{align}
\label{eq:thm-nonlinear-general-defining-eq-rate}  
\delrad = \frac{\localgauss (2 \delrad)}{\sqrt{\numobs}} +
\newvar(2\delrad) \sqrt{\frac{\log(1/\delta)}{\numobs}} +
\highorder_\numobs(\stepsize,
\delta).
\end{align}
\end{subequations}
As discussed below equation~\eqref{eqn:naive-bounds-explain} to
follow, equation~\eqref{eq:thm-nonlinear-general-defining-eq-rate} has
a non-empty and bounded set of non-negative solutions; let
$\delradstar$ be the largest such solution.

\begin{theorem}
\label{new--thm:main-nonlinear-general}
Suppose that
Assumptions~\ref{new--assume-pop-contractive}--~\ref{new--assume:linearization}
are in force, and that for some $\delta \in (0,1)$, we run
Algorithm~\ref{AlgRootSGD} using a stepsize $\stepsize$ in the
interval~\eqref{eqn:stepsize} and burn-in period \mbox{$\burnin =
  \tfrac{c}{(1 - \contraction)^2 \stepsize} \log\big(
  \tfrac{\numobs}{\delta} \big)$.}  Then the final iterate
$\theta_\numobs$ satisfies the bound
\begin{align}
  \norm{\theta_\numobs - \thetastar} \leq c \cdot \delradstar \quad
  \mbox{with probability at least $1 - \delta$.}
  \end{align}
\end{theorem}
\noindent See~\Cref{subsec:proof-thm-main-general-nonlinear} for the
proof of this theorem.

\medskip

Note that our contractivity assumption implies that functions
$\localgauss$ and $\localvar$ defined in
equation~\eqref{eqn:complexity-terms-estimation} are uniformly
bounded---viz.
\begin{align}
\label{eqn:naive-bounds-explain}
\localgauss(\delrad) & =\Exs \left[ \sup_{y \in \dualBall, A \in
    \mathcal{A}_\delrad} \inprod{\GaussNoise}{(I - A)^{-1} y} \right]
\leq \frac{\Exs\left[\norm{\GaussNoise} \right]}{1 - \contraction},
\quad \mbox{and} \notag \\
\newvar^2(\delrad) &\mydefn \sup_{y \in \dualBall, A \in
  \mathcal{A}_\delrad} \Exs \left[ \inprod{y}{(I - A)^{-1}
    \GaussNoise}^2 \right] \leq \frac{1}{(1 - \contraction)^2}\sup_{y
  \in \dualBall} \Exs \big[ \inprod{y}{\GaussNoise}^2 \big] .
\end{align}
These inequalities~\eqref{eqn:naive-bounds-explain}, in conjunction
with Theorem~\ref{new--thm:bellman-err-bound-general-nonlinear},
guarantee that the fixed-point
equation~\eqref{eq:thm-nonlinear-general-defining-eq-rate} has a
non-empty and bounded set of solutions; consequently, the maximum
solution $\delradstar$ is well-defined.  Moreover, this calculation
also reveals that the bound from
Theorem~\ref{new--thm:main-nonlinear-general} is always superior to
the naive bound~\eqref{new--eqn:sub-optimality-to-estimation-err}.

Note that only the high-order term $\highorder_\numobs (\stepsize,
\delta)$ depends on the stepsize. By taking the optimal stepsize
$\stepsize_\numobs = \Big \{ \Lip \DudGauss(\dualBall, \myrhotil) \log
\big(\tfrac{\numobs}{\delta} \big) \sqrt{\numobs} \Big \}^{-1}$, this
term becomes
\begin{align}\label{eq:high-order-in-est-err-bound}
      \highorder_\numobs(\stepsize_\numobs, \delta) & \mydefn
      \frac{\log (\tfrac{\numobs}{\delta})}{(1 - \contraction)^2} \Big
      \{ \frac{\sqrt{\Lip \constMGtype}}{\numobs^{3/4}} \cdot
      \Exs[\norm{W}] + \frac{\boundstar}{\numobs} \big[
        \DudExp(\dualBall, \myrhotil) + \log (\tfrac{\numobs}{\delta})
        \big] \Big \},
\end{align}
which consists of two terms: an $\myorder (\numobs^{-3/4})$ term
depending on the expected norm $\Exs[\norm{W}]$ that captures the
second moment of the noise, and an $\myorder (\numobs^{-1})$ term
depending on the worst-case upper bound on the noise, as well as the
Dudley integral. Under our stepsize choice, the high-order terms not
only decay at a faster rate with sample size $\numobs$, but also
capture the underlying complexity of the norm $\norm{\cdot}$, instead
of the ambient dimension of the space $\vecspace$.


\subsubsection{Asymptotic optimality}

Theorem~\ref{new--thm:main-nonlinear-general} provides a
non-asymptotic bound involving the Gaussian process $(I - \Amat)^{-1}
\GaussNoise$ for some $A \in \mathcal{A}_\delrad$.  It is natural to
ask whether or not this bound is improvable.  In certain cases it is
straightforward to address this question using local asymptotic
minimax theory
(cf.~\cite{lecam1953some,hajek1972local,van2000asymptotic}).

Let us suppose that $\vecspace$ is finite-dimensional, and the
operator $\hpop$ differentiable in an open neighborhood of the point
$\thetastar$.  In this case, we can use known results to state a lower
bound involving the random variable $(I - \Amat_0)^{-1} W$, where
$\Amat_0 = \nabla \hpop(\thetastar)$.  In order to state this lower
bound precisely, we consider problems indexed by distributions
$\mathbb{Q}$ in a local neighborhood of the target distribution
$\mathbb{P}$.  For any $\mathbb{Q}$, our goal is to solve the fixed
point equation $\theta = \Exs_{\Hstoch \sim \mathbb{Q}} \left[
  \Hstoch(\theta) \right]$ using i.i.d. samples of the random operator
$\Hstoch$.  Under suitable tail assumptions on the distributions
$\mathbb{P}$ and $\mathbb{Q}$, for any estimator
$\widetilde{\theta}_\numobs$ that maps a sequence of observed
operators $\{ \Hstoch_t \}_{t = 1}^\numobs$ to the vector space
$\vecspace$, an adaptation of Theorem 1 from the
paper~\cite{duchi2016asymptotic} (with loss function corresponding to
the Banach norm) yields the lower bound
\begin{align}
\label{eq:local-asymptotic-minimax}  
\liminf_{\lecamRad \rightarrow \infty}~ \liminf_{\numobs \rightarrow
  \infty} \sup_{\mathbb{Q} \, \mid \, \kull{\mathbb{Q}}{\Prob} \leq
  \tfrac{\lecamRad}{\numobs}} \Exs \left[ \sqrt{\numobs}
  \norm{\big(\widetilde{\theta}_\numobs - \thetastar (\mathbb{Q})
    \big) } \right] \geq \Exs \left[ \norm{(I - \Amat_0)^{-1} W}
  \right],
\end{align}
Thus, when estimating $\thetastar$ in the Banach norm $\norm{\cdot}$,
the asymptotic lower bound is given by $\Exs \left[ \norm{(I -
    \Amat_0)^{-1} \GaussNoise} \right]$.

Let us compare this fundamental limit to the behavior of the \rootSA
estimator.  We take a sequence of stepsizes $\{ \stepsize_\numobs
\}_{\numobs \geq 1}$ such that $\stepsize_\numobs \rightarrow 0^+$ and
$\numobs \stepsize_\numobs \rightarrow \infty$.  With this choice,
applying Theorem~\ref{new--thm:main-nonlinear-general} yields that the
\rootSA estimator $\theta_\numobs$ satisfies the bound
\begin{align}
\label{eq:asymptotic-upper-bound}  
\limsup_{\numobs \rightarrow \infty}~ \Prob \Big[ \|\theta_\numobs -
  \thetastar\| \geq c \cdot \Exs \left[ \norm{(I - \Amat_0)^{-1}
      \GaussNoise} \right] \Big] & \leq \tfrac{1}{3},
\end{align}
for some universal constant $c > 0$, showing its behavior is controlled
by the same functional that appears in the LAM lower bound.

\subsubsection{Semi-norm bounds on the estimation error}
\label{sec:estimation-error-anorm-bound}

Recall the setup of Section~\ref{sec:operator-defect-anorm-bound}.
We now refine these results by providing an upper bound on
$\anorm{\theta_\numobs - \thetastar}$, where $\anorm{\cdot}$ is a
semi-norm of the form~\eqref{eqn:semi-norm-defn}, assumed to satisfy
the domination condition~\eqref{eqn:norm-domination}. Moreover, we
assume the following modification of the local linearity condition
holds.
\paragraph{Assumption: Local linearity in semi-norm}
\begin{enumerate}[label=(A4)$^\prime$]
\item \label{new--assume:linearization-anorm}
For any $\delrad>0$, there is a set $\linearizationSet_\delrad$ of bounded
linear operators on $\vecspace$ such that
\begin{align}
\anorm{\theta - \thetastar} \leq \sup_{A \in \mathcal{A}_\delrad}
\norm{(I - A)^{-1} \big( \hpop(\theta) - \theta \big)} \qquad
\mbox{for all $\theta \in \ball(\thetastar, \delrad)$.}
 \end{align}
\end{enumerate}
As a refinement of the
definition~\eqref{eqn:complexity-terms-estimation}, we introduce the
complexity terms
\begin{align}
\label{eqn:complexity-terms-estimation-anorm}  
 \localgausssemi(\delrad) \mydefn \Exs \Big[ \sup_{ \substack{y \in
       \MySet \\ A \in \mathcal{A}_\delrad}} \inprod{\GaussNoise}{(I -
     A)^{-1} y} \Big], \quad \mbox{and} \quad \localvarsemi^2(\delrad)
 \mydefn \sup_{ \substack{y \in \MySet \\ A \in \mathcal{A}_\delrad}}
 \Exs \Big[ \inprod{y}{(I - A)^{-1} \GaussNoise}^2 \Big].
\end{align}
Given a stepsize $\stepsize$ satisfying the bound~\eqref{eqn:stepsize}
and a tolerance probability $\delta \in (0,
\tfrac{1}{\log(1/(1-\contraction))})$, we define $\delradastar > 0$ to
be the largest solution to the fixed-point equation
\begin{align}
  \delrad = \tfrac{\localgausssemi (2 \delrad)}{\sqrt{\numobs}}  +  \localvarsemi (2\delrad) \sqrt{\tfrac{\log (1/\delta)}{\numobs}}
  + \aUB \cdot \highorder_\numobs (\stepsize, \delta),\label{eq:fixed-point-eq-for-anorm-bound}
\end{align}
where the higher-order term $\highorder_\numobs (\stepsize, \delta)$
was previously defined~\eqref{eqn:higher-order-term}.

\begin{corollary}
\label{Cor:main-nonlinear-general-anrom}
Under
Assumptions~\ref{new--assume-pop-contractive}--~\ref{new--assume-noise-bound}
and~\ref{new--assume:linearization-anorm}, the estimate
$\theta_\numobs$ from Algorithm~\ref{AlgRootSGD} satisfies
\begin{align}
  \anorm{\theta_\numobs - \thetastar} \leq c \cdot \delradastar \quad
  \mbox{with probability at least $1 - \delta$.}
\end{align}
\end{corollary}
\noindent See~\Cref{subsec:proof-thm-main-general-nonlinear} for the
proof of this corollary.

\subsection{Linear operators with multi-step contraction}
\label{new--sec:multi-step-contraction}

In the special case where $\hpop$ is a bounded linear operator in
$\vecspace$, the contraction
assumption~\ref{new--assume-pop-contractive} can be significantly
weakened. In particular, it suffices to require that a multi-step
composition of the operator be contractive.

\paragraph{Assumption: Multi-step contraction}
\begin{enumerate}[label=(A1)$^\prime$]
\item \label{new--assume-compostep-contractive} For some integer
  $\compostep \geq 1$, the affine operator $\hpop(\theta) = \Amat
  \theta + \bvec$ is $\compostep$-stage contractive, meaning that
\begin{align}
\matsnorm{A}{\vecspace} \leq 1 \quad \mbox{ and} \quad \matsnorm{A^{\compostep}}{\vecspace} \leq
\tfrac{1}{2}.
\end{align}
\end{enumerate}
Note that assumption~\ref{new--assume-compostep-contractive} implies
that the linear operator $(I-A)$ is invertible; in particular, we have
the operator norm bound
\begin{align}
\label{eq:opnorm-bound-compostep}  
  \matsnorm{(I - A)^{-1}}{\vecspace} \leq \sum_{k = 0}^{\infty}
  \sup_{v \in \ball} \norm{A^k v}
  = \sum_{k = 0}^{\infty} \sum_{j =
    0}^{\compostep - 1} \matsnorm{A^{\compostep k+j}}{\vecspace}
  \leq \sum_{k = 0}^{\infty} \sum_{j =
    0}^{\compostep - 1} \matsnorm{A^\compostep}{\vecspace}^k \cdot
  \matsnorm{A^j}{\vecspace} \leq 2 \compostep.
\end{align}

As before, let $\GaussNoise$ be a centered Gaussian random variable in
$\vecspace$ with the same covariance structure as $\noise(\thetastar )
\defn \Hstoch(\thetastar) - \hpop(\thetastar)$; that is, \mbox{$\Exs
  \left[ \inprod{\GaussNoise}{y} \cdot \inprod{\GaussNoise}{z} \right]
  = \Exs \left[ \inprod{\noise(\thetastar)}{y} \cdot \inprod{\noise
      (\thetastar)}{z} \right]$} \mbox{for all $y, z \in \DualSpace$.}

\paragraph{Tuning parameters:}
Given a desired failure probability $\delta \in (0, 1)$, and a total
sample size $\numobs$, we run Algorithm~\ref{AlgRootSGD} with the
following choices of parameters: 
\begin{subequations}
\label{eqn:Tuning-pars-multistep}
\begin{align}
  \label{eqn:stepsize-multistep}   
  \texttt{Stepsize choice:} & \qquad
   \stepsize \leq
     \frac{c}{\compostep
     \Lip^2 \constMGtype^2 \cdot \log^2 \tfrac{\numobs}{\delta}} \\
  \label{eqn:burnin-multistep}  
  \texttt{Burn-in time:} & \qquad \burnin = \tfrac{c\compostep}{
    \stepsize} \log(\tfrac{\numobs}{\delta}),
\end{align} 
where $c$ is an universal constant.
\end{subequations}

\begin{theorem}
\label{thm:bellman-err-bound-general-linear}
Under Assumptions~\ref{new--assume-compostep-contractive},
~\ref{new--assume-sample-operator} and~\ref{new--assume-noise-bound},
and given a sample size $\numobs \geq 2 \burnin $, consider
Algorithm~\ref{AlgRootSGD} run using tuning parameters from
equation~\eqref{eqn:stepsize-multistep}
and~\eqref{eqn:burnin-multistep}.  Then for any given $t \in [\burnin,
  \numobs]$, the iterate $\theta_t$ satisfies
\begin{align}
\label{eq:main-thm-bellman-err-bound-general-linear}  
\norm{\hpop (\theta_\numobs) - \theta_\numobs} & \leq \frac{c}{\sqrt{\numobs}} \Big \{
\Exs[\norm{\GaussNoise}] + \sqrt{ \maxvar \log(\tfrac{1}{\delta})}
\Big \} +
 \frac{c \burnin}{\numobs} \norm{\thetainit - \hpop
   (\thetainit)} + \ThreeHot
\end{align}
with probability $1 - \delta$, where
\begin{align*}
\ThreeHot  \defn {c \boundstar} \Big \{ \tfrac{1}{\numobs} +
\tfrac{\stepsize \Lip \constMGtype}{ \sqrt{\numobs}}
\log(\tfrac{\numobs}{\delta}) \Big \} \Big \{ \DudExp(\dualBall,
\myrhotil) + \log(\tfrac{1}{\delta}) \Big \},
\end{align*}
\end{theorem}
\noindent See~\Cref{subsec:proof-theorem-thm-linear-compostep} for the
proof of this theorem.

\medskip

Observe that for a linear operator $\hpop(\theta) = \Amat \theta +
\bvec$ that satisfies the contractivity condition
(cf.\ Assumption~\ref{new--assume-compostep-contractive}), the inverse
$(\Id - \Amat)^{-1}$ exists, and we have
\begin{align*}
  \theta - \thetastar = (\Id - \Amat)^{-1}(\hpop(\theta) - \thetastar).
\end{align*}
Consequently, given any semi-norm $\anorm{\cdot}$ of the
form~\eqref{eqn:semi-norm-defn} satisfying
condition~\eqref{eqn:norm-domination}, an argument similar to
Corollary~\ref{Cor:main-nonlinear-general-anrom} yields the following
guarantee.  In stating it, we assume that the restarting scheme
from~\Cref{AppRestart} has been applied to remove dependence on the
initial condition.
\begin{corollary}
\label{thm:est-err-linear}
Under the conditions of
Theorem~\ref{thm:bellman-err-bound-general-linear}, running the
\rootSA algorithm with the restarting scheme yields an iterate
$\theta_\numobs$ such that
\begin{align}
\label{eq:est-error-linear-main}  
\anorm{\theta_\numobs - \thetastar} \leq \frac{c}{\sqrt{\numobs}} \Big
\{ \Exs \big[ \anorm{(I - \Amat)^{-1} W} \big] + \sqrt{\sup_{u \in
    \MySet} \Exs \big[ \inprod{u}{(I - \Amat)^{-1} W}^2 \big]
  \log(\tfrac{1}{\delta})} \Big \} + \SuperHot{\diamond}
\end{align}
with probability at least $1 - \delta$.
\end{corollary}
\noindent See~\Cref{subsec:proof-thm-est-err-linear} for the proof of
this corollary, along with the definition of $\SuperHot{\diamond}$.

Since the problem itself is linear, the class $\mathcal{A}_s$ of
linear operators is singleton, and the estimation error upper bounds
can be expressed directly through $\Exs \big[ \anorm{(I - \Amat)^{-1}
    W} \big] $, without resorting to fixed-point equations. Compared
with the high-order terms defined by
equation~\eqref{eqn:higher-order-term} in the general case, the high
order terms in equation~\eqref{eq:est-error-linear-main} (the second
and third line of the equation) save a factor of $\frac{1}{1 -
  \contraction}$ in the contractive case, while generalizing to the
multi-step contraction case. Furthermore, similar to the discussion in
Section~\ref{subsubsec:stepsize-choice-and-restarting}, the step
$\stepsize$ can be tuned based on the sample size $\numobs$ and
knowledge about other problem parameters, so as to minimize the
high-order terms $\ThreeHot$ and $\SuperHot{\diamond}$. The resulting
error bounds contain high-order terms similar to
equations~\eqref{eqn:higher-order-term}
and~\eqref{eq:high-order-in-est-err-bound}, the factor $(1 -
\contraction)^{-2}$ replaced by the integer $\compostep$.


\section{Consequences for specific use cases}
\label{new--sec:Applications}

Thus far, we have stated a number of general results.  In this
section, we discuss the consequences of these results for three
classes of problems that fall within the framework of this paper.  In
the main text, we discuss in detail the problem of stochastic shortest
paths in~\Cref{new--sec:SSP} and average-reward policy evaluation
in~\Cref{SecAveReward}.  We defer discussion of methods for solving
two-player zero-sum Markov games to~\Cref{sec:Markov-games}.

\subsection{Computing stochastic shortest paths}
\label{new--sec:SSP}

We begin with the problem of computing stochastic shortest
paths~\cite{yu2013q,bertsekas1991analysis}, or SSPs for short.  It
provides an illustration of the general theory using a Banach space
defined by a certain weighted $\ell_\infty$-norm.  On one hand, SSPs
can be formulated in terms of Markov decision process (MDP) with a
finite state space $\states$ and action space $\actions$.  Thus,
although they might appear to be a special case of an MDP, in fact,
they are sufficiently general to encompass both finite-horizon MDPs as
well as discounted MDPs.  Thus, the conclusions obtained in this
section apply to a fairly broad class of problems.

An MDP is defined by a collection of probability transition kernels
$\bracketMed{\TransMat_{\action}(\cdot \mid \state)}_{(\state,
  \action) \in \states \times \actions}$, where the transition kernel
$\TransMat_{\action}(\state' \mid \state)$ denotes the probability of
transition to the state $\state'$ when an action $\action$ is taken at
the current state $\state$. The MDP is equipped with a cost function
$\cost: \states \times \actions \mapsto \real$, and the value
$\cost(\state, \action)$ corresponds to cost incurred upon performing
the action $\action$ in state $\state$.  To formulate a
\emph{stochastic shortest path} (SSP) problem, we assume that state
$1$ is absorbing and cost-free, meaning that
\begin{align}
\label{eqn:cost-free-absorbing-state}
\cost(\state = 1, \action) = 0 \quad \text{and} \quad
\Prob_{\action}(\state' \mid \state = 1) = \bm{1}_{\state' = 1} \quad
\mbox{for all actions $\action \in \actions$.}
\end{align}
A stationary policy $\policy$ is a mapping $\states \mapsto \actions$
such that $\policy(\state) \in \actions$ denotes the action to be
taken in the state $\state$.  We assume that the total
infinite-horizon cost incurred by any stationary policy $\pi$ is
finite---viz.  $\Exs_{ \state_0 = \state} \Big[ \sum_{k = 1}^\infty
  \abss{\cost(\state_k, \policy(\state_k))} \Big] < \infty$ \mbox{for
  all $\state \in \states$.}  Such stationary policy $\pi$ is called a
\emph{proper policy}, and our goal is to obtain proper policy
$\optPolicy$ that minimizes the total cost.

Associated with any proper policy $\pi$ is its $Q$-function
\begin{align*}
\theta^{\pi}(\state, \action)& \defn \Exs \Big[\sum_{k=0}^{\infty}
  \cost (\state_k, \action_k ) \mid \state_0= \state, \action_0 =
  \action \Big], \quad \text { where } \action_k = \policy(\state_k)
\quad \mbox{for all $k = 1, 2, \ldots$.}
\end{align*}
An optimal policy can be obtained from the optimal $Q$-function, given
by \mbox{$\theta^\star(\state, \action) \defn \inf_{\pi \in \Pi}\;
  \theta^{\pi}(\state, \action)$.}

\subsubsection{Bellman operator and contractivity}

Observe that for any policy $\policy$, the cost-free absorbing state
property~\eqref{eqn:cost-free-absorbing-state} ensures that
$\theta^\policy(1, \action) = 0$, and as a result $\theta^\star(1,
\action) = 0$ for all actions $\action \in \actions$. In terms of the
shorthand $ \nonAbsStates \defn \states \setminus \{ 1\}$, classical
theory~\cite{bertsekas1991analysis,yu2013q} guarantees that the
optimal $Q$-function restricted to the set $\nonAbsStates \times
\actions$ is the unique fixed point of the Bellman operator
\begin{align}
\label{eqn:SSP-Bellman}
\SSPOp(\theta)(\state, \action) = \cost(\state, \action) +
\sum_{\state' \in \nonAbsStates} \TransMat_{\action}(\state' \mid
\state) \min_{\action' \in \actions} \theta(\state', \action') \qquad
(\state, \action) \in \states_{-1} \times \actions.
\end{align}

For SSP problems with finite state and action spaces, any $Q$-function
can be viewed an element of $\real^{|\states_{-1} \times \actions|}$,
in which case the Bellman operator $\SSPOp$ can be viewed as acting on
$\real^\SSPdim$ where $\SSPdim \defn |\nonAbsStates \times \actions|$.
For a vector $\NormWeight \defn \{w_1, \ldots, w_{\SSPdim}\}$ of
strictly positive weights, we define a weighted $\ell_\infty$-norm on
$\real^\SSPdim$ via $\wNorm{\theta} \defn \max \limits_{i = 1, \ldots,
  \SSPdim} \tfrac{|\theta_i|}{w_i}$.  From known results on SSP
problems~\cite{bertsekas1991analysis,tseng1990solving}, one can use a
hitting time analysis to define a weight vector $\NormWeight$ such
that, for any $\theta_1, \theta_2 \in \real^{\SSPdim}$, we have
\begin{align}
  \label{eqn:SSP-contraction}
  \wNorm{\SSPOp(\theta_1) - \SSPOp(\theta_2)} \leq \Big(1 -
  \frac{1}{\weightmax} \Big) \cdot \wNorm{\theta_1 - \theta_2}
\end{align}
where $\weightmax = \max \limits_{i= 1, \ldots, \SSPdim} w_i \geq 1$.
Thus, the Bellman operator $\SSPOp$ is $\big(1 - \tfrac{1}{\weightmax}
\big)$-contractive in the weighted $\ell_\infty$-norm, so that our
general theory can be applied with this choice of Banach space.


\subsubsection{Generative observation model}

We analyze the \rootSA~algorithm under a stochastic oracle known as
the \emph{generative observation model} for the SSP problem.  For any
state-action pair $(\state, \action)$, the generative model allows us
to draw next-state and cost samples from the MDP $(\rewardQ,
\TransMat)$.  More precisely, we have access to a collection of
$\numobs$ i.i.d.\ samples of the form $\left\{\left(\TranSample_{k},
\costRnd_{k}\right)\right\}_{k=1}^{\numobs}$, where both
$\TranSample_{k}$ and $\costRnd_{k}$ are random matrices in
$\real^{|\mathcal{X}_{-1}| \cdot|\mathcal{U}|}$. For each state-action
pair $(\state, \action)$, the entry $\TranSample_{k}(\state, \action)$
is drawn according to the transition kernel $\TransMat_{\action}(\cdot
\mid \state)$, whereas the entry $\costRnd_{k}(\state, \action)$ is a random variable with mean $\cost(\state, \action)$; this
corresponds to a noisy observation of the cost function.  We assume
that the random cost $\costRnd_{k}(\state, \action)$ is upper bounded
by $\costMax$ in absolute value.  Here the cost samples
$\left\{\costRnd_{k}(\state, \action)\right\}_{(\state, \action) \in
  \states \times \actions}$ are independent across all state-action
pairs, and the cost samples $\left\{\costRnd_{k}\right\}$ are
independent of the transition samples
$\left\{\TranSample_{k}\right\}$.

\paragraph{The empirical Bellman operator:}
Given a sample $(\TranSample, \costRnd)$ from our observation model,
we define the single-sample empirical Bellman operator
$\SSPstoch(\cdot)$ on the space of $Q$-functions, whose action on a
$Q$-function $\theta$ is given by
\begin{align}
\label{eqn:SSP-Bellman-Noisy}
\SSPstoch(\theta)(\state, \action) & \mydefn \costRnd(\state, \action)
+ \sum_{\state' \in \nonAbsStates} \TranSample_{\action}\left(\state'
\mid \state \right) \min_{\action' \in \actions} \theta\left(\state',
\action'\right).
\end{align}
Here we have introduced $\TranSample_{\action}\left(\state' \mid
\state\right) \mydefn \mathbf{1}_{\TranSample(\state,
  \action)=\state'}$. We are ready to state our guarantees for the
stochastic shortest path problem.


\subsubsection{Guarantees for stochastic shortest path}

It is easy to see that the operators $\hpop(\cdot)$ and
$\Hstoch(\cdot)$, defined respectively in
equations~\eqref{eqn:SSP-Bellman} and~\eqref{eqn:SSP-Bellman-Noisy},
satisfy
Assumptions~\ref{new--assume-pop-contractive}-~\ref{new--assume-noise-bound}
with the weighted $\ell_\infty$-norm $\wNorm{\cdot}$. In order to
obtain an optimal policy from an estimate $\theta_\numobs$ of the
optimal $Q$ function, it is natural to obtain performance bounds in
the $\infNorm{\cdot}$ norm, and we do so by invoking Corollaries~\ref{cor:anorm-suboptimality-bound}
and~\ref{Cor:main-nonlinear-general-anrom} with $\anorm{\cdot} = \|
\cdot \|_\infty$.

Accordingly, consider a Gaussian random vector $\GaussNoise$ with
\mbox{$\GaussNoise \sim \mathcal{N} \Big( 0, \mathrm{cov}(\Hstoch
  (\thetastar) - \thetastar) \Big)$}, and define
\begin{align}
\label{eqn:key-complexity-SSP}
\noisegauss = \Exs[\|W\|_\infty], \qquad \noisevar^2 \mydefn
\sup_{\state \in \nonAbsStates, \action \in \actions} \Exs
    [\GaussNoise_{\state, \action}^2], \quad \mbox{and} \quad
    \boundstar \mydefn \frac{\costMax}{\weightmin} +
    \wNorm{\thetastar}.
\end{align}
For a given failure probability $\delta \in (0,1)$, our result applies
to the algorithm with parameters
\begin{subequations}
\label{eqn:SSP-tuning}
\begin{align}
\label{EqnOmicron}  
\stepsize = c_1 \Big\{\sqrt{ \numobs \log \big( |\states| \cdot |\actions|
  \big)} \cdot \log (\numobs / \delta) \Big\}^{-1}, \quad \mbox{and}
\quad \burnin = \tfrac{c_2 \weightmax^2}{ \stepsize} \log
(\tfrac{\numobs}{\delta}),
\end{align}
We also choose the initialization $\theta_0$ and the number of restarts $\numRestarts$ such that
\begin{align}
\label{eqn:numrestarts-and-init-SSP}
 \log \left( \tfrac{\norm{\theta_0 - \hpop (\theta_0)}
   \sqrt{\numobs}}{\noisegauss} \right) \leq c_0 \log \numobs\qquad
 \text{and} \qquad \numRestarts \geq 2 c_0 \log \numobs,
\end{align}
\end{subequations}
where $c_0, c_1, c_2$ are appropriate universal constants.  We obtain
the following guarantee:
\begin{corollary}
    \label{cors:SSP-cors}
Given a sample size $\numobs$ such that $\tfrac{\numobs}{\log \numobs}
\geq {c' \log(\abss{\states} \cdot \abss{\actions})} \cdot
\weightmax^4 \log(1 / \delta)$, running Algorithm~\ref{AlgRootSGD}
with the tuning parameter choices~\eqref{eqn:SSP-tuning} yields an
estimate $\theta_\numobs$ such that
\begin{align*}
\vecnorm{\hpop(\theta_\numobs) - \theta_\numobs}{\infty} \leq
\frac{c}{\sqrt{\numobs}} \cdot \Big \{ \noisegauss + \noisevar \sqrt{
  \log(\tfrac{1}{\delta})} \Big \} + c \boundstar \weightmax^2 \tfrac{\log (|\states| \cdot
  |\actions|)}{\numobs} \log^2(\tfrac{\numobs}{\delta}),
\end{align*}
with probability at least $1 - \delta$.
\end{corollary}
Note that when we invoke Corollary~\ref{cor:anorm-suboptimality-bound}
to obtain this corollary, the second term is absorbed into the
leading-order term under the sample size lower bound
$\tfrac{\numobs}{\log \numobs} \geq {c' \log(\abss{\states} \cdot
  \abss{\actions})} \cdot \weightmax^4 \log(1 / \delta)$. In
particular, the semi-norm domination factor is $\aUB = \weightmax$ in
this case, and we have the following inequalities:
\begin{align*}
   \aUB \cdot \Exs \left[ \wNorm{W} \right] &\leq
   \frac{\weightmax}{\weightmin} \Exs \left[ \vecnorm{W}{\infty}
     \right] \leq \weightmax \noisegauss, \quad \mbox{and}\\ \aUB
   \cdot \sup_{ \vecnorm{y}{1/w} \leq 1} \sqrt{ \Exs \big[
       \inprod{y}{W}^2 \big] } & \leq \frac{\weightmax}{\weightmin}
   \sup_{\vecnorm{y}{\infty}\leq 1} \sqrt{\Exs \big[ \inprod{y}{W}^2}
     \big]\leq \weightmax \noisevar,
\end{align*}
which makes the second term of
equation~\eqref{eqn:anorm-bellman-bound} dominated by the first term.

Next, in order to obtain an upper bound on the estimation error
$\|\theta_\numobs - \thetastar\|_\infty$ we need a few more definitions.  For a given $Q$-function $\Q$, we say $\pi$ is a greedy policy of $\theta$ if and only if
\begin{align*}
  \pi (\state) = \arg\min_{\action} \Q(\state,
\action) \qquad \text{for all}
\quad \state
\in \states_{-1}, 
\end{align*}
and denote $\Pi^\theta$ as the set of all  greedy policies of $\theta$.
Note that the greedy policies of a given $Q$-function may not be unique.  Using
this greedy
policy, we can define the right-linear operator
\begin{align*}
  \TransMat^{\pi_{\Q}} \Q(\state, \action) = \sum_{\state'}
  \TransMat_{\action}(\state' \mid \state) \Q(\state',
  \pi_{\Q}(\state')). 
\end{align*}
We also define a set
$\mathcal{A}_{\delrad}$ of linear operators as
\begin{align}
\label{eqn:A_delta-SSP-2}
\mathcal{A}_{\delrad} = \{ \TransMat^{\pi_\Q} \;
\mid \; \pi_\Q \;\; \text{is a greedy policy of } \Q \text{ with } \Q
\in \ball(\Qstar, \delrad)\}.
\end{align}
Let $\ball(\thetastar, s) \defn \{ \theta \mid \|\theta -
\thetastar\|_\infty \leq s \}$ denote the $\ell_{\infty}$-ball of
radius $s$ around $\thetastar$.  We use $\pi_\star$ to denote the
greedy policy associated with the optimal $Q$-function $\Qstar$.  In
Appendix~\ref{sec:localLinSSPProof}, we show that the local linearity
assumption~\eqref{new--assume:linearization-anorm} is satisfied for
the Bellman operator~\eqref{eqn:SSP-Bellman} with the set of operators
$\mathcal{A}_\delrad$ from equation~\eqref{eqn:A_delta-SSP-2}, and
with $\anorm{\cdot} = \|\cdot\|_\infty$.

\vspace{10pt}

Given a tolerance probability $\delta \in (0, \tfrac{1}{\log(1/(1-\contraction)}))$, let $\delrad^*_\numobs$ denotes 
the largest positive solution to the fixed-point
equation
\begin{multline}
  \delrad_\numobs = \frac{1}{\sqrt{\numobs}} \Big \{ \Exs \Big[
    \sup_{\stackrel{\theta \in \ball_\infty(\thetastar, \delrad_\numobs)}{\pi \in \Pi^\theta}} \big
    \|\big( \Id - \TransMat^{\pi} \big)^{-1} \GaussNoise \big
    \|_\infty \Big] + \sup_{\stackrel{\theta \in
      \ball_\infty(\thetastar, \delrad_\numobs), \pi \in \Pi^\theta}{(\state, \action) \in \states_{-1} \times \actions } }
  \mybiggersqrt{\Exs \big[ \delta_{\state,\action}^\top (\Id -
      \TransMat^{\pi})^{-1} \GaussNoise \big] \log(1 / \delta)}
  \Big \} \\
  + {\weightmax^2 \log( \tfrac{\numobs}{\delta})} \left\{
  \tfrac{\weightmax}{\weightmin} \big( \tfrac{\log(\abss{\states_{-1}}
    \cdot \abss{\actions})}{\numobs} \big)^{3/4} \; \; \Exs
        [\|W\|_\infty] + \boundstar \weightmax \tfrac{\log
          (|\states_{-1}| \cdot |\actions|)}{\numobs}
        \right\}.\label{eq:fixed-point-eq-ssp}
\end{multline}
Here we have defined the indicator function $\delta_{\state, \action}
= \mathbf{1}_{(\state', \action') = (\state, \action)}$.
We obtain the following corollary:
\begin{corollary}
\label{cors:SSP-est-error}
Under the setup of Corollary~\ref{cors:SSP-cors}, the estimate
$\theta_\numobs$ satisfies the bound
\begin{align}
\|\theta_\numobs - \thetastar \|_\infty & \leq c \cdot
\delrad^*_\numobs \quad \mbox{with probability at least $1 - \delta$.}
\end{align}
\end{corollary}
A few remarks are in order. First, the bound depends on the size of
state-action space only poly-logarithmically, and depends on the
quantity $\weightmax$ through two sources: the contraction parameter
and the norm domination factor between $\infNorm{\cdot}$ and
$\wNorm{\cdot}$. Second, let $\Pi^*$ be the set of all optimal
policies for the SSP problem, for sample size $\numobs$ large
enough,\footnote{The sample size requirement may depend on the gap
between the value of optimal and sub-optimal actions, as in the prior
work~\cite{khamaru2021instance}.} the ball $\ball_\infty (\thetastar,
s_n)$ will eventually shrink to the singleton $\thetastar$, and the
supremum in the fixed-point equation~\eqref{eq:fixed-point-eq-ssp} is
taken over $\pi \in \Pi^*$. Therefore, using $\PLAINHOT$ to denote
higher-order terms, the solution $\delrad_n$ takes the form
\begin{align*}
 \delrad_n & = \frac{1}{\sqrt{\numobs}} \Big\{ \Exs \Big[ \sup_{\pi
     \in \Pi^*} \big \|\big( \Id - \TransMat^{\pi} \big)^{-1}
   \GaussNoise \big \|_\infty \Big] + \sup_{\stackrel{\pi \in
     \Pi^*}{(\state, \action) \in \states_{-1} \times \actions } }
 \mybiggersqrt{\Exs \big[ \delta_{\state,\action}^\top (\Id -
     \TransMat^{\pi})^{-1} \GaussNoise \big] \log( \tfrac{1}{\delta})}
 \Big\} + \PLAINHOT \\
 & \leq \frac{1}{\sqrt{\numobs}} \max_{\stackrel{\pi \in
     \Pi^*}{(\state, \action) \in \states_{-1} \times \actions } }
 \sqrt{\Exs \left[ \big( \delta_{\state, \action}^\top (\Id -
     \TransMat^\pi)^{-1} (\Hstoch (\thetastar) - \thetastar) \big)^2
     \right]} \cdot \sqrtlog{\tfrac{|\states|\cdot |\actions| \cdot
     |\Pi^*|}{\delta}} + \PLAINHOT.
\end{align*}
Up to a factor of $ \sqrtlog{\tfrac{|\states|\cdot |\actions| \cdot
    |\Pi^*|}{\delta}}$, this matches the two-point lower bound in the
paper~\cite{khamaru2021instance} (in the discounted MDP case). When
specializing to the cases where the optimal policy is unique, or
satisfies the Lipschitz-type assumptions in the
paper~\cite{khamaru2021instance}, the upper bound above also recovers
the leading-order term in that paper.  We conjecture that the
leading-order term of the solution $s_n$ to the fixed-point equation
is actually optimal for large $\numobs$. It would be interesting to
verify this conjecture, and establish some kind of optimality over
suitably defined problem classes.

When specialized to the $\contraction$-discounted MDPs, the sample
size requirement in Corollary~\ref{cors:SSP-est-error} scales as
$\order{(1 - \contraction)^{-4}}$.  This requirement is worse than
corresponding requirements in the paper~\cite{khamaru2021instance}, at
least in certain regimes. Intuitively, this is the price we pay when
moving to the general case where only the contraction of the
population-level operator is assumed, instead of the sample-level
contraction.


\subsection{Average cost policy evaluation}
\label{SecAveReward}

As a second illustration, we turn to a problem where the operator is
\emph{not} contractive, but does satisfy a form of multi-step
contractivity needed to apply the theory
from~\Cref{new--sec:multi-step-contraction}.  This example also
involves an error measure that is only a semi-norm in the original
space, but can converted to a norm in a Banach space by taking a
suitable quotient.

More specifically, consider an undiscounted Markov reward process
(MRP) with state space $\states$, probability transition kernel
$\TransMat \in \real^{\states \times \states}$ and cost function
$\cost: \states \rightarrow \real$.  When the Markov chain is
irreducible and ergodic, there is a unique stationary distribution
$\stationary$.  Letting $\avgcoststar \mydefn \Exs_{\state \sim
  \stationary} [\cost (\state)]$ denote the average cost under this
stationary distribution, our goal is to estimate the value function
\begin{align*}
\thetastar(\state) \mydefn \sum_{t = 0}^{\infty} \TransMat^t \big \{
\cost (\state) - \avgcoststar \big \}.
\end{align*}
It is known that the value function $\thetastar$ and average cost
$\avgcoststar$ jointly satisfy the Bellman equation
\begin{align}
\label{eq:bellman-average-reward}  
\avgcoststar + \thetastar(\state) - \TransMat \thetastar(\state) -
\cost(\state) = 0 \quad \mbox{for all $\state \in \states$.}
\end{align}
See the sources~\cite{derman1966denumerable,tsitsiklis1999average} for
more background.

In practical application of policy evaluations problems, of primary
interest are the relative differences between the value function at
different state-actions pairs.  Thus, the primary goal is to estimate
the function $\thetastar$, with the average cost $\avgcoststar$ being
a nuisance parameter. As shown in the sequel (see~\Cref{SecAveSemi}),
by considering the span semi-norm in an appropriate vector space
$\vecspace$, it is possible to estimate $\thetastar$ without
estimating $\avgcoststar$.

\paragraph{Observation models and relevant operators:}

As before, we consider a generative observation model, where we
observe a collection of $\numobs$ i.i.d. samples of the form
$\left\{\left(\TranSample_{k},
\costRnd_{k}\right)\right\}_{k=1}^{\numobs}$, where both
$\TranSample_{k} \in \real^{|\states| \times |\states|}$ and
$\costRnd_{k} \in \real^{|\states|}$. For each state $\state \in
\states$, the row $\state$ of the matrix $\TranSample_k$ is an
indicator vector $\bm{1}_{s'}$, where the state $s'$ is drawn
according to the transition kernel $\TransMat (\cdot \mid \state)$;
the entry $\costRnd_{k}(\state)$ is a random variable with mean
$\cost(\state)$ and uniformly bounded by $\sigma_{r}$, corresponding
to a noisy observation of the reward function.
 
The population and empirical Bellman operators for the average-cost policy evaluation can be written as follows:
\begin{align*}
    \hpop(\theta) \mydefn \TransMat \theta + \cost,\quad \mbox{and} \quad \Hstoch_k (\theta) \mydefn \TranSample_k \theta + \costRnd_k.
\end{align*}
It can be seen that both $\hpop$ and $\Hstoch_k$ are linear operators,
satisfying $\Exs[\Hstoch_k] = \hpop$.

In the rest of this subsection, we define a semi-norm and discuss the
multi-step contraction properties of the operator $\hpop$, and then
present the main consequences of
Theorem~\ref{thm:bellman-err-bound-general-linear} and
Corollary~\ref{thm:est-err-linear} for such models.


\subsubsection{The semi-norm and multi-step contraction}
\label{SecAveSemi}
Consider the Banach space $\vecspace$ given by
\begin{align}
\label{eqn:span-quotient-space}
\vecspace = \real^{|\states|} / \big\{\theta + \alpha \bm{1} \; \mid
\; \alpha \in \real \big\},
\end{align}
where each element of $\vecspace$ is an equivalence class of the form
$\big\{ \theta + \alpha \bm{1} : \alpha \in \real \big\}$, equipped
with the span norm
\begin{align*}
\vecnorm{\theta}{\myspan} \mydefn \max_{\state \in \states}
\theta(s) - \min_{\state \in \states} \theta(s) \quad \mbox{for all
  $\theta \in \vecspace$.}
\end{align*}
Note that $\vecnorm{\cdot}{\myspan}$ is a semi-norm on
$\real^{\states}$, but a norm on the quotient space $\vecspace$.  For
reinforcement learning problems, this choice is natural, since we
often care only about the relative advantages of state-action pairs,
in which case the average cost $\avgcoststar$ is irrelevant.

Under the norm $\vecnorm{\cdot}{\myspan}$ on $\vecspace$, the operator
$\hpop$ is non-expansive, but not necessarily a contraction. However,
under suitable conditions, it can be shown to contractive in a
multi-step sense.  In order to do so, we impose the following mixing
time condition. \\

\noindent {\bf{Assumption: Mixing time}}
\begin{enumerate}[label=(MT)] 
\item \label{assume:mixing-time} There exists a positive integer
  $\mixingtime$ such that
\begin{align*}
        \totalvariation (\delta_x^\top \TransMat^{\mixingtime},
        \delta_y^\top \TransMat^{\mixingtime}) \leq \tfrac{1}{2}
        \qquad \mbox{for any $x, y \in \states$.}
\end{align*}
Here the vector $\delta_x \in \real^{\states}$ is the unit basis
vector with a single one in entry $x \in \states$.
\end{enumerate}

\vspace*{0.02in}

\noindent Under Assumption~\ref{assume:mixing-time}, for any $\theta
\in \vecspace$, we have
\begin{multline}
    \vecnorm{\TransMat^{2 \mixingtime} \theta}{\myspan} =
    \max_{\state \in \states} \left\{ \delta_{\state}^\top \TransMat^{2
      \mixingtime} \theta \right\} - \min_{\state \in \states} \left\{
    \delta_{\state}^\top \TransMat^{2 \mixingtime} \theta \right\} \overset{(i)}{\leq} 2
    \max_{x \in \states} \abss{\delta_{\state}^\top \TransMat^{2
        \mixingtime} \theta - \stationary^\top \TransMat^{2
        \mixingtime} \theta}\\ \leq 2 
    \totalvariation (\delta_x \TransMat^{2 \mixingtime}, \stationary
    \TransMat^{2 \mixingtime}) \cdot \vecnorm{\theta}{\myspan}
    \overset{(ii)}{\leq} \tfrac{1}{2}
    \vecnorm{\theta}{\myspan}, \label{eq:compostep-contraction-in-average-cost-policy-evaluation}
\end{multline}
where step (i) is a direct consequence of triangle inequality, and in step (ii), we
exploits the bound \mbox{$\totalvariation (\delta_x \TransMat^{2
  \mixingtime}, \stationary \TransMat^{2 \mixingtime}) \leq
\frac{1}{2} \totalvariation (\delta_x \TransMat^{\mixingtime},
\stationary \TransMat^{\mixingtime}) \leq \frac{1}{4}$}, obtained by
applying the mixing time condition~\ref{assume:mixing-time} twice.
Consequently, we see that the multi-step contraction
assumption~\ref{new--assume-compostep-contractive} holds if the
operator is composed $\compostep = 2 \mixingtime$ times.


\subsubsection{Estimation error upper bounds}

Having defined the norm $\vecnorm{\cdot}{\myspan}$ and the established
the multi-step contraction
property~\eqref{eq:compostep-contraction-in-average-cost-policy-evaluation},
we are ready to derive a guarantee for average-cost policy
evaluation. This involves the Gaussian random variable
\begin{align*}
 W \sim \mathcal{N} \left(0, \mathrm{cov} \left( \Hstoch (\thetastar)
 - \thetastar \right) \right),
\end{align*}
as well as $\noisegauss \mydefn \Exs \left[
  \vecnorm{\GaussNoise}{\myspan} \right]$.  For a given failure
probability $\delta \in (0,1)$, our result applies to the algorithm
with parameters
\begin{subequations}
\label{eqn:tuning-avg-reward}
\begin{align}
\label{eq:requirement-cor-avg-cost}
\stepsize = c_1 \Big\{\sqrt{ \numobs \log \abss{\states}} \cdot \log
(\tfrac{\numobs}{\delta}) \Big\}^{-1}, \quad \mbox{and} \quad \burnin
= \tfrac{c \mixingtime}{ \stepsize} \log(\tfrac{\numobs}{\delta}).
\end{align}
We also choose the initialization $\theta_0$ and the number of
restarts $\numRestarts$ such that
\begin{align}
\label{eqn:numrestarts-and-init-avg-PE}
 \log \left(
\tfrac{\norm{\theta_0 - \hpop (\theta_0)} \sqrt{\numobs}}{\noisegauss} \right) \leq c_0 \log \numobs\qquad \text{and} \qquad \numRestarts \geq 2 c_0 \log \numobs,
\end{align}
\end{subequations}
where $c, c_0, c_1$ are appropriate universal constants.  We have the following guarantee:
\begin{corollary}
\label{cor:average-cost}
Suppose Assumption~\ref{assume:mixing-time} holds, and the sample size
$\numobs$ is lower bounded as \mbox{$\frac{\numobs}{\log^2 \numobs}
  \geq c' \mixingtime^2 \log (\abss{\states}) \cdot \log (1 /
  \delta)$}. Then the estimate $\theta_\numobs$ from
Algorithm~\ref{AlgRootSGD}, obtained using tuning parameters
satisfying conditions~\eqref{eqn:tuning-avg-reward}, satisfies the
bound
\begin{multline}
\vecnorm{\theta_\numobs - \thetastar}{\myspan} \leq
\tfrac{c}{\sqrt{\numobs}} \Big \{ \Exs \big[ \vecnorm{(\Id -
    \TransMat)^\dagger W}{\myspan} \big] + \sqrt{\sup_{\state_1,
    \state_2 \in \states} \Exs \big[ ( (\delta_{\state_1} -
    \delta_{\state_2}) (\Id - \TransMat)^\dagger W)^2 \big]
  \log(1/\delta)} \Big \} \\
+ c \mixingtime \Big\{ \big[ \tfrac{\log |\states|}{\numobs}
  \big]^{3/4} \; \noisegauss + \tfrac{ \log |\states|}{\numobs}
\big(\sigma_r + \vecnorm{\thetastar}{\myspan}\big) \Big \}
\log^2(\tfrac{\numobs}{\delta}),
\end{multline}
with probability at least $1 - \delta$.
\end{corollary}
A few remarks are in order. First, the linear operator $(\Id -
\TransMat)$ is not invertible in $\real^\states$, with the all-one
vector lying in its nullspace. However, it is invertible in the
quotient space $\vecspace$, with the pseudo-inverse $(\Id -
\TransMat)^\dagger$ being a representation of its inverse in the
coordinate system of $\real^\states$, which appears in the
bound. Second, as with the previous two cases, the bound depends on
the size of state space only poly-logarithmically; it depends
quadratically on the mixing time $\mixingtime$, as shown through the
required lower bound on $\numobs$. Taking the
$\contraction$-discounted MRP as a special case of the average-cost
framework,\footnote{This can be done by adding an absorbing state
$\perp$ to the state space. At a rate of $(1 - \contraction)$, the
Markov process is killed and moved to the absorbing state. In such
case, the unique stationary distribution is the singleton at $\perp$,
and the mixing time assumption is satisfied with $\mixingtime =
\tfrac{c}{1 - \contraction}$ for universal constant $c > 0$.}
Corollary~\ref{cor:average-cost} improves the results of previous
work~\cite{khamaru2020temporal} in two aspects:
\bcar
\item Corollary~\ref{cor:average-cost} is valid whenever sample size
  satisfies $\numobs \gtrsim (1 - \contraction)^{-2}$ up to log factors,
  which improves the
  previous $(1 - \contraction)^{-3}$
  dependency from the past work;
\item The instance-dependent quantity in the
  paper~\cite{khamaru2020temporal} is replaced with an optimal one
  matching the local asymptotic minimax limit.  \ecar
These improvements are made possible by making use of the linear
structure in policy evaluation problems. More importantly,
Corollary~\ref{cor:average-cost} applies to a more general class of
problems, where the mixing time $\mixingtime$ replaces the role of
effective horizon.

In terms of other related work, the quadratic mixing time dependence
(i.e., sample size scaling as $\myorder(\mixingtime^2)$) matches that
of the paper~\cite{jin2020efficiently}. On one hand, our results are
more refined in that we give instance-dependent guarantees.  On the
other hand, their results apply to Markov decision processes with
actions. Thus, an open and interesting direction of future work is to
extend our instance-dependent bounds to the case of average-cost MDPs
with policy optimization.


\section{Proofs}
\label{sec:Main-proofs}

This section is devoted to the proofs of our main results---namely,
Theorems~\ref{new--thm:bellman-err-bound-general-nonlinear}
and~\ref{new--thm:main-nonlinear-general}---along with the associated
corollaries. So as to facilitate reading of the proofs, we reproduce
here the two main recursions that define the algorithm:
\begin{subequations}
\begin{align}
\myV_t & = \Hstoch_{t}(\theta_{t-1}) - \theta_{t-1} + \tfrac{t - 1}{t}
\left(\myV_{t - 1} - \Hstoch_t(\theta_{t-2}) + \theta_{t-2} \right),
\quad \mbox{and} \\
\theta_t & = \theta_{t - 1} + \stepsize \myV_t.
\end{align}
\end{subequations}
Throughout the proofs, we make use of the shorthand $\noisegauss =
\Exs[\norm{W}]$, and $\noisevar = \maxstd(W)$.

\subsection{Proof of Theorem~\ref{new--thm:bellman-err-bound-general-nonlinear}}
\label{subsec:proof-prop-bellman-err-general-nonlinear}

Our proof is based on a bootstrapping argument, and can be broken down
into four steps:
\bcar
    \item First, we establish recursions that relate
      $\norm{\hpop(\theta_t) - \theta_t}$ and $\norm{\myV_t}$.
    \item Second, we prove coarse upper bounds on
      $\norm{\hpop(\theta_t) - \theta_t}$ and $\norm{\myV_t}$.
    \item Third, starting with the sub-optimal bounds from Step 2, we
      iteratively refine them using a bootstrapping argument and the
      recursions from Step 1.
    \item In the fourth step, we improve higher-order terms in the
      bounds.  \ecar

\medskip

\noindent For the purposes of analysis, it is useful to define the
auxiliary sequence
\begin{align}
\label{eqn:zt-defn}
 \myZ_t & \defn \big \{ \hpop(\theta_{t - 1}) - \theta_{t-1} \big\} -
 \myV_t, \quad \mbox{for $t = \burnin, \burnin + 1, \cdots$}.
\end{align}
Our strategy is to control $\norm{\hpop(\theta_{t}) - \theta_{t}}$ by
proving upper bounds on $\norm{\myZ_{t + 1}}$ and $\norm{\myV_{t +
    1}}$.

\vspace{10pt}

Let $\ratetheta(t)$ and $\ratev(t)$, respectively, denote high
probability bounds on the quantities $\norm{\hpop(\theta_t) -
  \theta_t}$ and $\norm{\myV_t}$.  It is useful to introduce the
notion of an \emph{admissible sequence}: for some $\admissiblePar \geq
0$, the sequence $\{r(t)\}_{t \geq \burnin}$ is said to be
$\admissiblePar$-\emph{admissible} if
\begin{enumerate}
\item[(i)] The sequence $\{r (t)\}_{t \geq \burnin}$ is
  non-increasing.
\item[(ii)] The sequence $\{t^\admissiblePar \cdot r (t)\}_{t \geq
  \burnin}$ is non-decreasing.
\end{enumerate}
We say that the sequence is admissible if it is
$\admissiblePar$-admissible for some $\admissiblePar \geq 0$. For
notational simplicity, we sometimes use the sequences with time index
less than $\burnin$, in such cases, we denote $\ratev (t) \mydefn
\ratev (\burnin)$ and $\ratev (t) \mydefn \ratetheta (\burnin)$ for $t
\in [1, \burnin]$.

\medskip

Observe that $\admissiblePar$-admissible sequences are also
$\beta$-admissible sequences for any $\beta > \admissiblePar$.  For
the sake of notational convenience, we use the shorthands $\ratetheta$
and $\ratev$ to denote the estimate sequences $\{\ratetheta(t)\}_{t
  \geq \burnin}$ and $\{\ratev(t)\}_{t \geq \burnin}$, respectively.
Given an admissible pair $(\ratetheta, \ratev)$ and an integer
$\numobs > 0$,  define the events
\begin{align}
\label{eqn:eventDefn}  
\bigthetaEvent \mydefn \left\{ \sup_{\burnin \leq t \leq \numobs}
\frac{\norm{\hpop(\theta_t) - \theta_t}}{\ratetheta(t)} \leq 1
\right\}, \quad \mbox{and} \quad \bigvEvent & \mydefn \left\{
\sup_{\burnin \leq t \leq \numobs} \frac{\norm{\myV_t}}{\ratev(t)}
\leq 1 \right\}.
\end{align}
A key portion of our proof involves ensuring that the estimate
sequences $\ratetheta$ and $\ratev$ are
\mbox{$\admissiblePar$-admissible} for carefully chosen values of
$\admissiblePar$.  With these concepts and notation in place, we are
now ready to start the main argument.

\subsubsection{Step 1: Relation between $\norm{\hpop(\theta_t) - \theta_t}$ and $\norm{\myV_t}$}

From the definition~\eqref{eqn:zt-defn}, we have the relation
\mbox{$\hpop(\theta_{t}) - \theta_{t} = \myZ_{t + 1} + \myV_{t +
    1}$}. As mentioned before, we prove an upper bound on
$\norm{\hpop(\theta_{t}) - \theta_{t}}$ by proving upper bounds on
$\norm{\myZ_{t + 1}}$ and $\norm{\myV_{t + 1}}$.  We do so using two
auxiliary lemmas, the first of which depends on a stepsize $\stepsize$
satisfying the bound~\eqref{eqn:stepsize}---namely:
\begin{align}
\label{eq:vt-recursion-stepsize-requirement}  
\stepsize \leq \frac{(1 -
    \contraction)^2}{c \Lip^2 \DudGauss^2 (\dualBall, \myrhotil)
    \log \big(\tfrac{\numobs}{\delta} \big)}.
\end{align}
\begin{lemma}
\label{lemma-vt-recursive-bound}
Suppose that Assumptions~\ref{new--assume-pop-contractive},~
\ref{new--assume-noise-bound} and~\ref{new--assume-sample-operator}
are in force, and that $(\ratetheta, \ratev)$ are
\mbox{$\admissiblePar$-admissible} sequences for some $\admissiblePar
\in [0, 2]$.  Then given a stepsize $\stepsize$ satisfying the
bound~\eqref{eq:vt-recursion-stepsize-requirement} and a burn-in
period $\burnin \geq \frac{100}{(1 - \contraction) \stepsize}$
conditioned on the event $\bigvEvent \cap \bigthetaEvent$, for each $t
\in [\burnin, \numobs]$ we have
\begin{multline}
\label{eq:vt-recursion-bound}  
\norm{\myV_t} \leq \tfrac{1 + \contraction}{2} \ratev(t) + \tfrac{8}{t
  \stepsize} \ratetheta(t) + \tfrac{c}{\sqrt{\stepsize}} \Big\{
\noisegauss + \noisevar \sqrt{\log(\tfrac{1}{\delta})} \Big\} +
\tfrac{c \boundstar}{t} \Big \{ \log(\tfrac{1}{\delta}) +
\DudExp(\dualBall, \myrhotil) \Big \} + 6 (1 - \contraction)
\big(\tfrac{\burnin}{t} \big)^2 \norm{\myV_{\burnin}},
\end{multline}
with probability at least $1 - \delta$.
\end{lemma}

\noindent See
Section~\ref{subsubsec:proof-of-lemma-vt-recursive-bound} for the
proof of this lemma.

\begin{lemma}
\label{lemma-zt-martingale-bound}
Under the same conditions as Lemma~\ref{lemma-vt-recursive-bound}, for
each $t \in [ \burnin, \numobs]$, we have:
\begin{multline}
\label{eq:zt-martingale-bound-refined}  
  \norm{\myZ_t} \leq \tfrac{c}{\sqrt{t}} \Big\{ \noisegauss +
  \noisevar \sqrt{\log(\tfrac{1}{\delta})} \Big \} + \tfrac{
    \boundstar}{t} \Big \{ \DudExp(\dualBall, \myrhotil) +
  \log(\tfrac{1}{\delta}) \Big \} \\
  + \tfrac{c \Lip}{t} \Big \{ \constMGtype +
  \sqrt{\log(\tfrac{1}{\delta})} \Big \} \; \Big \{ \stepsize
  \Big(\sum_{s = \burnin}^{t - 1} s^2 \ratev^2(s)\Big)^{1/2} +
  \tfrac{1}{1 - \contraction} \Big( \sum_{s = 1}^{t - 1}
  \ratetheta^2(s) \Big)^{1/2} \Big \},
\end{multline}
with probability $1 -
\delta$.
\end{lemma}
\noindent This lemma is a special case of
Lemma~\ref{lemma-zt-martingale-bound-anorm}, which is proved in
Section~\ref{subsubsec:proof-lemma-zt-martingale-bound-anorm}. \\

Note that although the two lemmas are for a single time index $t \in
[\burnin, \numobs]$, it is easy to transform them to guarantees that
are uniform over $t \in [\burnin, \numobs]$. In particular, applying
a union bound for $t = \burnin, \burnin + 1, \cdots, \numobs$, and
by replacing $\delta$ with $\delta' = \delta / \numobs$, the
bounds~\eqref{eq:vt-recursion-bound}
and~\eqref{eq:zt-martingale-bound-refined} are valid uniformly over $t
\in [\burnin, \numobs]$.

We use these two lemmas in our bootstrapping argument.  In particular,
beginning with the relation $\hpop(\theta_{t}) - \theta_{t} = \myZ_{t +
  1} + \myV_{t + 1}$, applying the triangle inequality yields the bound
$\norm{\hpop(\theta_{t}) - \theta_{t} } \leq \norm{ \myZ_{t + 1}} +
\ratev(t + 1)$ on the event~$\bigvEvent$.  Our analysis shows that by starting with an initial
estimate $(\ratetheta(t), \ratev(t))$, the
bounds~\eqref{eq:vt-recursion-bound}
and~\eqref{eq:zt-martingale-bound-refined} allow us to obtain an
improved estimate $(\ratetheta^{+}(t), \ratev^{+}(t))$ such that
\begin{align*}
\norm{\hpop(\theta_{t}) - \theta_{t}} \;\; &\leq \;\;
\ratetheta^{+}(t) \;\; < \;\; \ratetheta(t), \quad \mbox {and} \\
\norm{\myV_{t}} \;\; & \leq \;\; \ratev^{+}(t) \;\; < \;\; \ratev(t)
\end{align*}
with high probability.  We quantify the improvement in
$(\ratetheta^{+}(t), \ratev^{+}(t))$, and repeatedly apply this
argument so as to ``bootstrap'' the bound and ultimately obtain sharp
estimates for $\ratetheta(t)$ and $\ratev(t)$.


\subsubsection{Step 2: Setup for the bootstrapping argument}

Throughout this step, we require that the estimate sequences
$\ratetheta$ and $\ratev$ be $\tfrac{1}{2}$-admissible and
$1$-admissible, respectively.  As shown in this section, these choices
allow us to obtain upper bounds on $\norm{\hpop(\theta_t) -
  \theta_t}$ and $\norm{\myV_t}$ that decay at the rates $1/\sqrt{t}$ and
$1/t$, respectively.

We assume that the pair $(\ratev^+, \ratetheta^+)$ satisfy the
initialization condition
\begin{subequations}
  \begin{align}
\label{eq:bootstrap-iteration-phase-i-init}    
\ratev^+(\burnin) & \geq \norm{\myV_{\burnin}}, \quad \mbox{and}\quad
\ratetheta^+ (\burnin) \geq \norm{\hpop(\theta_0) - \theta_0},
  \end{align}
and for each integer $t \in [\burnin, \numobs]$, the bounds
\begin{multline}
 \label{eq:bootstrap-iteration-phase-i-ratev}  
\ratev^+(t) \geq \frac{1 + \contraction}{2} \ratev(t) +
\tfrac{8}{\stepsize t} \ratetheta(t) + \tfrac{c}{t \sqrt{\stepsize}}
\Big \{ \noisegauss + \noisevar \sqrt{\log(\tfrac{\numobs}{\delta})}
\Big \} \\
+ \tfrac{c \boundstar}{t} \Big \{
\log(\tfrac{\numobs}{\delta}) + \DudExp(\dualBall, \myrhotil) \Big \} +
6 (1 - \contraction)\big( \tfrac{\burnin}{t} \big)^2
\norm{\myV_{\burnin}},
\end{multline}
and
\begin{multline}
\label{eq:bootstrap-iteration-phase-i-ratetheta}
\ratetheta^+(t) \geq \Big \{ 1 + c \stepsize \sqrt{t} \Lip
\big[\constMGtype + \sqrtlog{ (\tfrac{\numobs}{\delta}) } \big] \Big
\} \ratev(t) + \tfrac{2 c \Lip }{(1 - \contraction)\sqrt{t}}
\constMGtype \log (\tfrac{\numobs}{\delta}) \cdot \ratetheta(t)\\
+ \tfrac{c}{\sqrt{t}} \Big \{ \noisegauss + \noisevar
\sqrtlog{(\tfrac{\numobs}{\delta})} \Big \} + \tfrac{c \boundstar}{t}
\Big \{ \log(\tfrac{\numobs}{\delta}) + \DudExp(\dualBall, \myrhotil)
\Big \}.
\end{multline}
\end{subequations}
Under these conditions, by combining the
bounds~\eqref{eq:vt-recursion-bound}
and~\eqref{eq:zt-martingale-bound-refined} and applying a union bound
over $t \in [\burnin, \numobs]$, we find that
\begin{align*}
\Prob \Big[ \Event^{(\theta)}_\numobs (\ratetheta^+) \cap
  \Event^{(v)}_\numobs (\ratev^+) \Big] \geq \Prob \Big[
  \bigthetaEvent \cap \bigvEvent \Big] - \delta,
\end{align*}
valid for any pair $(\ratev, \ratetheta)$ that are $\tfrac{1}{2}$ and
$1$-admissible, respectively.

We consider sequences of a particular form $\ratev^{(i)} (t) =
\frac{\psi_v^{(i)}}{t \sqrt{\stepsize}}$ and \mbox{$\ratetheta^{(i)}
  (t) = \frac{\psi_\theta^{(i)}}{\sqrt{t}}$}, for pairs of positive
reals $\big( \psi_v^{(i)}, \psi_\theta^{(i)} \big)$ independent of
$t$. Clearly, with such forms, the sequence $\ratetheta^{(i)}$ is
$\frac{1}{2}$-admissible, and the sequence $\ratev^{(i)}$ is
$1$-admissible. However, if we directly substitute the sequences
$\big( \ratev^{(i)} (t), \ratev^{(i)} (t) \big)$ of such forms into
the
relations~\eqref{eq:bootstrap-iteration-phase-i-ratev}-\eqref{eq:bootstrap-iteration-phase-i-ratetheta},
the resulting sequences $(\ratetheta^+, \ratev^+)$ are no longer be of
the desired form. So in order to unify the coefficients in
equations~\eqref{eq:bootstrap-iteration-phase-i-ratev}-\eqref{eq:bootstrap-iteration-phase-i-ratetheta}
into the same time scale, given $\stepsize > 0$, we define
the burn-in time
\begin{subequations}
\begin{align}
  \label{eq:burn-in-requirement-in-bootstrap-phase-i}  
  \burnin = \tfrac{c}{(1 - \contraction)^2 \stepsize} \log
  (\tfrac{\numobs}{\delta}).
\end{align}
For each $t = \burnin, \burnin +1 \ldots$, the coefficients in
\eqref{eq:bootstrap-iteration-phase-i-ratev} and~\eqref{eq:bootstrap-iteration-phase-i-ratetheta} then satisfy the
bounds
\begin{align}
\label{eq:coefficient-bounds-in-bootstrap-phase-i}  
 \tfrac{8}{\stepsize t} \leq \tfrac{1 - \contraction}{3} \cdot \tfrac{1}{
   \sqrt{\stepsize t}}, \quad
 \tfrac{1}{\sqrt{\stepsize t}} \leq \tfrac{1 - \contraction}{6},\quad
 \mbox{and} \quad
\tfrac{1}{(1 - \contraction)\sqrt{t}} \log (\tfrac{\numobs}{\delta})
\leq \sqrt{\stepsize},
\end{align}
\end{subequations}
Therefore, if we construct a two-dimensional vector sequence
$\psi^{(i)} = \begin{bmatrix} \psi_v^{(i)} & \psi_\theta^{(i)}
\end{bmatrix}^T$ satisfying the recursive relation $\psi^{(i+1)} = \Qmat \psi^{(i)} + \bvec$, where
\begin{subequations}
\label{eq:bootstrap-scalar}    
\begin{align}
 \Qmat & \mydefn \begin{bmatrix} \tfrac{1 + \contraction}{2} &
   \tfrac{1 - \contraction}{3}\\ \tfrac{1 - \contraction}{6} + c \Lip
   \constMGtype \sqrt{\stepsize} \log(\tfrac{\numobs}{\delta}) & 2 c
   \Lip \constMGtype \log (\tfrac{\numobs}{\delta}) \sqrt{\stepsize}
 \end{bmatrix}, \quad \mbox{and} \\
 \bvec & \mydefn \begin{bmatrix} c \big \{ \noisegauss + \noisevar
   \sqrt{\log( \tfrac{\numobs}{\delta})} \big \} + c \boundstar
   \sqrt{\stepsize} \big \{ \log(\tfrac{\numobs}{\delta}) +
   \DudExp(\dualBall, \myrhotil) \big \} + (1 - \contraction) \burnin
   \sqrt{\stepsize} \norm{\myV_{\burnin}} \\
\big \{ \noisegauss + \noisevar \sqrt{\log( \tfrac{\numobs}{\delta})}
\big \} + c \boundstar \sqrt{\stepsize} \big \{ \log(
\tfrac{\numobs}{\delta}) + \DudExp(\dualBall, \myrhotil) \big \} +
\sqrt{\burnin} \norm{\hpop(\theta_0) - \theta_0}
 \end{bmatrix}
\end{align}
\end{subequations}
 satisfy the
requirement~\eqref{eq:bootstrap-iteration-compostep-phase-i}.  Thus,
we are led to the probability bound
\begin{align}
\label{eqn:boot-strap-prob-recursion-Thm1}
  \Prob \Big[ \Event^{(\theta)}_\numobs \big( \ratetheta^{(i + 1)}
    \big) \cap \Event^{(v)}_\numobs \big( \ratev^{(i + 1)} \big)\Big]
  \geq \Prob \Big[ \Event^{(\theta)}_\numobs \big( \ratetheta^{(i)}
    \big) \cap \Event^{(v)}_\numobs \big( \ratev^{(i)} \big)\Big] -
  \delta,
\end{align}
for the sequences $\ratetheta^{(i)} (t) = \psi_\theta^{(i)} /
\sqrt{t}$ and $\ratev^{(i)} (t) = \psi_v^{(i)} / (\sqrt{\stepsize}
t)$. In order to initialize the argument, we need a coarse bound on
the pair $(\norm{\myV_t}, \norm{\hpop(\theta_t) - \theta_t})$; the
following lemma provides the requisite bound:
\begin{lemma}
\label{lemma:coarse-bound-gronwall}
Under Assumptions~\ref{new--assume-noise-bound}
and~\ref{new--assume-sample-operator}, we have
    \begin{align*}
     \norm{\theta_t - \thetastar} + \norm{\myV_t} \leq e^{1 + \Lip
       \stepsize t} \left( \boundstar + \norm{\theta_0 - \thetastar}
     \right),
    \end{align*}
almost surely for each $t = 0, 1, 2, \ldots$.
\end{lemma}
\noindent See
Appendix~\ref{subsubsec:proof-lemma-coarse-bound-gronwall} for the
proof of this claim. \\

Based on Lemma~\ref{lemma:coarse-bound-gronwall}, it follows that for
each integer $t \in [1, \numobs]$, we have (almost surely) the bound
\begin{align*}
\norm{\myV_t} & \leq \ratev^{(0)}(t) \mydefn {\tfrac{\numobs}{t}} e^{1 +
  \Lip \stepsize t} \big \{ \boundstar + \norm{\theta_0 - \thetastar}
\big \}, \quad \text{and} \\
\norm{\hpop(\theta_t) - \theta_t} & \stackrel{(i)}{\leq}
\ratetheta^{(0)} (t) \mydefn 2\cdot \sqrt{\tfrac{\numobs}{t}} e^{1 +
  \Lip \stepsize t} \Big \{ \boundstar + \norm{\theta_0 - \thetastar}
\Big \},
\end{align*}
where step (i) follows from the bound $\norm{\hpop(\theta_t) -
  \theta_t} \leq \norm{\theta_t - \thetastar} + \norm{\hpop(\theta_t)
  - \hpop(\thetastar)} \leq 2 \cdot \norm{\theta_t - \thetastar}$.

By construction, the sequences $\ratev^{(0)}$ and $\ratetheta^{(0)}$
are $1$-admissible and $\tfrac{1}{2}$-admissible, respectively, and by
Lemma~\ref{lemma:coarse-bound-gronwall}, the event
$\Event_\numobs^{(\theta)} \big( \ratetheta^{(0)} \big) \cap
\Event_\numobs^{(v)} \big( \ratev^{(0)} \big)$ happens almost
surely.

\subsubsection{Step 3: Bootstrapping step}

Recursing the bound~\eqref{eqn:boot-strap-prob-recursion-Thm1} for $i$
steps yields
\begin{align*}
\Prob \big[ \Event_\numobs^{(v)}(\ratev^{(i)}) \cap
  \Event_\numobs^{(\theta)}(\ratetheta^{(i)}) \big] & \geq \Prob \big[
  \Event_\numobs^{(v)}(\ratev^{(0)}) \cap
  \Event_\numobs^{(\theta)}(\ratetheta^{(0)}) \big] - i \delta = 1 - i
\delta.
\end{align*}
It remains to analyze
the 
sequence
$\psi^{(i)} = \begin{bmatrix} \psi_v^{(i)} & \psi_\theta^{(i)}
\end{bmatrix}^T$ as the number of bootstrap steps $i$
increases.  We do so by analyzing the recursion relation $\psi^{(i+1)}
= \Qmat \psi^{(i)} + \bvec$ with the matrix $\Qmat$ given in equation~\eqref{eq:bootstrap-scalar}.

Observe that the stepsize
condition~\eqref{eq:vt-recursion-stepsize-requirement} ensures that
\begin{align}
\label{eq:stepsize-requirement-bootstrap-phase-i}  
c \Lip \constMGtype \log (\tfrac{\numobs}{\delta}) \cdot
\sqrt{\stepsize} \leq \tfrac{ 1 - \contraction}{6}.
\end{align}
Consequently, the matrix $\Qmat$ from
equation~\eqref{eq:vt-recursion-stepsize-requirement} is entrywise
upper bounded by the matrix
\begin{align*}
\QmatTil = \begin{bmatrix} \tfrac{1 + \contraction}{2} & \tfrac{1 -
    \contraction}{3}\\ \tfrac{1 - \contraction}{3} & \tfrac{1}{2}
    \end{bmatrix}
\end{align*}
This fact implies that for any vector $u \in \real^2$ with
non-negative entries, we have the upper bound \mbox{$\Qmat u
  \preceqOrthant \QmatTil u$,} where $\preceqOrthant$ denotes the
orthant ordering. Straightforward calculation yields the bound
$\opnorm{\QmatTil} \leq 1 - \tfrac{1 - \contraction}{8}$.  Putting
together the pieces, we find that for each $N = 1, 2, \ldots$,
conditioned on the event $\Event_\numobs^{(v)}(\ratev^{(N)}) \cap
\Event_\numobs^{(\theta)}( \ratetheta^{(N)})$, we have
\begin{align*}
\psi^{(N)} &= \Big(\sum_{i = 0}^{N - 1} \Qmat^i \Big) b_\psi +
\Qmat^N \begin{bmatrix} \psi_v^{(0)} \\ \psi_\theta^{(0)}
\end{bmatrix} \preceqOrthant \Big(\sum_{i = 0}^{N - 1} \QmatTil^i \Big)
\bvec + \QmatTil^N \begin{bmatrix} \psi_v^{(0)} \\ \psi_\theta^{(0)}
    \end{bmatrix}\\
    & \preceqOrthant (I - \QmatTil)^{-1} \bvec + e^{\tfrac{(1 -
    \contraction)}{8} N} \big( \psi_v^{(0)} + \psi_\theta^{(0)} \big)
\bm{1}_2.
\end{align*}
We take $N = \lceil \tfrac{c \Lip \numobs}{1 - \contraction} \log
\numobs \rceil$.  Replacing $\delta$ with $\delta / N$ and
substituting into the above inequalities then yields
\begin{subequations}
\label{eq:bootstrap-phase-i-consequence}
\begin{multline}
 \label{eq:bootstrap-phase-i-consequence-v}   
    t \sqrt{\stepsize} \cdot \norm{\myV_t} \leq \psi_v^{(N)} \leq
    \tfrac{c}{1 - \contraction} \Big \{ \noisegauss + \noisevar
    \sqrt{\log (\tfrac{\numobs}{\delta})} \Big \} + \tfrac{c
      \boundstar \sqrt{\stepsize}}{1 - \contraction} \Big \{ \log
    (\tfrac{\numobs}{\delta}) + \DudExp(\dualBall, \myrhotil) \Big \} \\
 + c \burnin \sqrt{\stepsize} \norm{\myV_{\burnin}} + \sqrt{\burnin}
 \norm{\hpop(\theta_0) - \theta_0},
\end{multline}
and
\begin{multline}
\label{eq:bootstrap-phase-i-consequence-theta}  
 \sqrt{t} \norm{\hpop(\theta_t) - \theta_t} \leq \psi_\theta^{(N)}
 \leq c \Big\{ \noisegauss + \noisevar
 \sqrt{\log(\tfrac{\numobs}{\delta})} \Big \} + c \boundstar
 \sqrt{\stepsize} \Big \{ \log(\tfrac{\numobs}{\delta}) +
 \DudExp(\dualBall, \myrhotil) \Big \} \\
+ c \burnin (1 - \contraction) \sqrt{\stepsize} \norm{\myV_{\burnin}} +
\sqrt{\burnin} \norm{\hpop(\theta_0) - \theta_0},
\end{multline}
\end{subequations}
with probability at least $1 - \delta$, uniformly for each $t \in
\burnin, \burnin + 1, \cdots, \numobs$. \\

\noindent It remains to provide upper bounds on $\norm{\myV_{\burnin}}$.
\begin{lemma}
\label{lemma:burnin-period-bound}
Under Assumptions~\ref{new--assume-pop-contractive}
and~\ref{new--assume-noise-bound}, and a burn-in period given by
equation~\eqref{eq:burn-in-requirement-in-bootstrap-phase-i}, we have
\begin{align*}
  \norm{\myV_{\burnin}} \leq 2 \norm{\hpop(\thetainit) - \thetainit} +
\tfrac{c}{\sqrt{\burnin}} \Big \{ \noisegauss + \noisevar
\sqrt{\log(\tfrac{1}{\delta})} \Big \} + \tfrac{c\boundstar}{\burnin}
\Big \{ \DudExp(\dualBall, \myrhotil) + \log(\tfrac{1}{\delta}) \Big \}.
\end{align*}
with probability at least $1 - \delta$.
\end{lemma}
\noindent See~\Cref{sec:proof-of-lemma:burnin-period-bound} for the
proof. \\

\noindent Combining~\Cref{lemma:burnin-period-bound} and
bound~\eqref{eq:bootstrap-phase-i-consequence-v}, we find that
\begin{multline}
    \label{eqn:step-i-bound-h}   
  \norm{\hpop(\theta_t) - \theta_t} \leq \tfrac{c}{\sqrt{t}} \Big(
  \noisegauss + \noisevar \sqrt{\log(\tfrac{\numobs}{\delta})} \Big) +
  \tfrac{c \boundstar \sqrt{\stepsize}}{\sqrt{t}} \Big \{ \log
  (\tfrac{\numobs}{\delta}) + \DudExp(\dualBall, \myrhotil) \Big \} \\ +
  \norm{\hpop(\theta_0) - \theta_0} \tfrac{\sqrt{\burnin}}{\sqrt{t}}
  \log^{\tfrac{3}{2}}(\tfrac{\numobs}{\delta})
\end{multline}
with probability at least $1 - \delta$, uniformly for all integers $t
\in [\burnin, \numobs]$.

Although this bound has optimal dependence on $\noisegauss + \noisevar
\sqrt{\log(\tfrac{\numobs}{\delta})}$, its dependence on the terms
$\norm{\hpop(\theta_0) - \theta_0}$ and $\DudExp(\dualBall,
\myrhotil)$ and $\log(\numobs/\delta)$ in the
bound~\eqref{eqn:step-i-bound-h} can be sharpened. This motivates the
second phase of the bootstrap argument.


\subsubsection{Step 4: Improving higher-order terms}

Given the pair
$(\psi_v^{(N)}, \psi_\theta^{(N)})$ defined by\footnote{We redefine
$(\psi_v^{(N)}, \psi_\theta^{(N)})$ using the right-hand side of
\eqref{eq:bootstrap-phase-i-consequence}} the right-hand side of
\eqref{eq:bootstrap-phase-i-consequence}, conditioned on the event
$\bigthetaEvent \cap \bigvEvent$ with $\ratev(t) =
\tfrac{\psi_v^{(N)}}{t \sqrt{\stepsize}}$ and the sequence
$\ratetheta(t) = \tfrac{\psi_\theta^{(N)}}{\sqrt{t}}$,
invoking the bound~\eqref{eq:zt-martingale-bound-refined} from Lemma~\ref{lemma-zt-martingale-bound}
we have 
\begin{multline*}
 \norm{\hpop(\theta_t) - \theta_t} \leq \tfrac{c}{\sqrt{t}} \big(
 \noisegauss + \noisevar \sqrt{\log(\tfrac{1}{\delta})} \big) +
 \tfrac{c \boundstar}{t} \big( \DudExp(\dualBall, \myrhotil) + \log(
 \tfrac{1}{\delta}) \big) \\
 + \Big \{ \tfrac{1}{t} + c \tfrac{\stepsize \Lip \constMGtype}{
   \sqrt{t}} \cdot \log(\tfrac{\numobs}{\delta}) \Big \} \tfrac{
   \psi_v^{(N)}}{\sqrt{\stepsize}} + 2c \tfrac{ \Lip \constMGtype}{(1
   - \contraction)t} \log(\tfrac{\numobs}{\delta}) \cdot
 \psi_\theta^{(N)} \\
 \leq \tfrac{c}{\sqrt{t}} \Big \{ \noisegauss + \noisevar
 \sqrt{\log(\tfrac{1}{\delta})} \Big \} + c' \boundstar \Big \{
 \tfrac{1}{(1 - \contraction) t} + \tfrac{\stepsize \Lip
   \constMGtype}{ (1 - \contraction) \sqrt{t}} \cdot \log
 (\tfrac{\numobs}{\delta}) \Big \} \; \Big \{ \DudExp(\dualBall,
 \myrhotil) + \log (\tfrac{\numobs}{\delta}) \Big \} \\
+ c' \Big \{ \tfrac{1}{(1 - \contraction) t\sqrt{\stepsize}} +
\tfrac{\Lip \constMGtype \sqrt{\stepsize}}{(1 - \contraction)\sqrt{t}}
\cdot \log(\tfrac{\numobs}{\delta}) + \tfrac{ \Lip \constMGtype}{(1 -
  \contraction)t} \log(\tfrac{\numobs}{\delta}) \Big \} \Big \{
\noisegauss + \noisevar \sqrt{\log (\tfrac{\numobs}{\delta})} \Big \}
\\
+ c' \Big \{ \tfrac{1}{t \sqrt{\stepsize}} + \tfrac{\sqrt{\stepsize}
  \Lip \constMGtype}{ \sqrt{t}} \cdot \log(\tfrac{\numobs}{\delta}) +
\tfrac{ \Lip \constMGtype}{(1 - \contraction)t}
\log({\tfrac{\numobs}{\delta}}) \Big \}) \sqrt{\burnin}
\norm{\hpop(\theta_0) - \theta_0} \\
+ c' \Big \{ \tfrac{1}{t} + \tfrac{\stepsize \Lip \constMGtype}{
  \sqrt{t}} \cdot \log(\tfrac{\numobs}{\delta}) + \tfrac{
  \sqrt{\stepsize} \Lip \constMGtype}{t} \log(\tfrac{\numobs}{\delta})
\Big \} \burnin \norm{\myV_{\burnin}},
\end{multline*}
which holds with probability at least $1 - \delta$.  Given the burn-in
period satisfying
equation~\eqref{eq:burn-in-requirement-in-bootstrap-phase-i} and
stepsize satisfying
equation~\eqref{eq:stepsize-requirement-bootstrap-phase-i}, by
combining with the bound on $\norm{\myV_{\burnin}}$ from
Lemma~\ref{lemma:burnin-period-bound}, we find that
$\norm{\hpop(\theta_t) - \theta_t} \leq \widetilde{\ratetheta}(t)$
with at least probability $1 - \delta$, uniformly for any integer $t
\in [\numobs]$, where
\begin{multline}
\label{eq:thetat-rate-final}  
\widetilde{\ratetheta}(t) \mydefn \tfrac{c_1}{\sqrt{t}} \Big \{
\noisegauss + \noisevar \sqrt{\log (\tfrac{\numobs}{\delta})} \Big\} +
\tfrac{c_2 \boundstar}{1 - \contraction} \Big \{ \tfrac{1}{t} +
\tfrac{\stepsize \Lip }{\sqrt{t}} \cdot \constMGtype \log
(\tfrac{\numobs}{\delta}) \Big \} \; \Big \{ \DudExp(\dualBall,
\myrhotil) + \log (\tfrac{\numobs}{\delta}) \Big\} \\
+ c_2 \Big \{ \tfrac{\stepsize \burnin \Lip }{ \sqrt{t}} \cdot
\constMGtype \log(\tfrac{\numobs}{\delta}) + \tfrac{\burnin}{t} \Big\}
\norm{\hpop(\theta_0) - \theta_0}.
\end{multline}

By substituting our upper bound in terms of $\widetilde{\ratetheta}$
into equation~\eqref{eq:bootstrap-iteration-phase-i-ratev}, we obtain
a recursive inequality that takes an admissible sequence $\ratev$ and
generates a sequence $\ratev^+$ such
that
\begin{align*}
\Prob \big[ \Event_\numobs^{(v)} (\ratev^+) \big] \geq \Prob \big[
  \Event_\numobs^{(v)} (\ratev) \big] - \delta.
\end{align*}
For any positive integer $N_1$, we can apply the recursive inequality
for $N_1$ times with $\delta' = \delta / N_1$; doing so yields a
sharper bound for $\norm{\myV_t}$.  In particular, with probability at
least $1 - \delta$, we have
\begin{multline*}
\norm{\myV_t} \leq \tfrac{2}{1 - \contraction} \Big \{ \tfrac{c}{t
  \sqrt{\stepsize}} \Big[ \noisegauss + \noisevar
  \sqrt{\log(\tfrac{\numobs N_1}{\delta})} \Big] + \tfrac{c
  \boundstar}{t} \Big[ \log(\tfrac{\numobs N_1}{\delta}) +
  \DudExp(\dualBall, \myrhotil) \Big] \Big \} \\
+ \tfrac{8}{(1 - \contraction) \stepsize t} \widetilde{\ratetheta}(t)
+ \big( \tfrac{\burnin}{t} \big)^2 \norm{\myV_{\burnin}} + \big( \tfrac{1
  + \contraction}{2} \big)^{N_1} \cdot \tfrac{\psi_v^{(N)}}{t
  \sqrt{\stepsize}}.
\end{multline*}

We take $N_1 \mydefn \lceil\tfrac{10 \log \numobs}{1 - \contraction} \rceil$, and a
stepsize and burn-in period satisfying the
conditions~\eqref{eq:burn-in-requirement-in-bootstrap-phase-i}
and~\eqref{eq:stepsize-requirement-bootstrap-phase-i}. With these
choices, some algebra yields $\norm{\myV_t} \leq \widetilde{\ratev}(t)$
holds with probability at least $1 - \delta$, uniformly for each
integer $t \in [\burnin, \numobs]$, where
\begin{align}
  \label{eq:vt-rate-final}
  \widetilde{\ratev}(t) & \mydefn \tfrac{c'}{1 - \contraction} \Big\{
  \tfrac{1}{t \sqrt{\stepsize}} \big[ \noisegauss + \noisevar
    \sqrt{\log (\tfrac{\numobs}{\delta})} \big] +
  \tfrac{\boundstar}{t} \big[ \log (\tfrac{\numobs}{\delta}) +
    \DudExp(\dualBall, \myrhotil) \big] \Big \} + 2
  \big(\tfrac{\burnin}{t}\big)^2 \norm{\thetainit -
    \hpop(\thetainit)}.
\end{align}
%


It can be seen that the sequences $\widetilde{\ratev}$ and
$\widetilde{\ratetheta}$ are $2$-admissible.  Substituting their
definitions into the bound~\eqref{eq:zt-martingale-bound-refined} from Lemma~\ref{lemma-zt-martingale-bound}
we find that the inequality
\begin{multline*}
    \norm{\myZ_{t}} \leq \tfrac{c}{\sqrt{t}} \Big \{ \noisegauss +
    \noisevar \sqrt{\log(\tfrac{1}{\delta})} \Big \} + \tfrac{c
      \boundstar}{t} \Big \{ \DudExp(\dualBall, \myrhotil) +
    \log(\tfrac{1}{\delta}) \Big \} \\
+ \tfrac{c \Lip}{t} \Big \{ \constMGtype +
\sqrt{\log(\tfrac{1}{\delta})} \Big \} \; \Big \{ \stepsize
\mybiggersqrt{\sum_{s = \burnin}^{t - 1} s^2 \ratev^2(s)} + 
\frac{1}{1 -\contraction} \cdot \mybiggersqrt{\sum_{s = 1}^{t - 1} \ratetheta^2(s)}
\Big \}
\end{multline*}
holds with probability at least $1 - \delta$.

Under the stepsize and burn-in period 
conditions~\eqref{eq:burn-in-requirement-in-bootstrap-phase-i}
and~\eqref{eq:stepsize-requirement-bootstrap-phase-i},
some algebra yields:
\begin{multline*}
\norm{\myZ_{t}} \leq \tfrac{c}{\sqrt{t}} \big( \noisegauss + \noisevar
\sqrt{\log(\tfrac{1}{\delta})} \big) \\ + c \boundstar \Big \{
\tfrac{1}{t} + \tfrac{\stepsize \Lip }{ \sqrt{t}} \cdot \constMGtype
\log(\tfrac{\numobs}{\delta}) \Big \} \; \Big\{ \DudExp(\dualBall,
\myrhotil) + \log(\tfrac{1}{\delta}) \Big \} + c \tfrac{\burnin}{t}
\norm{\thetainit - \hpop(\thetainit)}.
\end{multline*}

Combining with equation~\eqref{eq:vt-rate-final} yields the upper
bound
\begin{multline*}
\norm{\hpop(\theta_t) - \theta_t} \leq \tfrac{c}{\sqrt{t}} \big(
\noisegauss + \noisevar \sqrt{\log(\tfrac{1}{\delta})} \big) + {c
  \boundstar} \Big \{ \tfrac{1}{t} + \tfrac{\stepsize \Lip }{
  \sqrt{t}} \constMGtype \log (\tfrac{\numobs}{\delta}) \Big \} \;
\Big \{ \DudExp(\dualBall, \myrhotil) + \log(\tfrac{1}{\delta}) \Big \}
\\
+ c\tfrac{\burnin}{t} \norm{\thetainit - \hpop(\thetainit)},
\end{multline*}
which completes the proof of the
Theorem~\ref{new--thm:bellman-err-bound-general-nonlinear}.

Besides, we also note that by taking a union bound over time steps $t \in \{ \burnin, \burnin +
1, \ldots, \numobs \}$, we have the lower bound
$\Prob \big[ \Event^{(\theta)}_\numobs(\ratetheta^*)
\big] \geq 1 - \delta$,
where
\begin{multline*}
    \ratetheta^*(t) \mydefn \tfrac{c}{\sqrt{t}} \Big \{ \noisegauss +
    \noisevar \sqrt{\log (\tfrac{\numobs}{\delta}) } \Big \} \\ + {c
      \boundstar} \Big \{ \tfrac{1}{t} + \tfrac{\stepsize \Lip
      \constMGtype}{ \sqrt{t}} \log (\tfrac{\numobs}{\delta}) \Big \}
    \; \Big \{ \DudExp(\dualBall, \myrhotil) + \log
    (\tfrac{\numobs}{\delta}) \Big \}  + c \tfrac{\burnin}{t}
    \norm{\thetainit - \hpop(\thetainit)}.
\end{multline*}

\renewcommand{\UpsSet}{\mathcal{S}}

\subsection{Proof of Corollary~\ref{cor:anorm-suboptimality-bound}}
\label{sec:proof-of-cor-anorm-suboptimality-bound}

The proof of this corollary is based on a modification of
Lemma~\ref{lemma-zt-martingale-bound}.  We introduce the shorthand
\begin{align*}
\noisegauss & \defn \Exs[\norm{\GaussNoise}], \qquad \qquad \qquad
\quad \noisegausssemi \defn \Exs[\anorm{\GaussNoise}], \\
\noisevar & \defn \sqrt{ \sup_{u \in \dualBall} \Exs \big[
    \inprod{u}{W}^2 \big]} \quad \text{and} \quad \noisevarsemi \defn
\sqrt{ \sup_{u \in \MySet} \Exs \big[ \inprod{u}{W}^2 \big]}.
\end{align*}
We begin by stating a lemma---a generalization of
Lemma~\ref{lemma-zt-martingale-bound}---that bounds the supremum of an
averaged process. In the proof of
Corollary~\ref{cor:anorm-suboptimality-bound}, we only use a special
case of Lemma~\ref{lemma-zt-martingale-bound-anorm}, but the
generality is useful later.

Recall the events $\bigthetaEvent$ and $\bigvEvent$ defined in
equation~\eqref{eqn:eventDefn}.  Given a bounded symmetric convex set
$\UpsSet \subseteq \vecspace^*$, we define the dimension factor
$\aUB_\UpsSet \mydefn \sup_{u \in \UpsSet} \dualnorm{u}$. Moreover, we
assume that there exists a constant $\MuCon > 0$ such that
\begin{align}
\label{EqnMucon}  
 \norm{\theta - \thetastar} \leq \tfrac{1}{\MuCon} \norm{\hpop(\theta) - \theta}  \;
  \mbox{for any $\theta
   \in \vecspace$.}
\end{align}
We point out that under assumption~\ref{new--assume-pop-contractive},
the last condition is satisfied for $\mu = 1 - \contraction$.  The
condition~\eqref{EqnMucon} also allows us to analysis behavior of
operators which satisfies a multi-step contraction assumption
~\ref{new--assume-compostep-contractive} (cf. the proof of
Theorem~\ref{thm:bellman-err-bound-general-linear}).

\begin{lemma}
\label{lemma-zt-martingale-bound-anorm}
Suppose that the Assumptions~\ref{new--assume-sample-operator} and
\ref{new--assume-noise-bound} are in force, the sequences $\ratetheta$
and $\ratev$ are \mbox{$\admissiblePar$-admissible} for some
$\admissiblePar \in (0,2]$, and condition~\eqref{EqnMucon} holds.
Then conditioned on the event $\Event_\numobs^{(\theta)}(\ratetheta)
\cap \Event_\numobs^{(v)}(\ratev)$, we have
\begin{multline}
\label{eq:zt-martingale-bound-refined-anorm}        
  \sup_{u \in \UpsSet} \inprod{u}{\myZ_t} \leq \tfrac{c}{\sqrt{t}}
  \Big \{ \Exs \big[ \sup_{u \in \UpsSet} \inprod{u}{W} \big] +
  \mybiggersqrt{\sup_{u \in \UpsSet} \Exs \big[ \inprod{u}{W}^2 \big]
    \log(\tfrac{1}{\delta})} \Big \} + \tfrac{\aUB_{\UpsSet}
    \boundstar}{t} \Big \{ \DudExp(\dualBall, \myrhotil) +
  \log(\tfrac{1}{\delta}) \Big \} \\
  + \tfrac{c \aUB_\UpsSet \Lip}{t} \Big \{ \constMGtype + \sqrt{\log
    (\tfrac{1}{ \delta})} \Big \} \Big \{ \stepsize
  \mybiggersqrt{\sum_{s = \burnin}^{t - 1} s^2 \ratev^2(s)} +
  \frac{1}{\MuCon} \cdot \mybiggersqrt{\sum_{s = 1}^{t - 1}
    \ratetheta^2(s)} \Big \},
\end{multline}
with probability at least $1 - \delta$, uniformly for all integers $t
\in [\burnin, \numobs]$.
\end{lemma}
\noindent See
Section~\ref{subsubsec:proof-lemma-zt-martingale-bound-anorm} for the
proof of this lemma.

Taking this lemma as given, we now proceed with proof of
Corollary~\ref{cor:anorm-suboptimality-bound}. As mentioned before,
under assumption~\ref{new--assume-compostep-contractive},
condition~\eqref{EqnMucon} is satisfied with $\MuCon = 1 -
\contraction$.  Applying Lemma~\ref{lemma-zt-martingale-bound-anorm}
with $\UpsSet = \MySet$ implies that
\begin{multline}
\anorm{\myZ_t} \leq \tfrac{c}{\sqrt{t}} \Big \{ \noisegausssemi +
\noisevarsemi \sqrt{\log(\tfrac{1}{\delta})} \Big \} + c \tfrac{\aUB
  \boundstar}{t} \Big \{ \DudExp(\dualBall, \myrhotil) +
\log(\tfrac{1}{\delta}) \Big \} \\
+ \tfrac{c \aUB \Lip}{t} \Big \{ \constMGtype +
\sqrt{\log(\tfrac{1}{\delta})} \Big \} \; \Big \{ \stepsize
\mybiggersqrt{\sum_{s = \burnin}^{t - 1} s^2 \ratev^2(s)} +
\tfrac{1}{1 - \contraction} \mybiggersqrt{\sum_{s = 1}^{t - 1}
  \ratetheta^2(s)} \Big \}.
  \label{eqn:simplified-zt-anorm-bd}
\end{multline}
Now all we have to do is substitute an appropriate value of the
sequences $\ratev$ and $\ratetheta$.  Note that the estimate sequences
$\widetilde{\ratev}$ and $\widetilde{\ratetheta}$ from
equations~\eqref{eq:vt-rate-final} and~\eqref{eq:thetat-rate-final},
respectively, are $2$-admissible; moreover, they provide upper bounds
on the quantities $\norm{\myV_t}$ and $\norm{\myZ_t}$ respectively. Next,
using the stepsize and burn-in
conditions~\eqref{eq:stepsize-requirement-bootstrap-phase-i}
and~\eqref{eq:burn-in-requirement-in-bootstrap-phase-i}, we find that
\begin{multline*}
\mybiggersqrt{\frac{1}{t}\sum_{s = \burnin}^{t - 1} s^2
  \widetilde{\ratev}^2(s)} \leq \tfrac{c}{(1 - \contraction)
  \sqrt{\stepsize}} \Big \{ \noisegauss + \noisevar \sqrt{\log
  (\tfrac{\numobs}{\delta})} \Big \} \\
+ \tfrac{c \cdot \boundstar}{(1 -
  \contraction)} \Big \{ \log (\tfrac{\numobs}{\delta}) +
\DudExp(\dualBall, \myrhotil) \Big \} + \tfrac{2 c \: \burnin^{3/2} \,
  \norm{\thetainit - \hpop(\thetainit)}}{t},
\end{multline*}
and
\begin{multline*}
\mybiggersqrt{\sum_{s = 1}^{t - 1} \widetilde{\ratetheta}^2(s)}
\leq c \cdot \Big \{ \noisegauss + \noisevar
\sqrt{\log(\tfrac{n}{\delta})} \Big \} \cdot \sqrt{\log t} \\
   + \tfrac{c \boundstar}{1 - \contraction} \Big \{
   \tfrac{1}{\sqrt{\burnin}} + \stepsize \LipCon \sqrt{\log t} \big[
     \constMGtype + \sqrt{\log (\tfrac{1}{\delta})} \big] \Big \} \;
   \Big \{ \DudExp(\dualBall, \myrhotil) + \log
   (\tfrac{\numobs}{\delta}) \Big \} \\
   + c \Big \{ \stepsize \burnin \LipCon \sqrt{\log(t)} \big[
     \constMGtype + \sqrt{\log(\tfrac{1}{\delta})} \big] +
   \sqrt{\burnin} \Big \} \cdot \norm{\hpop(\theta_0) - \theta_0};
\end{multline*}
both with probability at least $1 - \delta$.
Finally, substituting the last two bounds to
the bound~\eqref{eqn:simplified-zt-anorm-bd},
and applying the conditions on
stepsize~\eqref{eq:stepsize-requirement-bootstrap-phase-i} and burn-in
period~\eqref{eq:burn-in-requirement-in-bootstrap-phase-i}, and using the fact
\mbox{$\anorm{\hpop(\theta_t) - \theta_t} \leq \anorm{\myZ_{t}} + \aUB \cdot
\norm{\myV_{t}}$} we have 
\begin{multline*}
\anorm{\hpop(\theta_t) - \theta_t} \leq \tfrac{c}{\sqrt{t}} \Big\{
\noisegausssemi + \noisevarsemi \sqrt{\log(\tfrac{1}{\delta})} \Big\}
+ c \tfrac{\aUB \LipCon}{ (1 - \contraction)} \Big \{ \constMGtype
\log(\tfrac{\numobs}{\delta}) \sqrt{\tfrac{\stepsize}{t}} +
\tfrac{1}{t \sqrt{ \stepsize}} \Big \} \; \Big \{ \noisegauss +
\noisevar \sqrt{\log( \tfrac{\numobs}{\delta})} \Big \} \\
+ \tfrac{c \aUB \Lip \boundstar}{1 - \contraction} \Big \{
\tfrac{\sqrt{\stepsize}}{t} + \tfrac{\stepsize}{\sqrt{t}} \Big \}
\DudGauss(\dualBall, \myrhotil) \DudExp(\dualBall, \myrhotil) \log^2
(\tfrac{\numobs}{\delta}) + \tfrac{\aUB \burnin}{t} \cdot
\norm{\hpop(\theta_0) - \theta_0}.
\end{multline*}
This completes the proof of
Corollary~\ref{cor:anorm-suboptimality-bound}.

\subsection{Proofs of key Lemmas for Theorem~\ref{new--thm:bellman-err-bound-general-nonlinear}}
In this section, we provide a detailed proofs of
Lemmas~\ref{lemma-vt-recursive-bound}
and~\ref{lemma-zt-martingale-bound-anorm} , which play a central role
in the proofs of
Theorem~\ref{new--thm:bellman-err-bound-general-nonlinear} and
Corollary~\ref{cor:anorm-suboptimality-bound}.

\subsubsection{Proof of Lemma~\ref{lemma-vt-recursive-bound}}
\label{subsubsec:proof-of-lemma-vt-recursive-bound}

We recursively expand the update rule for $\myV_t$ from
Algorithm~\ref{AlgRootSGD}, and obtain the identity:
\begin{align*}
t \myV_t &= (t - 1) \left(\myV_{t - 1} - \theta_{t - 1} + \theta_{t - 2} -
\Hstoch_t (\theta_{t - 2}) + \Hstoch_t (\theta_{t - 1}) \right) +
\Hstoch_t (\theta_{t - 1}) - \theta_{t - 1} \\
& = (1 - \stepsize) (t - 1) \myV_{t - 1} + (t - 1) \big( \Hstoch_t
(\theta_{t - 1}) - \Hstoch_t (\theta_{t - 2}) \big) + \big(\Hstoch_t
(\theta_{t - 1}) - \theta_{t - 1}\big)\\ &= (1 - \stepsize)^{\tau} (t
- \tau) \myV_{t - \tau} \\
& \qquad + \sum_{j = 1}^{\tau} (1 - \stepsize)^{j - 1} \Big[ (t - j)
  \big( \Hstoch_{t - j + 1} (\theta_{t - j}) - \Hstoch_{t - j + 1}
  (\theta_{t - j - 1}) \big) + \Hstoch_{t - j + 1} (\theta_{t - j}) -
  \theta_{t - j} \Big],
\end{align*}
where the positive integer $\tau$ will be chosen later.

Consequently, we have the bound
\begin{align*}
t \norm{\myV_t} &\leq (1 - \stepsize)^\tau (t - \tau) \norm{\myV_{t - \tau}}
+ \sum_{j = 1}^{\tau} (1 - \stepsize)^{j - 1} (t - j)
\norm{\hpop(\theta_{t - j}) - \hpop(\theta_{t - j - 1})}\\ & \qquad+
\Big \| \sum_{j = 1}^{\tau} (1 - \stepsize)^{j - 1} \Big( (t - j)
(\noise_{t - j + 1} (\theta_{t - j}) - \noise_{t - j + 1} (\theta_{t -
  j})) + \Hstoch_{t - j + 1} (\theta_{t - j}) - \theta_{t - j} \Big)
\Big \| \\
& \leq (1 - \stepsize)^\tau (t - \tau) \norm{\myV_{t - \tau}} + \sum_{j =
  1}^{\tau} (1 - \stepsize)^{j - 1} \Big( (t - j) \contraction
\stepsize \norm{\myV_{t - j}} + \norm{\hpop(\theta_{t - j}) - \theta_{t -
    j}} \Big)\\ &\qquad + \Big \| \sum_{j = 1}^{\tau} (1 -
\stepsize)^{j - 1} \Big( (t - j) (\noise_{t - j + 1} (\theta_{t - j})
- \noise_{t - j + 1} (\theta_{t - j - 1})) + \noise_{t - j + 1}
(\theta_{t - j}) \Big) \Big \|.
\end{align*}

The estimate sequence is $\ratev$ $\admissiblePar$-admissible for some
$\admissiblePar \in [0, 2]$, so that the map $t \mapsto t^2 \cdot
\ratev(t)$ is non-decreasing. Thus, on the event $\bigvEvent$ for a
burn-in $\burnin \geq \tfrac{12 \tau}{1 - \contraction}$, we
have the upper bound
\begin{align*}
(t - j) \norm{\myV_{t - j}} \leq (t - j) \ratev(t - j) \leq \tfrac{t^2}{t
    - j} \ratev(t) \leq \tfrac{1}{1 - \tau / \burnin} \cdot t
  \ratev(t) \leq \Big \{1 + \tfrac{1 - \contraction}{6} \Big \} \cdot
  t \ratev(t),
\end{align*}
valid for each integer $j \in [\tau]$.

Therefore, on the event $\bigthetaEvent \cap \bigvEvent$, we have the
bound
\begin{align}
\label{eq:vt-recursion-error-decomposition}  
 t \norm{\myV_t} & \leq \big(1 + \tfrac{1 - \contraction}{6} \big) \cdot
 \Big\{ (1 - \stepsize)^\tau + \contraction \stepsize\sum_{j = 1}^\tau
 (1 - \stepsize)^{j - 1} \Big\} \cdot t \ratev(t) + \sum_{j = 1}^\tau
 \ratetheta(t - j) + \termOne + \termTwo,
\end{align}
where \mbox{$\termTwo \mydefn \mybignorm{\sum_{j = 1}^\tau (1 -
    \stepsize)^{j - 1} \noise_{t - j + 1} (\theta_{t - j})}$,} and
\begin{align*}
\termOne & \mydefn \norm{\sum_{j = 1}^\tau (1 - \stepsize)^{j - 1} (t
  - j) \big( \noise_{t - j + 1} (\theta_{t - j}) - \noise_{t - j + 1}
  (\theta_{t - j - 1}) \big)}.
\end{align*}
We simplify the first two terms on the right-hand side of
bound~\eqref{eq:vt-recursion-error-decomposition} by appropriately
choosing the triple $(\tau, \stepsize, \burnin)$. The later two terms
$\termOne$ and $\termTwo$ are norms of zero-mean random vectors in
Banach spaces. First, we provide upper bound on these two noise terms.

\paragraph{Upper bound on $\termOne$:}

First, we observe that the sum consists of the
$(1-\stepsize)$-weighted differences \mbox{$(1 - \stepsize)^{j - 1} (t
  - j) \big( \noise_{t - j + 1} (\theta_{t - j}) - \noise_{t - j + 1}
  (\theta_{t - j - 1}) \big)$} that form a martingale difference
sequence with respect to the natural filtration $(\filtration_t)_{t
  \geq 0}$. On the event $\bigvEvent$, we have that
\begin{align*}
    \norm{(1 - \stepsize)^{j - 1} (t - j) \big( \noise_{t - j + 1}
      (\theta_{t - j}) - \noise_{t - j + 1} (\theta_{t - j - 1})
      \big)} \leq (t - j) \stepsize \Lip \ratev(t - j) \\ \leq
    \tfrac{t^2}{t - j} \stepsize \Lip \ratev(t) \leq 2 t \stepsize
    \Lip \ratev(t), \quad \mbox{a.s.}
\end{align*}
The last inequality is due to the non-decreasing property of the
function $t \mapsto t^2 \ratev(t)$ and the fact that $t \geq \burnin >
2 \tau$.

Since $\obssubset$ is symmetric and convex by assumption, the
difference \mbox{$\noise_{t - j + 1} (\theta_{t - j}) - \noise_{t - j
    + 1} (\theta_{t - j - 1})$} belongs to the set $2 \obssubset$.
Conditioning on the event $\bigvEvent$ and invoking
Lemma~\ref{lemma:concentration-in-banach-space-martingale} yields
\begin{align*}
\Big \| \tfrac{1}{\tau}\sum_{j = 1}^\tau (1 - \stepsize)^{j - 1} (t -
j) \big( \noise_{t - j + 1} (\theta_{t - j}) - \noise_{t - j + 1}
(\theta_{t - j - 1}) \big) \Big \| & \leq \tfrac{c t \stepsize \Lip
  \ratev(t)}{\sqrt{\tau}} \Big \{ \constMGtype +
\sqrt{\log(\tfrac{1}{\delta})} \Big \},
\end{align*}
with probability at least $1 - \delta$.

\paragraph{Upper bound on $\termTwo$:}
In order to bound the last term in the decomposition~\eqref{eq:vt-recursion-error-decomposition}, we decompose it into two parts:
\begin{align*}
    \sum_{j = 1}^\tau (1 - \stepsize)^{j - 1} \noise_{t - j + 1} (\theta_{t - j}) = \sum_{j = 1}^\tau (1 - \stepsize)^{j - 1} \noise_{t - j + 1} (\thetastar) +  \sum_{j = 1}^\tau (1 - \stepsize)^{j - 1} \big( \noise_{t - j + 1} (\theta_{t - j}) - \noise_{t - j + 1} (\thetastar) \big). 
\end{align*}
The former term is sum of independent random variables, while the
latter is a martingale. Note that by
Assumption~\ref{new--assume-pop-contractive}, we have $\norm{
  \noise_{t - j + 1} (\theta_{t - j}) - \noise_{t - j + 1}
  (\thetastar)} \leq \tfrac{\Lip \ratev(t - j + 1)}{1 - \contraction}$
on the event $\bigthetaEvent$. Invoking
Lemma~\ref{lemma:bernstein-in-banach-space-iid} yields
\begin{align*}
\| \tfrac{1}{\tau} \sum_{j = 1}^\tau (1 - \stepsize)^{j - 1} \noise_{t
  - j + 1} (\thetastar) \| \leq \tfrac{c}{\sqrt{\tau}} \Big \{
\noisegauss + \noisevar \sqrt{\log(1/\delta)} \Big \} + \tfrac{c
  \boundstar}{\tau} \Big \{ \log(1/\delta) + \DudExp(\dualBall,
\myrhotil) \Big \},
\end{align*}
with probability at least $1 - \delta$.

Using the Lipschitz assumption~\ref{new--assume-sample-operator} and the contraction assumption~\ref{new--assume-pop-contractive}, we
have
\begin{align*}
\norm{\noise_{t - j + 1}(\theta_{t - j}) - \noise_{t - j +
    1}(\thetastar)} \leq \Lip \norm{\theta_{t - j} - \thetastar} \leq
\tfrac{\Lip}{1 - \contraction} \norm{\hpop(\theta_{t - j}) - \theta_{t
    - j}}.
\end{align*}
Furthermore, since $\obssubset$ is symmetric and convex, we have that
$\noise_{t - j + 1} (\theta_{t - j}) - \noise_{t - j + 1} (\thetastar)
\in 2 \obssubset$.  Conditioning on the event $\bigthetaEvent$ and
invoking Lemma~\ref{lemma:concentration-in-banach-space-martingale}
yields
\begin{align*}
  \mybignorm{\tfrac{1}{\tau} \sum_{j = 1}^\tau (1 - \stepsize)^{j - 1}
    \big(\noise_{t - j + 1} (\theta_{t - j}) - \noise_{t - j + 1}
    (\thetastar) \big)} \leq \tfrac{c}{\sqrt{\tau}} \cdot \tfrac{\Lip
    \ratetheta(t - \tau + 1)}{1 - \contraction} \Big( \constMGtype +
  \sqrt{\log( \tfrac{1}{\delta})} \Big),
\end{align*}
with probability at least $1 - \delta$.

\paragraph{Combining the pieces:} Substituting
the above concentration bounds into the
decomposition~\eqref{eq:vt-recursion-error-decomposition} yields
the upper bound
\begin{multline*}
t \cdot \norm{\myV_t} \leq \Big \{ (1 - \stepsize)^{\tau} + \contraction
\stepsize \sum_{j = 1}^{\tau} (1 - \stepsize)^j \Big \} \cdot \Big \{
1 + \tfrac{1 - \contraction}{6} \Big \} \cdot t \ratev(t) + \tau \cdot
\ratetheta(t - \tau + 1) \\
+ c \sqrt{\tau} \cdot \tfrac{\Lip \ratetheta(t - \tau + 1)}{1 -
  \contraction} \Big \{ \constMGtype + \sqrt{\log(\tfrac{1}{\delta})}
\Big \} + c \sqrt{\tau} \Big \{ \noisegauss + \noisevar
\sqrt{\log(\tfrac{1}{\delta})} \Big \} \\
+ c \boundstar \Big \{ \log(\tfrac{1}{\delta}) + \DudExp(\dualBall, \myrhotil) \Big
\} + c \sqrt{\tau} \cdot t \stepsize \Lip \ratev(t) \cdot \Big \{
\constMGtype + \sqrt{\log(\tfrac{1}{\delta})} \Big \}.
\end{multline*}
Re-arranging the terms in the last bound yields
\begin{align*}
t \cdot \norm{\myV_t} & \leq \big \{ \contraction + (1 - \contraction)
\big( (1 - \stepsize)^\tau + \tfrac{1}{3} \big) + c \Lip \stepsize
\sqrt{\tau} \big(\constMGtype + \sqrt{\log(1 / \delta)} \big) \big \}
t \ratev(t)\\
& \qquad + \big \{ \tau + c \tfrac{\Lip \sqrt{\tau} }{1 -
  \contraction} \big(\constMGtype + \sqrt{\log(1 / \delta)} \big) \big
\} \ratetheta(t - \tau + 1)\\ &\qquad \qquad + c \sqrt{\tau} \Big(
\noisegauss + \noisevar \sqrt{\log(\tfrac{1}{\delta})} \Big) + c \boundstar \big \{
\log(\tfrac{1}{\delta}) + \DudExp(\dualBall, \myrhotil) \big \}.
\end{align*}

\subsubsection*{Case I: $t \geq \burnin + \lceil 2 \stepsize^{-1} \rceil$}
Taking $\tau = \lceil 2 \stepsize^{-1} \rceil \leq t - \burnin$ and
given a stepsize $\stepsize$ satisfying the bound
\begin{align}
\label{eq:stepsize-requirement}  
6 c \Lip \sqrt{\stepsize} \cdot \big(\constMGtype + \sqrt{\log(
  \tfrac{\numobs}{\delta})} \big) < 1 - \contraction,
\end{align}
we have the upper bounds
\begin{align*}
    & \tfrac{c \Lip}{1 - \contraction} \sqrt{\tau} \big(\constMGtype +
  \sqrt{\log( \tfrac{1}{\delta})} \big) \leq \tfrac{4}{\stepsize},
  \quad \mbox{and} \\
& \contraction + (1 - \contraction)\cdot \Big \{ (1 - \stepsize)^\tau
  + \tfrac{1}{3} + c \Lip \stepsize \sqrt{\tau} \big(\constMGtype +
  \sqrt{\log( \tfrac{1}{\delta})}\big) \Big \} \leq \tfrac{1 +
    \contraction}{2}.
\end{align*}
Furthermore, since the function $t \mapsto t^2 \cdot \ratetheta(t)$ is
non-decreasing, for burn-in period $\burnin \geq 4 \tau$, we have
\begin{align*}
    \ratetheta(t - \tau + 1) \leq \tfrac{t^2}{(t - \tau + 1)^2}
    \ratetheta(t) \leq \tfrac{16}{9} \ratetheta(t) \qquad \text{for
      all} \qquad t \geq \burnin.
\end{align*}
Substituting the bounds yields
\begin{align}
    t \norm{\myV_t} \leq \tfrac{1 + \contraction}{2} t \cdot \ratev(t)
    + \tfrac{8}{\stepsize} \ratetheta(t) + \tfrac{c}{\sqrt{\stepsize}}
    \Big( \noisegauss + \noisevar \sqrt{\log(\tfrac{1}{\delta})} \Big)
    + c \left( \log(\tfrac{1}{\delta}) + \DudExp(\dualBall, \myrhotil)
    \right),
\end{align}
which completes the proof of Lemma~\ref{lemma-vt-recursive-bound} in
the case of $t \geq \burnin + 2 / \stepsize$.

\subsubsection*{Case II: $\burnin \leq t \leq  \burnin + \lceil 2 \stepsize^{-1} \rceil $:}
This case requires a special treatment, since the number $\tau$ of
recursive expansion steps cannot be taken as large as $\lceil 2 /
\stepsize \rceil$. Instead, we choose $\tau = t - \burnin$, and expand
the recursions backwards up to the beginning of the iterates. In this
case, following the same arguments as above, on the event
$\bigthetaEvent \cap \bigvEvent$, the error
decomposition~\eqref{eq:vt-recursion-error-decomposition} takes the
form
\begin{align}
    t \cdot \norm{\myV_t} \leq (1 - \stepsize)^\tau \burnin
    \norm{\myV_{\burnin}} + \big(1 + \tfrac{1 - \contraction}{6} \big)
    \cdot \Big\{ \contraction \stepsize\sum_{j = 1}^\tau (1 -
    \stepsize)^{j - 1} \Big\} \cdot t \ratev(t) + \sum_{j = 1}^\tau
    \ratetheta(t - j) + \termOne + \termTwo.
\end{align}
Substituting the upper bounds on the terms $T_1$ and $T_2$ yields
\begin{multline*}
t \cdot \norm{\myV_t} \leq (1 - \stepsize)^\tau \burnin
\norm{\myV_{\burnin}} + \contraction \stepsize \sum_{j = 1}^{\tau} (1
- \stepsize)^j \cdot \Big \{ 1 + \tfrac{1 - \contraction}{6} \Big \}
\cdot t \ratev(t) + \tau \cdot \ratetheta(t - \tau + 1) \\
+ c \sqrt{\tau} \cdot \tfrac{\Lip \ratetheta(t - \tau + 1)}{1 -
  \contraction} \Big \{ \constMGtype + \sqrt{\log(\tfrac{1}{\delta})}
\Big \} + c \sqrt{\tau} \Big \{ \noisegauss + \noisevar
\sqrt{\log(\tfrac{1}{\delta})} \Big \} \\
+ c \boundstar \Big \{ \log(\tfrac{1}{\delta}) + \DudExp(\dualBall, \myrhotil) \Big
\} + c \sqrt{\tau} \cdot t \stepsize \Lip \ratev(t) \cdot \Big \{
\constMGtype + \sqrt{\log(\tfrac{1}{\delta})} \Big \}.
\end{multline*}
For a time index $t \in [\burnin, \burnin + 2 / \stepsize]$, we have
the decomposition
\begin{align*}
(1 - \stepsize)^\tau \burnin \cdot \norm{\myV_{\burnin}} &\leq \left(
  (1 - \stepsize)^\tau - 3(1 - \contraction) \right) \cdot q \Big\{ 1
  + \frac{2}{\stepsize \burnin}\Big\} t \cdot \ratev (t) + 3 (1 -
  \contraction) \burnin \cdot \norm{\myV_{\burnin}}\\
& \leq \Big \{ (1 - \stepsize )^\tau - 2 (1 - \contraction) \Big\}
  \cdot t \ratev (t) + 6 (1 - \contraction) \frac{\burnin^2}{t} \cdot
  \norm{\myV_{\burnin}}
\end{align*}
Given a stepsize $\stepsize$ satisfying the
requirement~\eqref{eq:vt-recursion-stepsize-requirement}, choosing the
number of steps such that \mbox{$\tau = t - \burnin \leq 2 /
  \stepsize$} leads to the inequalities
\begin{align*}
\Big\{ (1 - \stepsize )^\tau - 2 (1 - \contraction) \Big\} +
\contraction \stepsize \sum_{j = 1}^{\tau} (1 - \stepsize)^j \cdot
\Big \{ 1 + \tfrac{1 - \contraction}{6} \Big\} + c \sqrt{\tau} \cdot
\stepsize \Lip \cdot \Big \{ \constMGtype +
\sqrt{\log(\tfrac{1}{\delta})} \Big \} &\leq \frac{1 +
  \contraction}{2}, \quad \mbox{and} \\
\frac{c \Lip}{1 - \contraction} \sqrt{\tau} \big(\constMGtype +
\sqrt{\log(1 / \delta)} \big) \leq \frac{4}{\stepsize} + \tau &\leq
\frac{8}{\stepsize}.
\end{align*}
Putting together these bounds completes the proof in the second case.

\subsubsection{Proof of Lemma~\ref{lemma-zt-martingale-bound-anorm}}
\label{subsubsec:proof-lemma-zt-martingale-bound-anorm}

Expanding the update rule for $\myZ_t$ from Algorithm~\ref{AlgRootSGD}
we obtain the three-term decomposition \mbox{$t \cdot \myZ_t = \burnin
  \cdot \myZ_{\burnin} + \Martin_t + \Psi_t$,} where
\begin{align*}
\Martin_t \mydefn \sum_{s = \burnin}^{t - 1} \noise_s (\theta_{s -
  1}), \quad \mbox{and} \quad \Psi_t \mydefn \sum_{s = \burnin}^{t -
  1} (s - 1) \Big \{ \noise_s (\theta_{s - 1}) - \noise_s (\theta_{s -
  2}) \Big \}.
\end{align*}
It suffices to control each of these three terms in the
semi-norm induced by the set $\UpsSet$.

Beginning with the martingale $\{\Martin_t\}_{t \geq \burnin}$, we
further break it down into two parts:
\begin{align*}
 \Martin_t = \sum_{s = \burnin}^{t - 1} \varepsilon_s (\thetastar) +
 \sum_{s = \burnin}^{t - 1} \big( \varepsilon_s (\theta_{s - 1}) -
 \varepsilon_s (\thetastar) \big) \mydefn \MartinStar_t +
 \widetilde{M}_t.
\end{align*}
The term $\MartinStar_t$ is sum of $\mathrm{i.i.d.}$ random
variables. Invoking Lemma~\ref{lemma:bernstein-in-banach-space-iid}
and using the fact that the set $\UpsSet$ is contained within
$\aUB_{\UpsSet} \dualBall$, we have the bound
\begin{multline}
\label{eq:zt-mg-bound-part1-rearranged-anorm}  
\sup_{u \in \UpsSet} \inprod{u}{\MartinStar (t)} \leq c \sqrt{t} \Big
\{ \Exs \big[ \sup_{u \in \UpsSet} \inprod{u}{W} \big] +
\mybiggersqrt{\sup_{u \in \upsilon} \Exs \big[ \inprod{u}{W}^2 \big]
  \cdot \log(\tfrac{1}{\delta})}  \Big\} \\
+ c \aUB_{\UpsSet} \boundstar \Big \{ \DudExp(\dualBall, \myrhotil) +
\log(\tfrac{1}{\delta}) \Big \},
\end{multline}
where $W$ is the centered Gaussian process with covariance matching
that of $\noise(\thetastar)$.

Next we bound the terms $\inprod{u}{\widetilde{M}(t)}$ and
$\inprod{u}{\Psi (t)}$. First, we claim that conditioned on the event
$\bigthetaEvent \cap \bigvEvent$, we have
\begin{subequations}
\begin{align}
\label{eq:zt-mg-bound-part2-rearranged}    
\mynorm{\widetilde{M} (t)} & \leq c \frac{\Lip}{\MuCon}  \Big \{ \constMGtype +
\sqrt{\log(\tfrac{\numobs}{\delta})} \Big \} \cdot \mysqrt{\sum_{k =
    \burnin}^{t - 1} \ratetheta^2(k)}, \quad \mbox{and} \\
\label{eq:zt-mg-bound-part3-rearranged}
\mynorm{\Psi(t)} & \leq c \stepsize \Lip \Big \{ \constMGtype +
\sqrt{\log(\tfrac{\numobs}{\delta})} \Big \} \cdot \mysqrt{\sum_{s =
    \burnin}^t s^2 \ratev^2(s)},
\end{align}
\end{subequations}
both bounds holding with probability at least $1 - \delta$.

The proof of these two inequalities can be found at the end of this
subsection.  Since the set $\UpsSet$ is contained within
$\aUB_{\UpsSet} \dualBall$, it follows that
\begin{align*}
\sup_{u \in \UpsSet} \inprod{u}{\widetilde{M}(t)} \leq \aUB_{\UpsSet}
\mynorm{\widetilde{M}(t)}, \quad \mbox{and} \quad \sup_{u \in \UpsSet}
\inprod{u}{\Psi(t)} \leq \aUB_{\UpsSet} \norm{\Psi(t)}.
\end{align*}
Finally, observe that 
$\burnin \myZ_\burnin = \sum_{t = 1}^{\burnin}
   \noise_t(\thetastar) + \sum_{t = 1}^{\burnin}\{
   \noise_t(\theta_0)  - \noise_t(\thetastar)\}$. By
Lemma~\ref{lemma:bernstein-in-banach-space-iid} we have
\begin{multline*}
 \sup_{u \in \UpsSet} \inprod{u}{\sum_{t = 1}^{\burnin}
   \noise_t(\thetastar)} \leq c \sqrt{\burnin} \Big \{ \Exs \big[
   \sup_{u \in \UpsSet} \inprod{u}{W} \big] + \mybiggersqrt{\sup_{u
     \in \upsilon} \Exs \big[ \inprod{u}{W}^2 \big] \cdot
   \log(\tfrac{1}{\delta})} \Big \} \\ + c \aUB_{\UpsSet} \boundstar \Big \{
 \DudExp(\dualBall, \myrhotil) + \log(\tfrac{1}{\delta}) \Big \},
\end{multline*}
with probability at least $1 - \delta$.  On the other hand, using
Lemma~\ref{lemma:concentration-in-banach-space-martingale}, we have
\begin{multline*}
\sup_{u \in \UpsSet} \inprod{u}{\sum_{t = 1}^{\burnin} \big(
  \noise_t(\theta_0) - \noise_t(\thetastar) \big)} \leq \aUB_\UpsSet
\mynorm{\sum_{t = 1}^{\burnin} \big \{ \noise_t(\theta_0) -
  \noise_t(\thetastar) \big \} } \\
\leq c \Lip \aUB_\UpsSet \mynorm{\theta_0 - \thetastar} \cdot
\sqrt{\burnin} \Big \{ \constMGtype + \sqrt{\log(\tfrac{1}{\delta})}
\Big \} \\
\leq {c\aUB_\UpsSet \cdot \frac{\LipCon}{\MuCon}
  \sqrt{\burnin} \ratetheta(\burnin)} \Big \{ \constMGtype +
\sqrt{\log (\tfrac{1}{\delta})} \Big \},
\end{multline*}
with probability at least $1 - \delta$.  Combining the two bounds, we
conclude that
\begin{multline}
\label{eq:zt-mg-bound-burnin-anorm}  
\sup_{u \in \UpsSet} \inprod{u}{\myZ_\burnin} \leq \tfrac{c}{\sqrt{\burnin}}
\Big \{ \Exs \big[\sup_{u \in \UpsSet} \inprod{u}{W} \big] +
\mybiggersqrt{\sup_{u \in \upsilon} \Exs \big[ \inprod{u}{W}^2 \big]
  \cdot \log(\tfrac{1}{\delta})} \Big \} \\
+ \tfrac{c \aUB_{\UpsSet} \boundstar}{\burnin} \big( \DudExp (\dualBall,
\myrhotil) + \log(\tfrac{1}{\delta}) \big) + \tfrac{c \aUB_\UpsSet \cdot \Lip
  \ratetheta(\burnin)}{ \MuCon \sqrt{\burnin}} \big( \constMGtype +
\mysqrt{\log(\tfrac{1}{\delta})} \big),
\end{multline}
again with at least probability $1 - \delta$.

We now put together the
bounds~\eqref{eq:zt-mg-bound-part1-rearranged-anorm},~\eqref{eq:zt-mg-bound-part2-rearranged}~\eqref{eq:zt-mg-bound-part3-rearranged},
and~\eqref{eq:zt-mg-bound-burnin-anorm}.  By doing so, we are
guaranteed that conditioned on the event $\bigthetaEvent \cap
\bigvEvent$, for each integer $t \in [\burnin, \numobs]$, we have
\begin{multline*}
  \sup_{u \in \UpsSet} \inprod{u}{\myZ_t} \leq \tfrac{c}{\sqrt{t}} \Big
  \{ \Exs \big[ \sup_{u \in \UpsSet} \inprod{u}{W} \big] +
  \mybiggersqrt{\sup_{u \in \UpsSet} \Exs \big[ \inprod{u}{W}^2 \big]
    \log(\tfrac{1}{\delta})} \Big \} + \tfrac{\aUB_{\UpsSet}
    \boundstar}{t} \Big \{ \DudExp(\dualBall, \myrhotil) +
  \log(\tfrac{1}{\delta}) \Big \} \\
  + c \tfrac{\aUB_{\UpsSet} \cdot \stepsize \Lip}{\sqrt{t}} \Big \{
  \constMGtype + \sqrt{\log (\tfrac{1}{\delta})} \Big \}
  \mybiggersqrt{\frac{1}{t} \sum_{s = \burnin}^{t - 1} s^2
    \ratev^2(s)} \\
+ 2 c \tfrac{\aUB_{\UpsSet} \cdot \Lip}{\MuCon t} \Big \{ \constMGtype
+ \sqrt{\log(\tfrac{1}{\delta})} \Big \} \cdot \mybiggersqrt{\sum_{s =
    \burnin}^{t-1} \ratetheta^2(s) + \burnin \ratetheta^2 (\burnin) }.
\end{multline*}
The claim of Lemma~\ref{lemma-zt-martingale-bound-anorm} now follows by noting $\mybiggersqrt{\sum_{s =
    \burnin}^{t-1} \ratetheta^2(s) + \burnin \ratetheta^2 (\burnin) } = \mybiggersqrt{\sum_{s = 1
    \burnin}^{t-1} \ratetheta^2(s)}$.
It remains to prove
inequalities~\eqref{eq:zt-mg-bound-part2-rearranged}
and~\eqref{eq:zt-mg-bound-part3-rearranged}.

\paragraph{Proof of the bound~\eqref{eq:zt-mg-bound-part2-rearranged}:}

Conditioned on the event $\bigthetaEvent$, we have the upper bounds
\begin{align*}
\mynorm{\noise_s (\theta_{s - 1}) - \noise_s (\thetastar)} =
\mynorm{\Hstoch_s (\theta_{s - 1}) - \Hstoch_s (\thetastar) -
  \hpop(\theta_{s - 1}) + \hpop(\thetastar)} & \leq \Lip
\norm{\theta_{s - 1} - \thetastar} \\
& \leq \frac{\Lip}{\MuCon} \ratetheta(s - 1),
\end{align*}
where the last inequality follows from the assumption
\mbox{$\norm{\theta_{s - 1} -
  \thetastar} \leq \frac{1}{\MuCon}  \norm{\hpop(\theta_{s
    - 1}) - \hpop(\thetastar)}$} (cf. assumption~\eqref{EqnMucon}).

On the event $\bigthetaEvent$, we apply
Lemma~\ref{lemma:concentration-in-banach-space-martingale} to the
martingale differences \mbox{$ \{\noise_s(\theta_{s - 1}) - \noise_s
  (\thetastar) \}_{s=\burnin+1}^t$,} and find that
\begin{align*}
\mybiggernorm{\sum_{s = \burnin + 1}^t \big( \noise_s (\theta_{s - 1})
  - \noise_s (\thetastar) \big)} \leq c \Big \{ \constMGtype +
\sqrt{\log(\tfrac{1}{\delta})} \Big \} \; \frac{\Lip}{\MuCon}  \;
\mybiggersqrt{\sum_{s = \burnin}^{t-1} \ratetheta^2(s)},
\end{align*}
with probability at least $1 - \delta$, as claimed in
inequality~\eqref{eq:zt-mg-bound-part2-rearranged}.


\paragraph{Proof of bound~\eqref{eq:zt-mg-bound-part3-rearranged}:}

We now control the martingale sequence $\{\Psi_t\}_{t \geq \burnin}$.
Conditioned on the event $\bigvEvent$, we have
\begin{align*}
\norm{(s - 1) \big \{ \noise_s(\theta_{s - 1}) - \noise_s(\theta_{s -
    2}) \big \} } \leq (s-1) \Lip \cdot \norm{\theta_{s - 1} -
  \theta_{s - 2}} & = (s-1) \stepsize \Lip \norm{\myV_{s - 1}} \\
& \leq (s
- 1) \stepsize \Lip \ratev(s - 1),
\end{align*}
valid for any integer $s \in [\burnin, t]$.  By
Lemma~\ref{lemma:concentration-in-banach-space-martingale}, on the
event $\bigvEvent$, we have
\begin{align*}
\mybignorm{\sum_{s = \burnin + 1}^t (s - 1) \Big \{ \noise_s
  (\theta_{s - 1}) - \noise_s (\theta_{s - 2}) \Big \} } \leq c
\stepsize \Big \{ \constMGtype + \sqrt{\log (\tfrac{1}{\delta})} \Big
\} \Lip \mybiggersqrt{\sum_{s = \burnin}^{t - 1} s^2 \ratev^2(s)},
\end{align*}
with probability at least $1 - \delta$, which establishes the claim. This completes the proof of Lemma~\ref{lemma-zt-martingale-bound}.


\section{Discussion}

In this paper, we have analyzed \rootSANoSpace, a variance-reduced
stochastic approximation procedure designed for solving contractive
fixed-point equations in Banach spaces.  This procedure builds upon
the \ROOTSGD $\,$ algorithm~\cite{li2020root} for stochastic
optimization, as studied in past work by a subset of the current
authors.  Our main contribution was to derive non-asymptotic upper
bounds on the error of the \rootSA iterates in any semi-norm.  We
showed that these bounds are sharp in the sense that the leading order
term matches the optimal risk characterized by local minimax theory.
Furthermore, the sample complexity needed for such an
instance-dependent optimal statistical behavior scale with the
intrinsic complexity of the norm (measured in Dudley integral of the
dual ball under certain metric), instead of the problem dimension.
Our main results, while formulated for general Banach spaces and
contractions, have interesting consequences for specific classes of
problems.  Here we illustrated with applications to dynamical
programming and game theory, including stochastic shortest path
problems, minimax Markov games, as well as average-cost policy
evaluation.  In terms of proof techniques, our analysis is rather
different than much other theory on stochastic approximation that
relies the inner product available in the Hilbert setting.  To the
best of our knowledge, our paper is the first to provide sharp
non-asymptotic bounds for stochastic approximation \emph{without}
requiring such inner product structure.

\smallskip

\noindent Our work leaves open a number of open questions, among them:
\bcar
\item \textbf{Optimal sample complexity for SA schemes:} One open
  question in our analysis concerns how the minimal sample size scales
  with the contraction factor $\contraction \in (0,1)$, Our main
  results in this paper
  (Theorem~\ref{new--thm:bellman-err-bound-general-nonlinear}
  and~\ref{new--thm:main-nonlinear-general}) have a scaling condition
  of the form $\numobs \gtrsim (1 - \contraction)^{-4}$.  These
  results are novel even under this quadratic scaling, and also
  optimal for large $\numobs$.  However, it is not yet clear whether
  this quadratic scaling is necessary, or rather an artifact of our
  proof technique.  In certain special cases, the quadratic scaling
  can be avoided; for example, in the special case of $\hpop$ being an
  affine operator, our results (see
  Theorem~\ref{thm:bellman-err-bound-general-linear} and
  Corollary~\ref{thm:est-err-linear}) show that the quadratic scaling
  $\myorder \big( (1 - \contraction)^{-2} \big)$ is sufficient.
  Furthermore, in the classical Euclidean setting and in the special
  case gradient-update operator $\hpop: \theta \mapsto \theta -
  \beta^{-1} \nabla f (\theta)$ for a $\mu$-strongly-convex and
  $\beta$-smooth function $f$, the paper~\cite{li2020root} establishes
  instance-optimal bounds that require only $\myorder \big( (1 -
  \contraction)^{-2} \big)$ samples.  An interesting open problem,
  therefore, is to determine the minimum sample size for which
  non-asymptotic bounds of the form stated in this paper hold in the
  general Banach space setting.
  \item \textbf{Online statistical inference procedures:} In this
    paper, we focused exclusively on computing point estimates of the
    fixed point.  However, a natural question is the construction of
    confidence sets for the solution $\thetastar$ to the fixed-point
    equation. Ideally, such confidence set should be efficiently
    computable, asymptotically exact, while capturing the desirable
    non-asymptotic properties satisfied by our estimator. Focusing
    stochastic optimization in the Euclidean setting and
    Polyak-Ruppert-averaged SGD, the paper~\cite{chen2020statistical}
    proposed an online estimator for the covariance that partly
    achieves these goals. In a concurrent piece of work involving a
    subset of the authors~\cite{xia2022instance}, confidence sets and
    early stopping rules are developed in the special case of policy
    evaluation and optimization for discounted MDPs.  It is an
    interesting direction for future research to construct confidence
    sets with improved guarantees in the general setting, based purely
    on the algorithm's trajectory alone.
    \item \textbf{General operator equations beyond the contractive
      setting:} In the Euclidean setting, stochastic approximation
      procedures for nonlinear equations share geometric structure,
      giving rise to key concepts such as monotonicity and
      smoothness. This story becomes more complex for Banach spaces,
      with there being at least two distinct methods of analysis
      depending on the set-up.  On the one hand, if the operator
      $\hpop$ is mapping from the space $\vecspace$ to itself, then
      convergence is governed by contraction properties of the
      operator. On the other hand, if $\hpop$ maps from the Banach
      space $\vecspace$ to its dual space $\vecspace^*$, then a
      monotonicity condition with respect to the Bregman divergence
      plays a key role (see
      e.g.~\cite{juditsky2011solving,kotsalis2020simple}).  This paper
      focuses on the former case, in which $\hpop$ maps the Banach
      space to itself, but it is an interesting direction of future
      research to provide instance-dependent guarantees for various
      stochastic approximation procedures in the latter case, and
      examine their optimality properties. Even more broadly, it is
      interesting to consider stochastic approximation procedures for
      solving general non-linear equations defined on pairs of Banach
      spaces.  \ecar

\section*{Acknowledgements}
This work was supported in part by NSF-FODSI grant 202350, and the
DOD- ONR Office of Naval Research N00014-21-1-2842 as well as NSF-DMS
grant 2015454 and NSF-IIS grant 1909365 to MJW and PLB. We gratefully
acknowledge the support of the NSF through grants DMS-2023505 and of
the ONR through MURI award N000142112431 to PLB. This work was also
supported in part by the Vannevar Bush Faculty Fellowship program
under grant number N00014-21-1-2941 to MIJ. The authors would like to
thank Chris Junchi Li for helpful discussion.

\AtNextBibliography{\small}
\printbibliography


\appendix

\section{Restarting procedure}
\label{AppRestart}

In this section, we describe a simple restarting procedure that allows
us to refine the dependency of all of our bounds on the initial
condition.  This restarting procedure requires $\myorder(\burnin
\log\numobs)$ additional samples.

\begin{algorithm}
\caption{~\rootSA with re-starting}
\label{AlgRootRestart}
\begin{algorithmic}[1]
\STATE Given (a) Initialization $\thetainit \in \vecspace$, (b)
stepsize $\stepsize > 0$, (c) number of restarting epochs
$\numRestarts$
\STATE Set $\thetainit^{(0)} = \thetainit$ \FOR{$i = 1, \ldots,
  \numRestarts$} \STATE Run Algorithm~\ref{AlgRootSGD} with initial
point $\thetainit^{(i - 1)}$, stepsize $\stepsize$, burn-in period
$\burnin = \frac{c}{(1 - \contraction^2) \stepsize} \log
(\tfrac{\numobs}{\delta})$, and sample size $t_0 = 2 c \burnin$;
generate the sequence $(\theta_t^{(i - 1)})_{t = 0,1, \cdots, t_0}$.
\STATE Set $\thetainit^{(i)} = \theta_{t_0}^{(i - 1)}$.  \ENDFOR
\STATE Run Algorithm~\ref{AlgRootSGD} with initial point
$\thetainit^\circ = \thetainit^{(\numRestarts)}$, stepsize
$\stepsize$, burn-in period $\burnin = \frac{c}{(1 - \contraction^2)
  \stepsize} \log (\tfrac{\numobs}{\delta})$, and sample size $T =
\numobs - 2 c \numRestarts \burnin$; generate the sequence
$(\theta^\circ_t)_{t = 0,1, \cdots, T}$.  \RETURN $\theta_T^\circ$.
\end{algorithmic}
\end{algorithm}

Given some fixed number $\numRestarts \geq 1$ of restarting epochs, we
can run the \rootSA algorithm for $\numRestarts$ consecutive short
epochs, each with length $2 c \burnin$, with the constant $c$ being
the one in
equation~\eqref{new--eq:main-prop-bellman-err-bound-general-nonlinear-optimal}. The
last iterate $\theta_{2 c \burnin}$ of each short epoch is used as the
initial point of the subsequent epoch, and in the end, the output of
last short epoch is used as the initial point $\widetilde{\theta}_0$
to run a final single-epoch instantiation of \rootSA on the rest of
the data stream. The detail of the re-starting procedure is described
in Algorithm~\ref{AlgRootRestart}. In total, this restarting procedure
uses an additional $2 c \burnin \numRestarts$ samples, and the
initialization of the last epoch satisfies the bound
\begin{align}
\label{eq:restarting-initial-condition-nonlinear}  
    \norm{\widetilde{\theta}_0 - \hpop(\widetilde{\theta}_0) } \leq
    \tfrac{c}{\sqrt{\burnin}} \big( \noisegauss + \noisevar \sqrt{\log
      (1/\delta)} \big) + \tfrac{c \boundstar}{\burnin (1 -
      \contraction)} \big( \DudExp(\dualBall, \myrhotil) + \log
    (1/\delta) \big) + \tfrac{\norm{\thetainit - \hpop(\thetainit)}
    }{2^{\numRestarts}}.
\end{align}
By choosing $\numRestarts \geq \log \left(\norm{\hpop(\theta_0)
    - \theta_0}\sqrt{\numobs} / \noisegauss \right)$ with a restarting
sample size $2c\burnin\numRestarts$, we can ensure that $ \tfrac{
  \norm{\thetainit - \hpop(\thetainit)} }{2^{\numRestarts}} \leq \noisegauss / \sqrt{\numobs}$.

Our standard restarting procedure is based on the following
conditions.  We assume that the initialization $\theta_0$ is such that
the number of restarts $\numRestarts$ satisfies
\begin{subequations}
\label{eqn:initialization-burnin}
\begin{align}
\label{eqn:initialization}
\texttt{Initialization:} \qquad \log \left( \frac{\norm{\theta_0 -
    \hpop (\theta_0)} \sqrt{\numobs}}{\overline{\mathcal{W}}} \right) \leq c_0
\log \numobs,
\end{align}
for a universal constant $c_0 > 0$.  In words, the condition ensures
that the operator defect $\| \hpop(\theta_0) - \theta_0 \|$ for the
initialization $\theta_0$ is not exponentially large compared to
$\noisegauss$.  We set the number of restarts $\numRestarts$ as
\begin{align}
\label{eqn:numrestarts}
\texttt{Number of restarts:} \qquad \numRestarts = 2 c_0 \log \numobs
\qquad \qquad \qquad \qquad
\end{align}
\end{subequations}
These conditions ensure that performing $\numRestarts$ many restarts
requires at most
\begin{align*}
    2 c \burnin\log \left(\frac{\norm{\hpop(\theta_0) - \theta_0}
      \sqrt{\numobs}}{\overline{\mathcal{W}}} \right) \lesssim 4c_0
    \burnin \log(\numobs) \;\;\;\;
\end{align*}
additional samples, assuming that the original sample size is lower
bounded as $\numobs \gtrsim \tfrac{\LipCon^2 \constMGtype^2 }{(1 -
  \contraction)^4}$.  Substituting this bound back to the bounds from
Theorem~\ref{new--thm:bellman-err-bound-general-nonlinear} with the
optimal stepsize
choice~\eqref{eq:main-prop-bellman-err-bound-stepsize-choice}, we find
that
\begin{align}
\label{eq:main-thm-bellman-err-bound-general-nonlinear-restarting}  
 \norm{\hpop(\theta_\numobs) - \theta_\numobs} \leq
 \tfrac{c}{\sqrt{\numobs}} \Big \{\Exs[\norm{\GaussNoise} + \sqrt{
     \maxvar(\GaussNoise) \log(\tfrac{1}{\delta})} \Big\}
+ \tfrac{c \boundstar}{(1 - \contraction) \numobs} \cdot \big \{
\DudExp(\dualBall, \myrhotil) + \log(\tfrac{1}{\delta}) \big \}.
\end{align}


\section{Proofs of Theorem~\ref{new--thm:main-nonlinear-general}
and Corollary~\ref{Cor:main-nonlinear-general-anrom}}
\label{subsec:proof-thm-main-general-nonlinear}
In this section, we prove
Theorem~\ref{new--thm:main-nonlinear-general} and
Corollary~\ref{Cor:main-nonlinear-general-anrom}.  In fact,
Corollary~\ref{Cor:main-nonlinear-general-anrom} is actually a
generalization of Theorem~\ref{new--thm:main-nonlinear-general}; the
theorem follows from the corollary by setting $\anorm{\cdot} =
\norm{\cdot}$.  Accordingly, we devote our effort to proving the
corollary.


\subsection{Proof of Corollary~\ref{Cor:main-nonlinear-general-anrom}}
Define the pair
\begin{subequations}
\begin{multline}
\label{eq:copied-rate-theta-in-main-thm-proof}
\ratetheta^*(t) \mydefn \tfrac{c}{\sqrt{t}} \big( \noisegauss +
\noisevar \sqrt{\log(\tfrac{\numobs}{\delta})} \big) + \tfrac{c
  \boundstar}{1 - \contraction} \Big \{ \tfrac{1}{t} +
\tfrac{\stepsize \Lip }{\sqrt{t}} \cdot \constMGtype
\log(\tfrac{\numobs}{\delta}) \Big \} \; \Big \{\DudExp(\dualBall,
\myrhotil) + \log(\tfrac{1}{\delta}) \Big \} \\
+ c \tfrac{\burnin}{t} \norm{\thetainit - \hpop(\thetainit)}, \quad
\mbox{and}
\end{multline}
\begin{align}
\label{eq:copied-rate-v-in-main-thm-proof}    
  \ratev^*(t) \mydefn \tfrac{c}{1 - \contraction} \Big \{ \tfrac{1}{t
    \sqrt{\stepsize}} \Big( \noisegauss + \noisevar
  \sqrt{\log(\tfrac{\numobs}{\delta})} \Big) + \tfrac{\boundstar}{t}
  \Big \{ \log(\tfrac{\numobs}{\delta}) + \DudExp(\dualBall, \myrhotil)
  \Big \} \Big \} + 2 \big(\tfrac{\burnin}{t} \big)^2 \norm{\thetainit
    - \hpop(\thetainit)}.
\end{align}
\end{subequations}

Invoking Theorem~\ref{new--thm:bellman-err-bound-general-nonlinear}
and applying a union bound over the iterates, we have the pair of
bounds $\norm{\hpop(\theta_t) - \theta_t} \leq \ratetheta^*(t)$, and
and $\norm{\myV_t} \leq \ratev^*(t)$, uniformly for $t = \burnin,
\burnin + 1, \cdots, \numobs$ with probability at least $1 -
\delta$. Using the restarting scheme with parameter
choice~\eqref{eqn:numrestarts}, we can guarantee that the initial
operator defect $\norm{\hpop (\thetainit) - \thetainit}$ satisfies the
upper bound:
\begin{align*}
    \norm{\theta_0 - \hpop(\theta_0) } \leq
    \tfrac{c}{\sqrt{\burnin}} \big( \noisegauss + \noisevar \sqrt{\log
      (1/\delta)} \big) + \tfrac{c \boundstar}{\burnin (1 -
      \contraction)} \big( \DudExp(\dualBall, \myrhotil) + \log
    (1/\delta) \big).
\end{align*}

By the linearization condition~\ref{new--assume:linearization-anorm}, for
any $\theta \in \ball(\thetastar, \delrad_0)$, we have
\begin{align*}
  \delrad \mydefn \anorm{\theta - \thetastar} \; \leq \; \sup_{A \in
    \mathcal{A}_\delrad} \anorm{(I - A)^{-1} (\hpop(\theta) -
    \theta)}.
\end{align*}
In order to obtain an upper bound on~$\anorm{\theta_\numobs -
  \thetastar}$, it suffices to provide an bound for the quantity
$\sup_{A \in \mathcal{A}_\delrad} \anorm{(I - A)^{-1}
  (\hpop(\theta_{\numobs - 1}) - \theta_{\numobs - 1})}$ for any given
$\delrad > 0$. Recall that $h (\theta_{\numobs - 1}) - \theta_{\numobs
  - 1} = \myV_{\numobs} - \myZ_{\numobs}$, by definition, and in the
rest of this section we provide upper bounds on
$\anorm{\myV_{\numobs}}$ and $\anorm{\myZ_{\numobs}}$

\subsubsection*{Upper bound on $\anorm{\myV_{\numobs}}$}
Observe that $\anorm{(I - A)^{-1} \myV_{\numobs}} \leq \tfrac{\aUB}{1 - \contraction}
\norm{\myV_\numobs}$.  Thus, if we invoke the bound $\norm{\myV_\numobs}
\leq \ratev^*(t)$, where $\ratev^*$ is defined in
\eqref{eq:copied-rate-v-in-main-thm-proof}, we are guaranteed
that
\begin{multline*}
\anorm{(I - A)^{-1} \myV_{\numobs}} \leq \tfrac{c' \aUB}{(1 -
  \contraction)^2} \Big \{ \tfrac{1}{\numobs \sqrt{\stepsize}} \Big[
  \noisegauss + \noisevar \sqrt{\log(\tfrac{\numobs}{\delta})} \Big] +
\tfrac{\boundstar}{\numobs} \Big[ \log(\tfrac{\numobs}{\delta}) +
  \DudExp(\dualBall, \myrhotil) \Big] \Big \} \\
+ \tfrac{c \aUB}{1 - \contraction} \Big \{ \tfrac{\burnin}{\numobs}
\Big \}^2 \norm{\thetainit - \hpop(\thetainit)}.
\end{multline*}
with probability at least $1 - \delta$. 


\subsubsection*{Upper bound on $\anorm{\myZ_{\numobs}}$}

In order to establish a sharp upper bound on the term $\sup_{\Amat \in
  \mathcal{A}_\delrad} \norm{(I - \Amat)^{-1} \myZ_{\numobs + 1}}$, we
define the class of test functions $\UpsSet \mydefn \left\{ (I -
\Amat^*)^{-1} u:~ \Amat \in \mathcal{A}_\delrad, u \in \MySet
\right\}$. Substituting the the
bounds~\eqref{eq:copied-rate-theta-in-main-thm-proof}
and~\eqref{eq:copied-rate-v-in-main-thm-proof} in
Lemma~\ref{lemma-zt-martingale-bound-anorm} with we find that for any
given $\delrad > 0$, the quantity $\sup_{\Amat \in
  \mathcal{A}_\delrad} \anorm{(I - \Amat)^{-1} \myZ_\numobs}$ is upper
bounded as
\begin{multline}
\label{eq:est-err-thm-fixed-delta-bound}  
\tfrac{c}{\sqrt{\numobs}} \Big \{ \Exs \big[
  \sup_{\Amat \in \mathcal{A}_\delrad} \anorm{(I - \Amat)^{-1} W}
  \big] + \localvar(\delrad) \sqrt{\log(\tfrac{1}{\delta})} \Big \} +
\tfrac{c \boundstar \aUB}{\numobs (1 - \contraction)} \big \{ \DudExp(\dualBall, \myrhotil)
+ \log(\tfrac{1}{\delta}) \big \} \\
+ \tfrac{c \Lip \aUB}{(1 - \contraction)^2 \numobs} \big \{ \constMGtype +
\sqrt{\log(\tfrac{1}{\delta})} \big \} \Big \{ \stepsize
\mybiggersqrt{\sum_{s = \burnin}^{\numobs - 1} s^2 \ratev^{*2}(s)} +
\tfrac{1}{1 - \contraction} \mybiggersqrt{\sum_{s = 1}^{\numobs - 1}
  \ratetheta^{*2}(s)} \Big \},
\end{multline}
with probability at least $1 - \delta$.

\subsubsection*{Putting together the pieces}
The last two bounds are valid for a fixed value of $\delrad$. In order to derive the fixed-point condition in
Theorem~\ref{new--thm:main-nonlinear-general} and Corollary~\ref{Cor:main-nonlinear-general-anrom}, however, we need a
bound that holds uniformly over $\delrad$ in a suitable range, which
we now do.  Define the quantity
\begin{align*}
\overline{R}_\numobs \mydefn \tfrac{c}{(1 - \contraction) \sqrt{\numobs}} \Big \{ \Exs \big[
  \anorm{\GaussNoise} \big] + \mybiggersqrt{ \sup_{u \in \MySet}
  \Exs \big[ \inprod{u}{W}^2 \big] \log(\tfrac{1}{\delta})} \Big \} + \aUB \cdot \highorder_\numobs (\stepsize, \delta),
\end{align*}
and let $\underline{R} \mydefn \frac{1 - \contraction}{1 + \contraction} \overline{R}$.

It can be seen that the solutions to
equation~\eqref{eq:fixed-point-eq-for-anorm-bound} all belong to the
interval $[\underline{R}, \overline{R}]$. In particular, contraction
assumption~\ref{new--assume-pop-contractive}, we find that
\begin{align*}
    \frac{1}{1 + \contraction} \leq \opnorm{(I - \Amat)^{-1}} \leq
    \frac{1}{1 - \contraction}, \qquad \mbox{valid for any $\Amat \in
      \mathcal{A}_{\delrad}$,}
\end{align*}
which leads to the bounds $\tfrac{1}{1 + \contraction} \noisegausssemi
\leq \localgausssemi (\delrad) \leq \tfrac{1}{1 -
  \contraction} \noisegausssemi$ and $\tfrac{\noisevarsemi}{1 + \contraction} \leq
\localvarsemi (\delrad) \leq \tfrac{\noisevarsemi}{1 - \contraction} $
for any $\delrad > 0$.

By Theorem~\ref{new--thm:bellman-err-bound-general-nonlinear} and the
contractive assumption~\ref{new--assume-pop-contractive}, we have the
upper bound
\begin{align*}
\Prob \big[ \anorm{\theta_\numobs - \thetastar} \geq
  \overline{R}_\numobs \big] & \leq \delta.
\end{align*}
Consider the sequence $\delrad_\ell = 2^{\ell-1}
\underline{R}_\numobs$ for $\ell = 1, 2, \ldots, k$, where $k \mydefn
\log_2(\lceil {\overline{R}_\numobs}/{\underline{R}_\numobs}\rceil)$.
It forms a doubling grid $\mathcal{M}_\numobs \mydefn \{\delrad_1,
\delrad_2, \cdots, \delrad_k\}$ on the interval
$[\underline{R}_\numobs, \overline{R}_\numobs]$, and it can be seen
that $k$ satisfies the upper bound
\begin{align*}
 k & \leq \log (\tfrac{1 + \contraction}{1 - \contraction}) \leq 1 +
 \log(\tfrac{1}{1 - \contraction}).
\end{align*}

Taking a union bound over $\delrad \in \mathcal{M}_\numobs$, we find
that the bound~\eqref{eq:est-err-thm-fixed-delta-bound} holds with
probability at least $1 - k \delta$, uniformly over $\delrad \in
\mathcal{M}_{\numobs}$.  For any $\delrad \in \big[
  \underline{R}_\numobs, \overline{R}_\numobs \big]$, define the index
\mbox{$\ell(\delrad) \mydefn \max \{\ell \, \mid \, \delrad_\ell \leq
  \delrad\}$.} On the event above, we can conclude that
\begin{multline*}
\sup_{\Amat \in \mathcal{A}_\delrad} \anorm{(I - \Amat)^{-1}
  \myZ_\numobs} \leq \sup_{\Amat \in \mathcal{A}_{\delrad_{\ell (\delrad)
      + 1}}} \anorm{(I - \Amat)^{-1} \myZ_\numobs} \\
 \leq \tfrac{c}{\sqrt{\numobs}} \Big\{ \localgausssemi(2 \delrad) +
 \localvarsemi(2 \delrad) \sqrt{\log(\tfrac{1}{\delta}) + \log (\log
   \numobs) } \Big \} + \tfrac{c \aUB \boundstar}{(1 - \contraction) \numobs} \Big \{
 \DudExp(\dualBall, \myrhotil) + \log(\tfrac{1}{\delta}) + \log(\log
 \numobs) \Big \} \\
+ \tfrac{c \aUB \Lip}{(1 - \contraction)^2 \numobs} \Big \{ \constMGtype +
\sqrt{\log(\tfrac{1}{\delta})} \Big \} \; \Big \{ \stepsize
\mybiggersqrt{\sum_{s = \burnin}^{\numobs - 1} s^2 \ratev^{*2}(s)} +
\tfrac{1}{1 - \contraction} \mybiggersqrt{\sum_{s = 1}^{\numobs - 1}
  \ratetheta^{*2}(s)} \Big \},
\end{multline*}
for $\delrad \in[ \underline{R}_\numobs, \overline{R}_\numobs ]$.
Here we have used the facts that $\localgausssemi(\cdot)$ and $\localvarsemi(\cdot)$ are
non-decreasing functions. We now substitute our expressions for $\ratetheta^*$ and $\ratev^*$,
and conclude that conditioned on the event $\Event_\numobs^{(\theta)}
\cap \Event_\numobs^{(v)} \cap \left\{ \norm{\theta_\numobs -
  \thetastar} \leq \overline{R}_\numobs \right\}$, we have
\begin{align*} 
\delrad_\numobs & \leq \sup_{\Amat \in \mathcal{A}_{\delrad_\numobs}}
\anorm{(I - \Amat)^{-1} \myZ_\numobs} + \tfrac{\aUB}{1 - \contraction}
\ratev^*(\numobs) \\
& \leq \tfrac{c}{\sqrt{\numobs}} \left( \localgauss(2 \delrad_\numobs)
+ \localvar(2 \delrad_\numobs) \sqrt{\log(\tfrac{1}{\delta})} \right) +
\underline{R}_\numobs + \tfrac{c \aUB \boundstar}{(1 -
  \contraction)\numobs} \left( \DudExp(\dualBall, \myrhotil) +
\log(\tfrac{1}{\delta}) \right) \\
    & \qquad + \tfrac{c \aUB \Lip}{(1 - \contraction)^2
  \sqrt{\numobs}} \cdot \constMGtype \log(\tfrac{\numobs}{\delta})
\cdot \left( \sqrt{\stepsize} \noisegauss + \stepsize \boundstar \Big
\{ \DudExp(\dualBall, \myrhotil) + \log(\tfrac{\numobs}{\delta}) \Big \}
+ \sqrt{\tfrac{\burnin}{\numobs}} \norm{\thetainit - \hpop
  (\thetainit)} \right) \\
& \qquad \qquad + \tfrac{c' \aUB}{(1 - \contraction)^2} \Big \{
\tfrac{1}{\numobs \sqrt{\stepsize}} \noisegauss \sqrt{\log
  \tfrac{\numobs}{\delta}} + \tfrac{\boundstar}{\numobs} \Big[ \log
  \tfrac{\numobs}{\delta} + \DudExp(\dualBall, \myrhotil) \Big] \Big\} +
\tfrac{c \aUB}{1 - \contraction} \left( \tfrac{\burnin}{\numobs}
\right)^2 \norm{\thetainit - \hpop(\thetainit)}\\
& \leq \tfrac{c}{\sqrt{\numobs}} \Big \{ \localgausssemi (2
\delrad_\numobs) + \localvarsemi (2 \delrad_\numobs)
\sqrt{\log(\tfrac{1}{\delta})} \Big \} + \tfrac{c \aUB \log(\numobs /
  \delta)}{(1 - \contraction)^2} \Big \{ \constMGtype \Lip
\sqrt{\tfrac{\stepsize}{\numobs}} + \tfrac{1}{\numobs
  \sqrt{\stepsize}} \Big \} \cdot \noisegauss\\
& \qquad + \tfrac{c \aUB \boundstar \log(\tfrac{\numobs}{\delta})}{(1
  - \contraction)^2} \Big \{ \constMGtype \Lip
\tfrac{\stepsize}{\sqrt{\numobs}}+ \tfrac{1}{\numobs} \Big \} \cdot
\Big \{ \DudExp(\dualBall, \myrhotil) + \log(\tfrac{\numobs}{\delta})
\Big\} + \tfrac{\aUB \burnin}{ (1 - \contraction) \numobs} \cdot
\norm{\thetainit - \hpop(\thetainit)},
\end{align*}
with probability at least $1 -\delta$, valid for any $\delta \in (0,
1/k)$, where $k = 1 + \log \tfrac{1}{1 - \contraction}$.
  
Finally, noting that $\Prob \big[ \Event_\numobs^{(\theta)} \cap
  \Event_\numobs^{(v)} \cap \big \{ \norm{\theta_\numobs - \thetastar}
  \leq \overline{R}_\numobs \big \} \big] \geq 1 - \delta$, and using
the initialization conditions~\eqref{eqn:initialization-burnin}, we
obtain the bound that was claimed in
Corollary~\ref{Cor:main-nonlinear-general-anrom}.


\section{Proofs for multi-step contractions}

This section is devoted to the proofs of our results on multi-step
contractions, with~\Cref{thm:bellman-err-bound-general-linear} proved
in~\Cref{subsec:proof-theorem-thm-linear-compostep}
and~\Cref{thm:est-err-linear}
in~\Cref{subsec:proof-thm-est-err-linear}.


\subsection{Proof of Theorem~\ref{thm:bellman-err-bound-general-linear} }
\label{subsec:proof-theorem-thm-linear-compostep}

The proof of this theorem is similar to that of
Theorem~\ref{new--thm:bellman-err-bound-general-nonlinear}, but is
based on an improved version of Lemma~\ref{lemma-vt-recursive-bound},
stated as Lemma~\ref{lemma:vt-recursive-bound-compostep}.  At a high
level, there are three main steps:
\begin{itemize}
\item[1.] First, we use Lemma~\ref{lemma-zt-martingale-bound} and
  Lemma~\ref{lemma:vt-recursive-bound-compostep} to establish a
  relation between $\norm{\hpop(\theta_t) - \theta_t}$ and
  $\norm{\myV_t}$.
\item[2.] Second, starting with the coarse bound on
  $\norm{\hpop(\theta_t) - \theta_t}$ and $\norm{\myV_t}$ from
  Lemma~\ref{lemma:coarse-bound-gronwall}, we iteratively refine our
  bounds using the relation from Step 1.
\item[3.] Finally, we improve the higher-order terms in these bounds.
\end{itemize}

\subsubsection{Step 1: Relating $\norm{\hpop(\theta_t)-\theta_t}$ and $\norm{\myV_t}$}

We first state a sharpening of Lemma~\ref{lemma-vt-recursive-bound}
that holds for a multi-step contractive linear operator (see
Assumption~\ref{new--assume-compostep-contractive}).
\begin{lemma}
\label{lemma:vt-recursive-bound-compostep}
Under
assumptions~\ref{new--assume-compostep-contractive},~\ref{new--assume-noise-bound},
and~\ref{new--assume-sample-operator}, there exists a universal
constant $c > 0$ such that for stepsize $\stepsize$ satisfying the bound
\begin{subequations}
\begin{align}
  c \sqrt{\compostep \stepsize} \cdot \Lip \constMGtype
  \cdot \log \tfrac{\numobs}{\delta} \leq
  \tfrac{1}{3},\label{eq:vt-recursion-stepsize-requirement-compostep}
\end{align}
and the burn-in period $\burnin \geq \frac{c \compostep}{\stepsize} $,
given any $\admissiblePar$-admissible sequences $\ratetheta(t)$ and
$\ratev(t)$ with $0< \admissiblePar \leq 2$, on the event $\bigvEvent
\cap \bigthetaEvent$, the following bound holds uniformly with respect
to $t \in [ \burnin, \numobs]$, with probability $1 - \delta$:
\begin{multline}
\label{eq:vt-recursion-bound-compostep}  
\norm{\myV_t} \leq \tfrac{2\ratev(t)}{3} + \tfrac{c \compostep
  \ratetheta(t)}{t \stepsize} + \tfrac{c}{t} \sqrt{\tfrac{
    \compostep}{\stepsize} } \Big \{ \noisegauss + \noisevar
\sqrt{\log(\tfrac{\numobs}{\delta})} \Big \}+ \tfrac{c \boundstar}{t}
\big \{ \DudExp(\dualBall, \myrhotil) + \log (\tfrac{\numobs}{\delta})
\big \} + 4\big(\tfrac{\burnin}{t} \big)^2 \norm{\myV_{\burnin}}.
\end{multline}
\end{subequations}
\end{lemma}

\noindent See
Section~\ref{subsubsec:proof-lemma-vt-recursive-bound-compostep} for
the proof of this lemma.

\vspace{8pt}

In addition, by Lemma~\ref{lemma-zt-martingale-bound-anorm} and the
operator norm bound on $(I - \Amat)^{-1}$, conditioned on the event
$\bigthetaEvent \cap \bigvEvent$, we have the following bound
uniformly for $t \in [\burnin, \numobs]$,
\begin{multline}
\label{eq:zt-recursive-bound-compostep}  
\norm{\myZ_t} \leq \frac{c}{\sqrt{t}} \Big \{ \noisegauss + \noisevar
\sqrt{\log(\tfrac{\numobs}{\delta})} \Big \} + \frac{c \boundstar}{t} \Big
\{ \DudExp(\dualBall, \myrhotil) + \log(\tfrac{\numobs}{\delta}) \Big \} \\
+ \frac{c \Lip}{t} \Big \{ \constMGtype +
\sqrt{\log(\tfrac{\numobs}{\delta})} \Big \} \; \Big \{ \stepsize
\mybiggersqrt{\sum_{s = \burnin}^{t - 1} s^2 \ratev^2(s)} +
\compostep \mybiggersqrt{\sum_{s = 1}^{t - 1} \ratetheta^2(s)}
\Big \},
\end{multline}
with probability at least $1 - \delta$.
   
\subsubsection{Step 2: Bounds using bootstrapping}
   
Akin to the proof of
Theorem~\ref{new--thm:bellman-err-bound-general-nonlinear}, we impose
the restrictions that the estimate sequences $(\ratetheta, \ratev)$
are $\tfrac{1}{2}$- and $1$-admissible, respectively.

Consider a new pair $(\ratev^+, \ratetheta^+)$ satisfying the initial
bounds \mbox{$\ratev^+(\burnin) \geq \norm{\myV_{\burnin}}$} and
\mbox{$\ratetheta^+(\burnin) \geq \norm{\hpop(\theta_0) - \theta_0}$,}
and such that
\begin{subequations}
\label{eq:bootstrap-iteration-compostep-phase-i}
\begin{align}
\label{eq:bootstrap-iteration-compostep-phase-i-ratev}  
\ratev^+(t) & \geq \tfrac{2}{3} \ratev(t) + \tfrac{1}{t
  \sqrt{\stepsize}} \cdot \frac{c \compostep}{\sqrt{t \stepsize}}
\cdot \sqrt{t}\ratetheta(t) \nonumber \\ & \qquad+ \frac{c}{t}
\sqrt{\tfrac{ \compostep }{\stepsize} } \Big \{ \noisegauss +
\noisevar \sqrt{\log(\tfrac{\numobs}{\delta})} \Big \} + c \frac{
  \boundstar}{t} \, \Big \{ \DudExp(\dualBall, \myrhotil) +
\log(\tfrac{\numobs}{\delta}) \Big \} + 4\big(\tfrac{\burnin}{t}
\big)^2 \norm{\myV_{\burnin}}, \quad \mbox{and} \\
\label{eq:bootstrap-iteration-compostep-phase-i-ratetheta}
\ratetheta^+ (t) & \geq \frac{c}{\sqrt{t}} \Big \{
\frac{1}{\sqrt{\stepsize t}} + {\sqrt{\stepsize} \Lip \constMGtype}
\log \tfrac{\numobs}{\delta} \Big \} \cdot \Big \{ t
\sqrt{\stepsize}\ratev(t) \Big \} + 2 c \tfrac{ \Lip
  \compostep}{\sqrt{t}} \constMGtype \log{\tfrac{\numobs}{\delta}}
\cdot \ratetheta(t), \nonumber \\
& \qquad + \tfrac{c}{\sqrt{t}} \Big \{ \noisegauss + \noisevar
\sqrt{\log(\tfrac{\numobs}{\delta})} \Big \} + \tfrac{c \boundstar}{t}
\Big \{ \DudExp(\dualBall, \myrhotil) + \log(\tfrac{\numobs}{\delta})
\Big \}.
\end{align}
\end{subequations}
for each integer $t \in [\burnin, \numobs]$.

By combining the bounds~\eqref{eq:vt-recursion-bound-compostep}
and~\eqref{eq:zt-recursive-bound-compostep}, we are guaranteed that
\begin{align*}
\Prob \Big[ \Event^{(\theta)}_\numobs (\ratetheta^+) \cap
  \Event^{(v)}_\numobs (\ratev^+) \Big] \geq \Prob \Big[
  \bigthetaEvent \cap \bigvEvent \Big] - \delta,
\end{align*}

Our goal is to construct two series of admissible sequences $\big(
\ratev^{(i)}, \ratev^{(i)} \big)$ with $i = 0,1, \cdots$, such that
the pair $(\norm{\myV_t}, \norm{\hpop (\theta_t) - \theta_t})_{t \geq
  \burnin}$ are dominated by $\big( \ratev^{(i)} (t), \ratev^{(i)} (t)
\big)_{t \geq 0}$, with high probability. Concretely, we consider
sequences of a particular form $\ratev^{(i)} (t) =
\frac{\psi_v^{(i)}}{t \sqrt{\stepsize}}$ and $\ratetheta^{(i)} (t) =
\frac{\psi_\theta^{(i)}}{\sqrt{t}}$, for pairs of positive reals
$\big( \psi_v^{(i)}, \psi_\theta^{(i)} \big)$ independent of
$t$. Apparently, with such forms, the sequence $\ratetheta^{(i)}$ is
$\frac{1}{2}$-admissible, and the sequence $\ratev^{(i)}$ is
$1$-admissible. However, if we directly substitute the sequences
$\big( \ratev^{(i)} (t), \ratev^{(i)} (t) \big)$ of such forms into
the iteration~\eqref{eq:bootstrap-iteration-compostep-phase-i}, the
resulting sequences $(\ratetheta^+, \ratev^+)$ will no longer be of
the desired form. So in order to unify the coefficients in
equation~\eqref{eq:bootstrap-iteration-compostep-phase-i} into the
same time scale, given a stepsize $\stepsize > 0$, we define the
burn-in time
\begin{subequations}
\begin{align}
\label{eq:burn-in-requirement-compostep-in-bootstrap-phase-i}  
\burnin = \tfrac{c  \compostep}{\stepsize} \log
(\tfrac{\numobs}{\delta}).
\end{align}
For each $t = \burnin, \burnin +1 \ldots$, the coefficients in
\eqref{eq:bootstrap-iteration-compostep-phase-i} then satisfy the
bounds
\begin{align}
\label{eq:coefficient-bounds-compostep-in-bootstrap-phase-i}  
\tfrac{c \compostep}{\sqrt{\stepsize  t}} \leq
\tfrac{1}{6} \sqrt{\compostep}, \quad
\tfrac{c}{\sqrt{\stepsize t}} \leq \tfrac{1}{12\sqrt{\compostep}}, \quad \mbox{and} \quad \tfrac{
  \compostep}{\sqrt{t}} \log (\tfrac{\numobs}{\delta}) \leq
\sqrt{\stepsize  \compostep}.
\end{align}
\end{subequations}
Therefore, if we construct a two-dimensional vector sequence
$\psi^{(i)} = \begin{bmatrix} \psi_v^{(i)} & \psi_\theta^{(i)}
\end{bmatrix}^T$ satisfying the recursive relation $\psi^{(i+1)} = \Qmat \psi^{(i)} + \bvec$, where
\begin{align}
\Qmat &\mydefn \begin{bmatrix} 2/3 & \frac{\sqrt{\compostep}}{6}
  \\ \tfrac{1}{12 \sqrt{\compostep}} + c \Lip \constMGtype
  \sqrt{\stepsize} \cdot \log (\tfrac{\numobs}{\delta}) & 2 c \Lip
  \constMGtype \sqrt{\stepsize \compostep} \log
  (\tfrac{\numobs}{\delta})
    \end{bmatrix}, \quad \mbox{and} \notag \\
    \bvec &\mydefn c \cdot \begin{bmatrix} \sqrt{\compostep} \Big(
      \noisegauss + \noisevar \sqrt{\log(\numobs / \delta)} \Big) +
      \boundstar \sqrt{\stepsize} \left( \log(\tfrac{1}{\delta}) +
      \DudExp(\dualBall, \myrhotil) \right) + \burnin \sqrt{\stepsize}
      \norm{\myV_{\burnin}}\\ \Big( \noisegauss + \noisevar
      \sqrt{\log(\numobs / \delta)} \Big) + \boundstar
      \sqrt{\tfrac{\stepsize}{\compostep}} \left(\log(\numobs /
      \delta)+ \DudExp(\dualBall, \myrhotil) \right) + \sqrt{\burnin}
      \norm{\hpop(\theta_0) - \theta_0}
    \end{bmatrix},
    \label{eqn:boot-recursion-Pars-multistep}
\end{align}
they will satisfy the
requirement~\eqref{eq:bootstrap-iteration-compostep-phase-i}, leading
to the probability bound:
\begin{align}
\label{eqn:multisep-prob-recursion}
  \Prob \Big[ \Event^{(\theta)}_\numobs \big( \ratetheta^{(i + 1)}
    \big) \cap \Event^{(v)}_\numobs \big( \ratev^{(i + 1)} \big)\Big]
  \geq \Prob \Big[ \Event^{(\theta)}_\numobs \big( \ratetheta^{(i)}
    \big) \cap \Event^{(v)}_\numobs \big( \ratev^{(i)} \big)\Big] -
  \delta,
\end{align}
for the sequences $\ratetheta^{(i)} (t) = \psi_\theta^{(i)} /
\sqrt{t}$ and $\ratev^{(i)} (t) = \psi_v^{(i)} / (\sqrt{\stepsize}
t)$.

It remains to specify an initial condition for the recursion
above. Note that Lemma~\ref{lemma:coarse-bound-gronwall} implies that
we have
\begin{align*}
\norm{\theta_t - \thetastar} + \norm{\myV_t} \leq e^{1 + \Lip \stepsize
  t} \big( \boundstar + \norm{\theta_0 - \thetastar} \big)
\end{align*}
almost surely. So we can take the initialization:
\begin{align*}
\psi_v^{(0)} \mydefn \numobs \sqrt{\stepsize} e^{1 + \Lip \stepsize \numobs}
(\boundstar + \norm{\theta_0 - \thetastar}), \quad \mbox{and} \quad
\psi_\theta^{(0)} \mydefn \sqrt{\numobs} e^{1 + \Lip \stepsize \numobs}
(\boundstar + \norm{\theta_0 - \thetastar}),
\end{align*}
for which the bounds $\norm{\myV_t} \leq \tfrac{\psi_v^{(0)}}{t
  \sqrt{\stepsize}}$ and $\norm{\theta_t - \hpop (\theta_t)} \leq
\tfrac{\psi_\theta^{(0)}}{\sqrt{t}}$ hold almost surely.

Given such an initial condition and the
recursion~\eqref{eqn:multisep-prob-recursion}, we find that
\begin{align*}
\Prob \Big[ \Event_\numobs^{(v)}(\ratev^{(i)}) \cap
  \Event_\numobs^{(\theta)}(\ratetheta^{(i)}) \Big] \geq \Prob
\Big[\Event_\numobs^{(v)}(\ratev^{(0)}) \cap
  \Event_\numobs^{(\theta)}(\ratetheta^{(0)}) \Big] - i \delta = 1 - i
\delta.
\end{align*}
It remains to understand the behavior of $\psi^{(i)}$
for large values of the index $i$,
i.e. the after $i$ iterations of the bootstrapping argument.
We do so by solving the recursion
  $\psi^{(i+1)} = \Qmat \psi^{(i)} + \bvec$.
Let us define a new matrix
\begin{align*}
\QmatTil & \defn
\begin{bmatrix} 2/3 &
  \tfrac{\sqrt{\compostep}}{6}\\ \tfrac{1}{6 \sqrt{\compostep}} & 2/3
\end{bmatrix} \; \stackrel{(i)}{=}
  \begin{bmatrix} \sqrt{\compostep} &
  \sqrt{\compostep}\\ 1 & -1
    \end{bmatrix} \cdot \begin{bmatrix}
     \tfrac{5}{6} &0\\0&\tfrac{1}{2}
    \end{bmatrix} \cdot  \begin{bmatrix}
     \sqrt{\compostep} & \sqrt{\compostep}\\ 1 & -1
    \end{bmatrix}^{-1},
\end{align*}
where the equivalence (i) follows by a direct calculation.  Note that the stepsize condition~\eqref{eqn:stepsize-multistep} ensures that 
\begin{align}
\label{eq:stepsize-requirement-bootstrap-compostep-phase-i}  
c \Lip \constMGtype \log \tfrac{\numobs}{\delta} \cdot
\sqrt{\compostep  \stepsize} \leq \tfrac{1}{12},
\end{align}
then the matrix $\QmatTil$ is coordinate-wise larger than the matrix
$\Qmat$ from equation~\eqref{eqn:boot-recursion-Pars-multistep}, and
consequently we are guaranteed that $\Qmat u \preceqOrthant \QmatTil
u$ for any $2$-dimensional vector $u \succeqOrthant 0$.  Thus, for
each integer $N = 1, 2, \ldots$, we have the upper bounds
\begin{align*}
\begin{bmatrix}
\psi_v^{(N)} \\ \psi_\theta^{(N)}
    \end{bmatrix}  = \left(\sum_{i = 0}^{N - 1} \Qmat^i \right) b + \Qmat^N \begin{bmatrix}
     \psi_v^{(0)} \\ \psi_\theta^{(0)}
    \end{bmatrix} &  \preceqOrthant \left(\sum_{i = 0}^{N - 1} \QmatTil^i \right) b_\psi + \QmatTil^N \begin{bmatrix}
     \psi_v^{(0)} \\ \psi_\theta^{(0)}
\end{bmatrix} \\
& \preceqOrthant (I - \QmatTil)^{-1} b + e^{-  N/6}
\sqrt{\compostep} \big( \psi_v^{(0)} + \psi_\theta^{(0)} \big)
\bm{1}_2.
\end{align*}
By taking $N = c \Lip \numobs \log \numobs$, replacing $\delta$ with
$\delta / N$ and substituting with the above bounds, we find that
\begin{subequations}
\label{eq:bootstrap-phase-i-consequence-compostep}
\begin{multline}
 \label{eq:bootstrap-phase-i-consequence-compostep-v}   
t \sqrt{\stepsize} \cdot \norm{\myV_t} \leq \psi_v^{(N)} \leq c 
\sqrt{\compostep} \big \{ \noisegauss + \noisevar \sqrt{\log
  (\tfrac{\numobs}{\delta})} \big \} \\ + {c  \boundstar
  \sqrt{\stepsize}} \left( \log (\tfrac{\numobs}{\delta}) +
\DudExp(\dualBall, \myrhotil) \right) + c \burnin \sqrt{\stepsize}
\norm{\myV_{\burnin}} + \sqrt{\burnin \compostep} \norm{\hpop(\theta_0) -
  \theta_0},
\end{multline}
along with
\begin{multline}
\label{eq:bootstrap-phase-i-consequence-compostep-theta}  
\sqrt{t} \norm{\hpop(\theta_t) - \theta_t} \leq \psi_\theta^{(N)} \leq
c \big \{ \noisegauss + \noisevar \sqrt{\log
  (\tfrac{\numobs}{\delta})} \big \} + c \boundstar \sqrt{\stepsize}
\big \{ \log (\tfrac{\numobs}{\delta}) + \DudExp(\dualBall, \myrhotil)
\big \} \\
+ c \burnin \sqrt{\stepsize / \compostep} \norm{\myV_{\burnin}} +
\sqrt{\burnin} \norm{\hpop(\theta_0) - \theta_0},
\end{multline}
\end{subequations}
valid uniformly over $t \in \{ \burnin, \burnin + 1, \cdots, \numobs
\}$ with probability at least $1 - \delta$.

The latter bound, when combined with the
Lemma~\ref{lemma:burnin-period-bound} yields an upper bound on
$\norm{\hpop(\theta_t) - \theta_t}$ which has the correct leading-order term, i.e., the correct dependence on the term $\noisegauss +
\noisevar \sqrt{\log \tfrac{\numobs}{\delta}}$.  In order to refine
the dependence on the terms $\norm{\hpop(\theta_0) - \theta_0}$ and
$\log (\tfrac{\numobs}{\delta}) + \DudExp(\dualBall, \myrhotil)$, we
need do another round of bootstrapping.


\subsubsection{Step 3: Improving the higher-order terms}

With a slight abuse of notation, let the $2$-vector $\psi^{(N)}
\mydefn (\psi_v^{(N)}, \psi_\theta^{(N)})$ be defined by the
right-hand side of
equation~\eqref{eq:bootstrap-phase-i-consequence-compostep}, and
consider the choices \mbox{$\ratev(t) \mydefn \tfrac{\psi_v^{(N)}}{t
    \sqrt{\stepsize}}$} and \mbox{$\ratetheta(t) \mydefn
  \tfrac{\psi_\theta^{(N)}}{\sqrt{t}}$.}  Conditioned on the event
$\bigthetaEvent \cap \bigvEvent$, we have
\begin{align*}
\norm{\hpop(\theta_t) - \theta_t} & \leq \tfrac{c}{\sqrt{t}} \Big \{
\noisegauss + \noisevar \sqrt{\log(\tfrac{1}{\delta})} \Big \} +
\tfrac{c \boundstar}{t} \Big \{ \DudExp(\dualBall, \myrhotil) +
\log(\tfrac{1}{\delta}) \Big \} \nonumber \\
& \qquad + \Big \{ \tfrac{1}{t} + c \tfrac{\stepsize \Lip
  \constMGtype}{ \sqrt{t}} \cdot \log (\tfrac{\numobs}{\delta}) \Big
\} \tfrac{ \psi_v^{(N)}}{\sqrt{\stepsize}} + 2c \tfrac{ 
  \compostep\Lip \constMGtype}{t} \log(\tfrac{n}{\delta}) \cdot
\psi_\theta^{(N)} \\
& \leq \Big \{ \tfrac{c}{\sqrt{t}} + 
\sqrt{\tfrac{\compostep}{\stepsize}} \big[ \tfrac{1}{t} + c
  \tfrac{\stepsize \Lip \constMGtype}{ \sqrt{t}} \cdot \log
  (\tfrac{\numobs}{\delta}) \big] \Big \} \Big \{ \noisegauss +
\noisevar \sqrt{\log(\tfrac{\numobs}{\delta})} \Big \} \\
& \qquad+ c' \boundstar \Big \{ \tfrac{1}{t} + \tfrac{\stepsize \Lip
  \constMGtype}{ \sqrt{t}} \cdot \log(\tfrac{\numobs}{\delta}) \Big \}
\; \Big \{ \DudExp(\dualBall, \myrhotil) + \log(\tfrac{\numobs}{\delta})
\Big \} \\
& \qquad + c' \Big \{ \tfrac{1}{t}
\sqrt{\tfrac{\compostep}{\stepsize}} + \tfrac{\sqrt{\stepsize
    \compostep} \Lip \constMGtype}{ \sqrt{t}} \cdot \log
(\tfrac{\numobs}{\delta}) + \tfrac{ \compostep \Lip
  \constMGtype}{t} \log({\tfrac{\numobs}{\delta}}) \Big \}
\sqrt{\burnin} \norm{\hpop(\theta_0) - \theta_0} \\
& \qquad + c' \Big \{ \tfrac{1}{t} + \tfrac{\stepsize \Lip
  \constMGtype}{ \sqrt{t}} \cdot \log(\tfrac{\numobs}{\delta}) +
\tfrac{ \sqrt{\stepsize / \compostep} \Lip \constMGtype}{t}
\log(\tfrac{\numobs}{\delta})\Big \} \burnin \norm{\myV_{\burnin}},
\end{align*}
with probability at least $1 - \delta$.

Given a burn-in period $\burnin$ satisfying
\eqref{eq:burn-in-requirement-compostep-in-bootstrap-phase-i} and step
size satisfying
\eqref{eq:stepsize-requirement-bootstrap-compostep-phase-i}, using the
bound on $\norm{\myV_{\burnin}}$ from
Lemma~\ref{lemma:burnin-period-bound}, we have the upper bound
$\norm{\hpop(\theta_t) - \theta_t} \leq \widetilde{\ratetheta}(t)$,
with probability at least $1 - \delta$, uniformly over all integers $t
\in [\numobs]$, where
\begin{multline*}
\widetilde{\ratetheta}(t) \mydefn \tfrac{c_1}{\sqrt{t}} \Big \{
\noisegauss + \noisevar \sqrt{\log(\tfrac{\numobs}{\delta})} \Big \} +
            {c_2 \boundstar} \Big \{ \tfrac{1}{t} + \tfrac{\stepsize
              \Lip \constMGtype}{\sqrt{t}} \cdot \log^3
            (\tfrac{\numobs}{\delta}) \Big \} \; \Big \{
            \DudExp(\dualBall, \myrhotil) +
            \log(\tfrac{\numobs}{\delta}) \Big \} \\
+ c_2 \Big \{ \tfrac{\stepsize \burnin \Lip \constMGtype}{ \sqrt{t}}
\cdot \log(\tfrac{\numobs}{\delta}) + \tfrac{\burnin}{t} \Big \} \;
\norm{\hpop(\theta_0) - \theta_0}.
\end{multline*}

By substituting the upper bound $\widetilde{\ratetheta}$ into
equation~\eqref{eq:bootstrap-iteration-compostep-phase-i-ratev}, we
obtain a recursive inequality that takes as input an admissible
sequence $\ratev(t)$, and generates as output a new sequence
$\ratev^+(t)$ such that
\begin{align*}
\Prob \Big[ \Event_\numobs^{(v)} (\ratev^+) \Big] & \geq \Prob \Big[
    \Event_\numobs^{(v)} (\ratev) \Big] - \delta.
\end{align*}
Taking any integer $N_1 > 0$, by applying the recursive inequality for
$N_1$ times with $\delta' = \delta / N_1$, we get a sharper bound for
$\norm{\myV_t}$ with probability $1 - \delta$:
\begin{align*}
    \norm{\myV_t} &\leq 3 c \left[ \tfrac{ \sqrt{\compostep}}{t
        \sqrt{\stepsize}} \Big( \noisegauss + \noisevar
      \sqrt{\log(\tfrac{\numobs N_1}{\delta})} \Big) +
      \tfrac{\boundstar}{t} \left( \log(\tfrac{\numobs N_1}{\delta}) +
      \DudExp(\dualBall, \myrhotil) \right) \right]\\ &\qquad +\tfrac{c
      \compostep }{ \stepsize t} \widetilde{\ratetheta} (t) + c
    \left(\tfrac{\burnin}{t} \right)^2 \norm{\myV_{\burnin}} +
    \left(\tfrac{1 + \contraction}{2} \right)^{N_1} \cdot
    \tfrac{\psi_v^{(N)}}{t \sqrt{\stepsize}}.
\end{align*}

Taking $N_1 \mydefn 10 \log \numobs$, for stepsize and burn-in period
satisfying the
conditions~\eqref{eq:burn-in-requirement-compostep-in-bootstrap-phase-i}
and~\eqref{eq:stepsize-requirement-bootstrap-compostep-phase-i}, some
algebra yields
that
$\norm{\myV_t} \leq \widetilde{\ratev}(t)$
with probability at least $1 - \delta$,
uniformly for each integer $t \in [\burnin, \numobs]$, where
\begin{align}
 \label{eq:vt-rate-final-compostep}  
\widetilde{\ratev}(t) \mydefn c' \Big \{ \tfrac{1}{t}
\sqrt{\tfrac{\compostep}{\stepsize}} \big[ \noisegauss + \noisevar
  \sqrt{\log(\tfrac{\numobs}{\delta})} \big] + \tfrac{\boundstar}{t}
\big[ \log(\tfrac{\numobs}{\delta}) + \DudExp(\dualBall, \myrhotil)
  \big] \Big \} + 2 c' \big(\tfrac{\burnin}{t} \big)^2
\norm{\thetainit - \hpop(\thetainit)}
\end{align}
for a universal constant $c' > 0$. \\

It can be seen that the sequences $\widetilde{\ratev}$ and
$\widetilde{\ratetheta}$ are $2$-admissible.  Substituting their
definitions into the bound~\eqref{eq:zt-recursive-bound-compostep}.
we find that the inequality
\begin{multline*}
    \norm{\myZ_{t}} \leq \tfrac{c}{\sqrt{t}} \Big \{ \noisegauss +
    \noisevar \sqrt{\log(\tfrac{1}{\delta})} \Big \} + \tfrac{c
      \boundstar}{t} \Big \{ \DudExp(\dualBall, \myrhotil) +
    \log(\tfrac{1}{\delta}) \Big \} \\
+ \tfrac{c \Lip}{t} \Big \{ \constMGtype +
\sqrt{\log(\tfrac{1}{\delta})} \Big \} \; \Big \{ \stepsize
\mybiggersqrt{\sum_{s = \burnin}^{t - 1} s^2 \ratev^2(s)} + 
\compostep \mybiggersqrt{\sum_{s = 1}^{t - 1} \ratetheta^2(s)}
\Big \}
\end{multline*}
holds with probability at least $1 - \delta$.

Under the
conditions~\eqref{eq:stepsize-requirement-bootstrap-compostep-phase-i}
and~\eqref{eq:burn-in-requirement-compostep-in-bootstrap-phase-i},
some algebra yields:
\begin{multline*}
\norm{\myZ_{t}} \leq \tfrac{c}{\sqrt{t}} \Big \{ \noisegauss + \noisevar
\sqrt{\log(\tfrac{1}{\delta})} \Big \} \\ + c \boundstar \Big \{
\tfrac{1}{t} + \tfrac{\stepsize \Lip }{ \sqrt{t}} \cdot \constMGtype
\log(\tfrac{\numobs}{\delta}) \Big \} \; \Big \{\DudExp(\dualBall,
\myrhotil) + \log(\tfrac{1}{\delta}) \Big \} + c \tfrac{\burnin}{t}
\norm{\thetainit - \hpop(\thetainit)}.
\end{multline*}
Combining with equation~\eqref{eq:vt-rate-final-compostep} yields the
upper bound
\begin{multline}
\norm{\hpop(\theta_t) - \theta_t} \leq
\tfrac{c}{\sqrt{t}} \big( \noisegauss + \noisevar
\sqrt{\log(\tfrac{1}{\delta})} \big) \\ + {c \boundstar} \Big \{
\tfrac{1}{t} + \tfrac{\stepsize \Lip }{ \sqrt{t}}
\constMGtype \log (\tfrac{\numobs}{\delta}) \Big \} \; \Big \{
\DudExp(\dualBall, \myrhotil) + \log(\tfrac{1}{\delta}) \Big \} 
+ c\tfrac{\burnin}{t} \norm{\thetainit - \hpop(\thetainit)},\label{eq:bellman-err-linear-without-restart}
\end{multline}
which completes the proof of
equation~\eqref{eq:main-thm-bellman-err-bound-general-linear}. \\

Besides, by taking a union bound over time steps $t \in \{ \burnin, \burnin +
1, \ldots, \numobs \}$, we have the lower bound
$\Prob \big[ \Event^{(\theta)}_\numobs(\ratetheta^*)
\big] \geq 1 - \delta$,
where
\begin{multline*}
    \ratetheta^*(t) \mydefn \tfrac{c}{\sqrt{t}} \Big \{ \noisegauss +
    \noisevar \sqrt{\log (\tfrac{\numobs}{\delta}) } \Big \} \\ + {c
      \boundstar} \Big \{ \tfrac{1}{t} + \tfrac{\stepsize \Lip
      \constMGtype}{ \sqrt{t}} \log (\tfrac{\numobs}{\delta}) \Big \}
    \; \Big \{ \DudExp(\dualBall, \myrhotil) + \log
    (\tfrac{\numobs}{\delta}) \Big \} + c \tfrac{\burnin}{t}
    \norm{\thetainit - \hpop(\thetainit)}.
\end{multline*}


\subsubsection{Proof of Lemma~\ref{lemma:vt-recursive-bound-compostep}}
\label{subsubsec:proof-lemma-vt-recursive-bound-compostep}

Starting with recursion satisfied by
$\myV_t$, we have
\begin{align*}
t \cdot \myV_t & = (t - 1) \big \{ \myV_{t - 1} + \theta_{t - 2} - \Hstoch_t
(\theta_{t - 1}) - \theta_{t - 1} + \Hstoch_t (\theta_{t - 2}) \big \}
+ \big \{ \Hstoch (\theta_{t - 1}) - \theta_{t - 1} \big \} \\
& = \big \{ (1 - \stepsize) I + \stepsize \Amat \big \} \cdot (t - 1)
\myV_{t - 1} - (t - 1) \big \{ \noise_t(\theta_{t - 1}) -
\noise_t(\theta_{t - 2}) \big \} + \noise_t(\theta_{t - 1}) + \big \{
\hpop(\theta_{t - 1}) - \theta_{t - 1} \big \}.
\end{align*}
For any positive integer $\tau$, we can expand the above expression
for $\tau$ steps so as to obtain
\begin{multline}
  \label{EqnKento}
t \cdot \myV_t = \big( (1 - \stepsize) I + \stepsize \Amat \big)^\tau (t
- \tau) \myV_{t - \tau} - \sum_{j = 1}^\tau (t - j) \big( (1 - \stepsize)
I + \stepsize \Amat \big)^{j - 1} \left( \noise_{t - j + 1} (\theta_{t
  - j}) - \noise_{t - j + 1} (\theta_{t - j - 1}) \right) \\
+ \sum_{j = 1}^\tau \big( (1 - \stepsize) I + \stepsize \Amat \big)^{j
  - 1} \noise_{t - j + 1} (\theta_{t - j}) + \sum_{j = 1}^\tau \big(
(1 - \stepsize) I + \stepsize \Amat \big)^{j - 1} \big(
\hpop(\theta_{t - j}) - \theta_{t - j} \big).
\end{multline}
In addition, our analysis makes use of the following auxiliary bound
\begin{align}
\label{eq:operator-norm-compostep}  
  \matsnorm{\big( (1 - \stepsize) I + \stepsize \Amat
    \big)^t}{\vecspace} &\leq \min \Big \{  1 , 2  \left(1 - \tfrac{\stepsize}{2
    \compostep} \right)^t \Big \},
\end{align}
valid for all $t = 1, 2, \ldots$.  See the end of this subsection for
the proof of this claim.

Taking this bound as given, we proceed with the proof of this
lemma. First, substituting the
bound~\eqref{eq:operator-norm-compostep} into the
decomposition~\eqref{EqnKento} yields the bound
\begin{align}
\label{eq:compostep-vt-decomposition} 
 t \cdot \norm{\myV_t} \leq 2  \left(1 - \tfrac{\stepsize}{2
   \compostep} \right)^{\tau} (t - \tau) \norm{\myV_{t - \tau}} +
 \norm{\Psi_{t - \tau, \tau}} + \norm{\Martin_{t - \tau, \tau }} +
  \sum_{j = 1}^\tau \norm{ \hpop(\theta_{t - j}) - \theta_{t -
     j} }.
\end{align}
where we define the terms
\begin{subequations}
  \begin{align}
 \Psi_{t -\tau, \tau} &\defn \sum_{j = 1}^\tau (t - j) \big( (1 -
 \stepsize) I + \stepsize \Amat \big)^{j - 1} \left( \noise_{t - j +
   1} (\theta_{t - j}) - \noise_{t - j + 1} (\theta_{t - j - 1})
 \right), \quad \mbox{and} \\
\Martin_{t - \tau, \tau} &\defn \sum_{j = 1}^\tau \big( (1 - \stepsize)
I + \stepsize \Amat \big)^{j - 1} \noise_{t - j + 1} (\theta_{t -
  j}).
  \end{align}
\end{subequations}

Now we bound the terms in the
decomposition~\eqref{eq:compostep-vt-decomposition}.  On the event
$\bigvEvent$, each term in the summation defining $\Psi_{t -\tau,
  \tau}$ satisfies an almost-sure upper bound:
\begin{align*}
\big( (1 - \stepsize) I + \stepsize \Amat \big)^{j - 1} \left(
\noise_{t - j + 1} (\theta_{t - j}) - \noise_{t - j + 1} (\theta_{t -
  j - 1}) \right) \leq (t - j)  \Lip \stepsize \norm{\myV_{t - j}}
\leq (t - j)  \Lip \stepsize \ratev(t - j).
\end{align*}
Since the sequence $\ratev$ is admissible, for burn-in time $\burnin
\geq 2 \tau$, we have that $(t - j) \ratetheta(t - j) \leq
\tfrac{t^2}{(t - j)} \ratev(t) \leq 2 t \ratev(t)$. Note that the
terms in $\Psi_{t - \tau, \tau}$ form a martingale difference
sequence, adapted to the natural filtration $(\filtration_{t})_{t \geq
  0}$. Invoking the martingale concentration inequality from
Lemma~\ref{lemma:concentration-in-banach-space-martingale} yields the
bound
\begin{align}
\norm{\Psi_{t - \tau, \tau}} \leq c \sqrt{\tau} \big \{ \constMGtype +
\sqrt{\log(1 / \delta)} \big \} \cdot  \Lip \stepsize t
\ratev(t),
\end{align}
which holds with probability at least $1 - \delta$.

As for the term $\Martin_{t - \tau, \tau}$, we use a decomposition
similar to the one used in the proof of
Lemma~\ref{lemma-vt-recursive-bound}:
\begin{align*}
\Martin_{t - \tau, \tau} &= \sum_{j = 1}^\tau \big( (1 - \stepsize) I
+ \stepsize \Amat \big)^{j - 1} \noise_{t - j + 1} (\thetastar) +
\sum_{j = 1}^\tau \big( (1 - \stepsize) I + \stepsize \Amat \big)^{j -
  1} \big( \noise_{t - j + 1} (\theta_{t - j}) - \noise_{t - j + 1}
(\thetastar) \big) \\
& =: \MartinStar_{t - \tau, \tau} + \widetilde{M}_{t - \tau, \tau}.
\end{align*}
The term $\MartinStar_{t - \tau, \tau}$ is sum of independent random
variables in $\vecspace$, with each term satisfying the conditions
\begin{align*}
\norm{ \big( (1 - \stepsize) I + \stepsize \Amat \big)^{j - 1}
  \noise_{t - j + 1} (\thetastar)} \leq  \cdot \norm{\noise_{t
    - j + 1} (\thetastar)}, \quad \mbox{and} \quad (1 - \stepsize) I +
\stepsize \Amat \big)^{j - 1} \noise_{t - j + 1} (\thetastar) \in
 \obssubset.
\end{align*}
Invoking the concentration inequality from
Lemma~\ref{lemma:bernstein-in-banach-space-iid} yields the bound
\begin{align*}
  \norm{\MartinStar_{t - \tau, \tau}} \leq c  \sqrt{\tau}
  \left( \noisegauss + \noisevar \sqrt{\log(\tfrac{1}{\delta})} \right) + c 
  \boundstar \left( \DudExp(\dualBall, \myrhotil) + \log(\tfrac{1}{\delta})
  \right),
\end{align*}
which holds with probability at least $1 - \delta$.

For the excess noise term $\widetilde{M}_{t - \tau, \tau}$, we note
that conditioned on the event $\bigthetaEvent$, we have the upper
bound
\begin{align*}
\Big \| \big( (1 - \stepsize) I + \stepsize \Amat \big)^{j - 1} \big(
\noise_{t - j + 1} (\theta_{t - j}) - \noise_{t - j + 1} (\thetastar)
\big) \Big\| & \leq 2 \Lip \compostep \ratetheta(t - j).
\end{align*}
For an admissible sequence $\ratetheta$ and burn-in period $\burnin
\geq 2 \tau$, we have that $\ratetheta(t - j) \leq \tfrac{t^2}{(t -
  j)^2} \ratetheta(t) \leq 4 \ratetheta(t)$ for any $j \in
     [\tau]$. Furthermore, the terms in $\widetilde{M}_{t - \tau,
       \tau}$ form a martingale difference sequence adapted to the
     natural filtration. By
     Lemma~\ref{lemma:concentration-in-banach-space-martingale}, on
     the event $\bigthetaEvent$, we have the martingale concentration
     inequality:
\begin{align*}
    \norm{\widetilde{M}_{t - \tau, \tau}} \leq c \sqrt{\tau}  \big( \constMGtype + \sqrt{\log(1 / \delta)} \big) \cdot \Lip \compostep \ratetheta(t).
\end{align*}
Finally, for the last term in the decomposition~\eqref{eq:compostep-vt-decomposition}, we note that on the event $\bigthetaEvent$, we have the bounds:
\begin{align*}
    \norm{ \hpop(\theta_{t - j}) - \theta_{t - j} } \leq \ratetheta(t - j) \leq \tfrac{t^2}{(t - j)^2} \ratetheta(t) \leq 4 \ratetheta(t).
\end{align*}
In order to prove the final results, as with the proof of Lemma~\ref{lemma-vt-recursive-bound}, we consider the cases of $t \geq \burnin + 2\compostep / \stepsize$ and $t \leq \burnin +2 \compostep / \stepsize$ separately.

When $t \geq \burnin + 2 \compostep / \stepsize$, collecting above bounds, by taking $\tau =
2\compostep / \stepsize $, we find that
\begin{multline*}
    t \cdot \norm{\myV_t} \leq \Big \{ \tfrac{1}{3} + c \sqrt{\tau}
    \big[ \constMGtype + \sqrt{\log(\tfrac{1}{\delta})} \big] \cdot
     \Lip \stepsize \Big \} t \ratev(t) \\
+  \Big \{ c \sqrt{\tau}  \big[ \constMGtype +
  \sqrt{\log(\tfrac{1}{\delta})} \big] \cdot \Lip \compostep + \tau
\Big \} \ratetheta(t) \\
+ c  \sqrt{\tau} \Big \{ \noisegauss +
\noisevar \sqrt{\log(\tfrac{1}{\delta})} \Big \} + c  \boundstar \Big \{
\DudExp(\dualBall, \myrhotil) + \log(\tfrac{1}{\delta}) \Big \}.
\end{multline*}
Given a stepsize $\stepsize$ such that
\begin{align}
  c \sqrt{\compostep \stepsize} \cdot \Big \{
  \constMGtype + \sqrt{\log( \tfrac{1}{\delta})} \Big \} \cdot 
  \Lip \leq \tfrac{1}{3},\label{eq:vt-recursive-lemma-compostep-stepsize-requirement}
\end{align}
the above inequality implies that
\begin{align*}
t \cdot \norm{\myV_t} \leq \tfrac{2}{3} t \ratev(t) + \tfrac{c \compostep
}{\stepsize} \ratetheta(t) + c 
\sqrt{\tfrac{ \compostep}{\stepsize} } \Big \{
\noisegauss + \noisevar \sqrt{\log( \tfrac{1}{\delta})} \Big \} + c 
\boundstar \Big \} \DudExp(\dualBall, \myrhotil) +
\log(\tfrac{1}{\delta}) \Big \},
\end{align*}
which completes the proof of the first case.

On the other hand, when $t \leq \burnin + 2 \compostep / \stepsize$, we let $\tau = t - \burnin$, and find that:
\begin{multline*}
        t \cdot \norm{\myV_t} \leq 2 \burnin \cdot \norm{\myV_{\burnin}} + c \sqrt{\tau}
    \big[ \constMGtype + \sqrt{\log(\tfrac{1}{\delta})} \big] \cdot
     \Lip \stepsize  t \ratev(t) \\
+  \Big \{ c \sqrt{\tau}  \big[ \constMGtype +
  \sqrt{\log(\tfrac{1}{\delta})} \big] \cdot \Lip \compostep + \tau
\Big \} \ratetheta(t) \\
+ c  \sqrt{\tau} \Big \{ \noisegauss +
\noisevar \sqrt{\log(\tfrac{1}{\delta})} \Big \} + c  \boundstar \Big \{
\DudExp(\dualBall, \myrhotil) + \log(\tfrac{1}{\delta}) \Big \}.
\end{multline*}
Note that for $t \in [\burnin, \burnin + 2 \compostep / \stepsize]$, we have that $ 2 \burnin \cdot \norm{\myV_{\burnin}} \leq 4 \frac{\burnin^2}{t} \norm{\myV_{\burnin}}$. Assuming the stepsize condition~\eqref{eq:vt-recursive-lemma-compostep-stepsize-requirement}, we conclude the inequality:
\begin{align*}
t \cdot \norm{\myV_t} \leq \tfrac{2}{3} t \ratev(t) + \tfrac{c \compostep
}{\stepsize} \ratetheta(t) + c 
\sqrt{\tfrac{ \compostep}{\stepsize} } \Big \{
\noisegauss + \noisevar \sqrt{\log( \tfrac{1}{\delta})} \Big \} + c 
\boundstar \Big \} \DudExp(\dualBall, \myrhotil) +
\log(\tfrac{1}{\delta}) \Big \} + \frac{\burnin^2}{t} \norm{\myV_{\burnin}},
\end{align*}

\paragraph{Proof of equation~\eqref{eq:operator-norm-compostep}:}

Applying the triangle inequality yields
\begin{align}
\matsnorm{\big( (1 - \stepsize) I + \stepsize \Amat
  \big)^t}{\vecspace} \leq \sum_{k = 0}^t \binom{t}{k} (1 -
\stepsize)^k \stepsize^{t - k} \matsnorm{\Amat^{t -
    k}}{\vecspace}.\label{eq:proof-of-operator-norm-compostep-decomp}
\end{align}
Since $\matsnorm{\Amat^t}{\vecspace} \leq \matsnorm{\Amat}{\vecspace}^t \leq 1$ for each $t = 0, 1,
2, \ldots$, we have
\begin{align*}
\matsnorm{\big( (1 - \stepsize) I + \stepsize \Amat
  \big)^t}{\vecspace} \leq \sum_{k = 0}^t \binom{t}{k} (1 -
\stepsize)^k \stepsize^{t - k}   \leq 1.
\end{align*}
  On the other hand, we note that for any time index $i \in \mathbb{N}_+$, using the $\compostep$-step contraction
  condition~\ref{new--assume-compostep-contractive}, we have that:
\begin{align*}
\matsnorm{\Amat^i}{\vecspace} \leq
\matsnorm{\Amat^\compostep}{\vecspace}^{\lfloor \tfrac{i}{\compostep} \rfloor} \cdot
\matsnorm{\Amat^{i - \compostep \lfloor \tfrac{i}{\compostep} \rfloor}}{\vecspace} \leq 2^{- \lfloor \frac{i}{\compostep} \rfloor}   = 2^{1 - i / \compostep} .
\end{align*}
Applying this inequality with $i = t - k$ and substituting into
equation~\eqref{eq:proof-of-operator-norm-compostep-decomp}, we have
that:
\begin{align*}
\matsnorm{\big( (1 - \stepsize) I + \stepsize \Amat
  \big)^t}{\vecspace} \leq \sum_{k = 0}^t \binom{t}{k} (1 - \stepsize)^k
\stepsize^{t - k} \cdot 2^{1 - \tfrac{t - k}{\compostep}}
\leq 2  \left( 1 - \stepsize + \stepsize \cdot \big( 1 -
\tfrac{1}{2 \compostep} \big) \right)^t = 2 \left( 1 -
\tfrac{\stepsize}{2 \compostep} \right)^t.
\end{align*}


\subsection{Proof of Corollary~\ref{thm:est-err-linear}}
\label{subsec:proof-thm-est-err-linear}

In this section, we prove the stated claim with the higher-order
term defined as
\begin{multline*}
\SuperHot{\diamond} = c \compostep \aUB \Big \{ \Lip \constMGtype
\log(\tfrac{\numobs}{\delta}) \sqrt{\tfrac{\stepsize
    \compostep}{\numobs}} + \tfrac{1}{\numobs}
\sqrt{\tfrac{\compostep}{\stepsize}} \Big \} \; \Big\{ \noisegauss +
\noisevar \sqrt{\log(\tfrac{\numobs}{\delta})} \Big \} \\
+ c \compostep
\boundstar \aUB \Big \{ \tfrac{1}{\numobs} + \tfrac{\stepsize \Lip
  \constMGtype}{\sqrt{\numobs}} \log(\tfrac{\numobs}{\delta}) \Big \}
\; \Big\{ \DudExp(\dualBall, \myrhotil) +
\log(\tfrac{\numobs}{\delta}) \Big\}
\end{multline*}

Recall that by Theorem~\ref{thm:bellman-err-bound-general-linear} and
a union bound, for the restarting procedure described
in~\Cref{AppRestart}, the event $\Event^{(\theta)}_\numobs
(\ratetheta^*) \cap \Event^{(v)}_\numobs (\ratev^*)$ occurs with
probability $1 - \delta$, for the function pair $(\ratetheta^*,
\ratev^*)$ given by
\begin{subequations}
\begin{align}
\label{eq:ratev-copied-linear} 
\ratev^*(t) & \mydefn c \Big[ \tfrac{1}{t}
  \sqrt{\tfrac{\compostep}{\stepsize}} \Big( \noisegauss + \noisevar
  \sqrt{\log (\tfrac{\numobs}{\delta})} \Big) + \tfrac{\boundstar}{t}
  \Big \{ \log (\tfrac{\numobs}{\delta}) + \DudExp(\dualBall, \myrhotil)
  \Big \} \Big] \\
\label{eq:ratetheta-copied-linear}
\ratetheta^*(t) & \mydefn \tfrac{c}{\sqrt{t}} \Big \{( \noisegauss +
\noisevar \sqrt{\log(\tfrac{\numobs}{\delta})} \Big \} + {c
  \boundstar} \Big \{ \tfrac{1}{t} + \tfrac{\stepsize \Lip
  \constMGtype}{ \sqrt{t}} \log^3 (\tfrac{\numobs}{\delta}) \Big \} \;
\Big \{ \DudExp(\dualBall, \myrhotil) + \log(\tfrac{1}{\delta}) \Big \}.
\end{align}
\end{subequations}
Since $\hpop$ is an affine operator, we have the decomposition
\begin{align*}
\anorm{\theta_\numobs - \thetastar} \leq \anorm{(I - \Amat)^{-1}
  \myV_{\numobs + 1}} + \anorm{(I - \Amat)^{-1} \myZ_{\numobs + 1}}.
\end{align*}
By the operator norm bound~\eqref{eq:opnorm-bound-compostep} and the
bound~\eqref{eq:ratev-copied-linear} on the norm $\norm{\myV_t}$, we have
\begin{align*}
\anorm{(I - \Amat)^{-1} \myV_{t + 1}} \leq c' \aUB \compostep \Biggr[
  \tfrac{1}{t} \sqrt{\tfrac{\compostep}{\stepsize}} \Big\{\noisegauss
  + \noisevar \sqrt{\log (\tfrac{\numobs}{\delta})} \Big\} +
  \tfrac{\boundstar}{t} \big \{ \log (\tfrac{\numobs}{\delta}) +
  \DudExp(\dualBall, \myrhotil) \big \} \Biggr].
\end{align*}

For the term $\norm{(I - \Amat)^{-1} \myZ_{\numobs + 1}}$, we consider
the class of test functions $\UpsSet \mydefn \big \{ (I -
\Amat^*)^{-1} u \, \mid \, u \in \MySet \big \}$.  Invoking
Lemma~\ref{lemma-zt-martingale-bound-anorm} with $\MuCon = (1-
\contraction)$ yields that $\anorm{(I - \Amat^{-1}) \myZ_\numobs}$ is at
most
\begin{multline*}
\tfrac{c}{\sqrt{\numobs}} \Big\{ \Exs \big[ \anorm{(I - \Amat)^{-1} W}
  \big] + \mysqrt{\sup_{u \in \MySet} \Exs \big[ \inprod{u}{(I -
      \Amat) W}^2 \big] \log(\tfrac{1}{\delta})} \Big \} + c \tfrac{
  \aUB \compostep \boundstar}{\numobs} \Big \{ \DudExp(\dualBall,
\myrhotil) + \log(\tfrac{1}{\delta}) \Big \} \\
+ \tfrac{c \aUB \compostep\Lip}{t} \Big \{ \constMGtype +
\mybiggersqrt{\log(\tfrac{1}{\delta})} \Big \} \; \Big \{ \stepsize
\mysqrt{\sum_{s = \burnin}^{\numobs - 1} s^2 \ratev^* (s)^2} +
\compostep \mysqrt{\sum_{s = 1}^{\numobs - 1} \ratetheta^*(s)^2} \Big
\},
\end{multline*}
with probability at least $1 - \delta$. Combining above results, some
algebra yields that
\begin{multline*}
\anorm{\theta_\numobs - \thetastar} \leq \tfrac{c}{\sqrt{\numobs}}
\Big\{ \Exs \big[ \anorm{(I - \Amat)^{-1} W} \big] +
\mybiggersqrt{\sup_{u \in \MySet} \Exs \big[ \inprod{u}{(I - \Amat)
      W}^2 \big] \log(\tfrac{1}{\delta})} \Big\} \\
+ c \compostep \aUB \Big\{ \Lip \constMGtype
\log(\tfrac{\numobs}{\delta}) \sqrt{\tfrac{\stepsize
    \compostep}{\numobs}} + \tfrac{1}{\numobs}
\sqrt{\tfrac{\compostep}{\stepsize}} \Big \} \; \Big \{ \noisegauss +
\noisevar \sqrt{\log(\tfrac{\numobs}{\delta})} \Big \} \\
+ c \compostep \boundstar \aUB \Big \{ \tfrac{1}{\numobs} +
\tfrac{\stepsize \Lip \constMGtype}{\sqrt{\numobs}}
\log(\tfrac{\numobs}{\delta}) \Big \} \; \Big \{ \DudExp(\dualBall,
\myrhotil) + \log(\tfrac{\numobs}{\delta}) \Big \},
\end{multline*}
with probability at least $1 - \delta$.  This completes the proof of
Corollary~\ref{thm:est-err-linear}.


\section{Two-player zero-sum Markov games}
\label{sec:Markov-games}

In this section, we explore the consequences of our general theory for
two-player zero-sum Markov games.  This class of problems results from
a marriage between MDPs and two player zero-sum games: it is used the
model two agents who play multiple rounds of a zero-sum game, and each
has the goal to maximize their expected long-term reward.  Markov
games are characterized by a six-tuple $\{\states, \actionsOne,
\actionsTwo, \TranMG, \reward, \contraction\}$. Let $\states$ denote
the state space, and let $\actionsOne$ and $\actionsTwo$ denote the
action sets for players one and two, respectively. Here we focus on
games with finite state and action space, i.e., \mbox{$|\states \times
  \actionsOne \times \actionsTwo| < \infty$}.

The probability transition kernel $\bracketMed{\TranMG_{\action_1,
    \action_2}(\state' \mid \state) \mid (\state, \action_1,
  \action_2) \in \states \times \actionsOne \times \actionsTwo}$,
encodes the transition to the next state given the actions of the
players. In particular, the scalar $\TransMat_{\action}(\state' \mid
\state)$ denotes the probability of transition to the state $\state'$,
when at state $\state$ player 1 takes the action $\action_1$ and
player 2 takes the action $\actions_2$.  The MDP is equipped with a
reward function $\reward: \states \times \actionsOne \times
\actionsTwo \mapsto \real$ such that the scalar $\reward(\state,
\action_1, \action_2)$ denotes the cost received at state $\state$
when player 1 takes the action $\action_1$ and player 2 takes the
action $\action_2$. Finally, the scalar $\contraction \in (0, 1)$ is a
parameter reflecting the discounting of future rewards.

For each player $i \in \{1, 2\}$, a stationary policy $\policy_i$ is a
mapping $\states \rightarrow \ProbDistr{\actions_i}$, where
$\ProbDistr{\actions_i}$ denotes the set of probability distributions
over the finite action set $\actions_i$.  In other words, the actions
taken by the players can be random, and for any state $\state \in
\states$, the distribution $\policy_i(\cdot \mid \state)$ is a
probability distribution on the set of actions $\actions_i$ to be
taken by player $i$. We use $\Pi_1$ and $\Pi_2$ to denote the set of
all policies for players $1$ and $2$, respectively.

Assuming player $1$ is following policy $\policy_1$, and player $2$ is
following policy $\policy_2$, the value $\valueMG(\cdot \mid
\policy_1, \policy_2) : \real^{|\states|} \mapsto \real$ of player $1$
is defined as the expected sum of discounted rewards in an infinite
sample path:
\begin{align}
\label{eqn:value-for-player-one}
\valueMG(\state \mid \policy_1, \policy_2) = \Exs \Big[ \sum_{k =
    1}^{\infty} \reward(\state_k, \action_{1k}, \action_{2k} \mid
  \state_0 = \state) \Big], \quad \mbox{where $\action_{1k}\sim
  \policy_1(\cdot \mid \state_k)$ and $\action_{2k} \sim
  \policy_2(\cdot \mid \state_k)$.}
\end{align}

Given that the game is zero-sum, the reward for player $2$ with initial
state $\state$ is $-\valueMG(\state \mid \policy_1,
\policy_2)$. Players $1$ and $2$ want to choose their policies
$\policy_1$ and $\policy_2$ that maximize their respective reward for all values
of initial state $\state$.

\paragraph{Nash equilibrium:} A natural notion of equilibrium
in two-player zero-sum Markov games is the Nash equilibrium.  A policy
pair $(\policy_1^\star, \policy_2^\star)$ is called a \emph{Nash
equilibrium} if for all initial states $\state \in \states$
\begin{align}
\label{eqn:Nash-eqb-MG}
\valueMG(\state \mid \policy_1^\star, \policy_2^\star) &\geq
\valueMG(\state \mid \policy_1, \policy_2^\star) \quad \text{for all
  policies} \quad \policy_1 \in \policySet_1, \qquad \text{and} \notag
\\
- \valueMG(\state \mid \policy_1^\star, \policy_2^\star) & \geq
- \valueMG(\state \mid \policy_1^\star, \policy_2)
\quad \text{for all policies} \quad \policy_2 \in \policySet_2.
\end{align}
In words, the policy $\policy_1^\star$ is the best response for player
$1$ assuming player $2$ is playing policy $\policy_2^\star$, and the
policy $\policy_2^\star$ is the best response for player $2$ assuming
player $1$ is playing policy $\policy_1^\star$. Thus, neither player
has any incentive to deviate from the policy pair $(\policy_1^\star,
\policy_2^\star)$. In two-player zero-sum Markov games, a Nash
equilibrium always exists, and it is equivalent to the minimax
solution~\cite{perolat2015approximate,patek1997stochastic}. Concretely,
there exist policies $(\policy_1^\star, \policy_2^\star)$ such that
\begin{align}
  \label{eqn:MG-minimax-eqb}
  \valueMG^\star(\state) = \valueMG(\state \mid \policy_1^\star,
  \policy_2^\star) = \min_{\policy_1} \max_{\policy_2} \valueMG(\state
  \mid \policy_1, \policy_2) = \max_{\policy_1} \min_{\policy_2}
  \valueMG(\state \mid \policy_1, \policy_2) \;\; \text{for all} \;\;
  \state \in \states.
\end{align}
The function $\valueMG^\star$ is known as the \emph{value} of the game.


\subsection{$Q$-function and the Bellman fixed-point equation}
\label{sec:Bellman-MG}

One method for finding a pair of policies $(\policy_1^\star,
\policy_2^\star)$ that achieves the
equilibrium~\eqref{eqn:MG-minimax-eqb} is by computing the optimal
state-action value functions or the optimal $Q$-function $\QfunMG$.
It is known~\cite{patek1997stochastic,perolat2015approximate} to be
the fixed point of the Bellman operator
\begin{multline}
\label{eqn:Bellman-MG}
\MGOp(\theta)(\state, \action_1, \action_2) = \cost(\state,
\action_1, \action_2) \\
 + \contraction \cdot \sum_{\state' \in
  \states} \TranMG_{\action_1, \action_2}(\state' \mid \state)
\max_{\policy_1} \min_{\policy_2} \sum_{\action_1^\prime,
  \action_2^\prime} \policy_1(\action_1^\prime \mid \state') \cdot
\policy_2(\action_2^\prime \mid \state') \cdot \theta(\state',
    \action_1^\prime, \action_2^\prime).
\end{multline}
Notably, when the number of states and actions are finite, the minimax
problem on the right-hand side of equation~\eqref{eqn:Bellman-MG} can
be computed by solving the two-player zero-sum matrix game with the
payoff matrix $\{\theta(\state', \action_1, \action_2 ) \mid \action_1
\in \actions_1, \action_2 \in \actions_2\}$.  Finally, for Markov
games with finite state and action spaces, the $Q$-function $\theta$
can be conveniently represented as an element of
$\mathbf{R}^{|\states| \times|\actions_1| \times |\actions_2|}$, and
the Bellman operator $\MGOp$ is an operator on $\real^{|\states|
  \times|\actions_1| \times |\actions_2|}$.

A simple calculation yields that the Bellman operator is
$\contraction$-contractive in the
$\ell_\infty$-norm~\cite{patek1997stochastic,perolat2015approximate},
and as a result, the optimal $Q$-function is the unique fixed point of
the operator $\MGOp$.  We can thus apply our general Banach space
theory to derive bounds on the \rootSA procedure.

\subsection{The generative model and empirical Bellman operator}
\label{sec:MG-generative-model}

We analyze the behavior of the \rootSA algorithm under a stochastic
oracle known as the generative model.  A sample from this model
consists of a pair of real-valued tensors $(\TranSample,
\randReward)$, each with dimensions $|\states| \times |\actions_1|
\times |\actions_2|$. For each triple $(\state, \action_1,
\action_2)$, the entry $\TranSample_(\state, \action_1, \action_2)$ is
drawn according to the transition kernel $\TransMat_{\action_1,
  \action_2}(\cdot \mid \state)$, whereas the entry
$\randReward(\state, \action_1, \action_2)$ is a zero-mean random
variable with mean $\reward(\state, \action_1, \action_2)$,
corresponding to a noisy observation of the reward function.  The
transition and reward samples across entries of the tensors are
independently sampled, and we assume that the rewards are bounded in
absolute value by $\rewardMax$.

Given a sample $(\TranSample, \randReward)$ from our observation
model, we can define the single-sample empirical Bellman operator
\begin{multline}
\label{eqn:Games-noisy-Bellman}
\SSPstoch(\theta)(\state,\action_1, \action_2) \defn
\randReward(\state,\action_1, \action_2) \\ + \sum_{\state' \in \states}
\TranSample_{\action_1, \action_2}\left(\state' \mid \state \right)
\max_{\policy_1} \min_{\policy_2} \sum_{\action_1^\prime,
  \action_2^\prime} \policy_1(\action_1^\prime \mid \state') \cdot
\policy_2(\action_2^\prime \mid \state') \cdot \theta(\state',
\action_1^\prime, \action_2^\prime),
\end{multline}
where we have introduced the notation $\TranSample_{\action_1,
  \action_2}\left(\state' \mid \state\right) \defn
\mathbf{1}_{\TranSample(\state, \action_1, \action_2)=\state'}$.  With
these definitions in hand, we are now ready to state our guarantees
for two-player zero-sum Markov games.


\subsection{Guarantees for two-player zero-sum Markov games}

Let $\GaussNoise$ be a zero-mean Gaussian random vector with
covariance $\mathrm{cov}(\Hstoch (\thetastar) - \thetastar)$, and
define
\begin{align}
\label{eqn:key-complexity-games}
\noisegauss = \Exs[\|W\|_\infty], \qquad \noisevar^2 \mydefn
\sup_{\state \in \states, \action_1 \in \actions_1, \action_2 \in
  \actions_2} \Exs [\GaussNoise_{\state, \action_1, \action_2}^2],
\quad \mbox{and} \quad \boundstar \mydefn \rewardMax +
\infNorm{\thetastar}.
\end{align}
For a given failure probability $\delta \in (0,1)$, our result applies
to the algorithm with parameters
\begin{subequations}
\label{eqn:Game-tuning}
\begin{align}
\label{eqn:Games-tuning}  
\stepsize = c_1 \Big\{\sqrt{ \numobs \log \abss{\states \times
    \actionsOne \times \actionsTwo}} \cdot
\log(\tfrac{\numobs}{\delta}) \Big\}^{-1}, \quad \mbox{and} \quad
\burnin = \tfrac{c_2}{(1 - \contraction)^2 \stepsize}
\log(\tfrac{\numobs}{\delta}).
\end{align}
We also choose the initialization $\theta_0$ and the number of
restarts $\numRestarts$ such that
\begin{align}
\label{eqn:numrestarts-and-init-Game}
 \log \left( \tfrac{\norm{\theta_0 - \hpop (\theta_0)}
   \sqrt{\numobs}}{\noisegauss} \right) \leq c_0 \log \numobs \qquad
 \mbox{and} \qquad \numRestarts \geq 2 c_0 \log \numobs
\end{align}
\end{subequations}
for appropriate universal constants $c_0$, $c_1$ and $c_2$.  With this
setup, a direct application of
Theorem~\ref{new--thm:bellman-err-bound-general-nonlinear} yields the
following:
\begin{corollary}
\label{cors:game-cors}
Given a sample size $\numobs$ such that $\tfrac{\numobs}{\log \numobs}
\geq \tfrac{c' \log(\abss{\states} \cdot \abss{\actionsOne} \cdot
  \abss{\actionsTwo})}{(1 - \contraction)^4} \log(
\tfrac{1}{\delta})$, running Algorithm~\ref{AlgRootSGD} with the
tuning parameter choices~\eqref{eqn:Game-tuning} yields an estimate
$\theta_\numobs$ such that
\begin{align*}
\myinfnorm{\hpop(\theta_\numobs) - \theta_\numobs} & \leq
\tfrac{c}{\sqrt{\numobs}} \cdot \Big \{ \noisegauss + \noisevar \sqrt{
  \log(\tfrac{1}{\delta})} \Big \} + \tfrac{c \boundstar}{1 -
  \contraction} \cdot \tfrac{\log (|\states| \cdot \abss{\actionsOne}
  \cdot \abss{\actionsTwo})}{\numobs} \log^2(\tfrac{\numobs}{\delta}).
\end{align*}
with probability at least $1 - \delta$.
\end{corollary}
Note that the bound in Corollary~\eqref{cors:game-cors} depends on the
size of state-action space $|\states| \cdot \abss{\actionsOne} \cdot
\abss{\actionsTwo}$ only poly-logarithmically. Moreover, one can
obtain an upper bound on the estimation error $\infNorm{\theta_\numobs
  - \thetastar}$ using the
bound~\eqref{new--eqn:sub-optimality-to-estimation-err}.

A special case of interest is when the set of actions for player two
is a singleton, i.e., $|\actionsTwo| = 1$. Observe that in this case
the optimal state-action value estimation problem for the two-player
zero-sum Markov game reduces to the optimal value estimation problem
of an appropriate MDP in the discounted
setting~\cite{bertsekas2019reinforcement,watkins1992q,wainwright2019variance}. In
Appendix~\ref{sec:localLinQProof}, we show that the Bellman operator
associated with the optimal value estimation problem of an MDP in the
discounted setting satisfies the local linearity
assumption~\ref{new--assume:linearization}.  Consequently, an argument
similar to Corollary~\ref{cors:SSP-cors} yields an upper bound on the
estimation error $\infNorm{\theta_\numobs - \thetastar}$ which matches
the instance dependent lower bound (up to logarithmic terms) from the
paper~\cite{khamaru2021instance} for large $\numobs$.\footnote{ The
sample size requirement for achieving the lower
bound~\cite{khamaru2021instance} may depend on the gap between the
value of optimal and sub-optimal actions.  } Finally, it is an
important direction of future work to investigate whether the local
linearity assumption~\ref{new--assume:linearization} holds when
$|\actionsTwo| > 1$.


\section{Some Concentration Inequalities in Banach Spaces}

Our analysis makes use of some concentration inequalities for
Banach-space-valued random variables, which we state and prove here.

\subsection{Statement of the results}

We begin with a bound for a sequence $\{X_i\}_{i=1}^\numobs$ of
i.i.d. zero-mean random elements.  Our bound involves a zero-mean
Gaussian random variable $W$ in $\vecspace$ such that
\begin{align*}
  \Exs \left[ \inprod{W}{y} \cdot \inprod{W}{z} \right] = \Exs \left[
    \inprod{X_1}{y} \cdot \inprod{X_1}{z} \right] \qquad \mbox{for all
    $y, z \in \vecspace^*$.}
\end{align*}

\begin{lemma}
\label{lemma:bernstein-in-banach-space-iid}
Let $\{X_i\}_{i=1}^\numobs$ be independent zero-mean random elements
taking values in $\obssubset \subseteq \vecspace$ with $\norm{X_i}
\leq 1$ almost surely for each $i = 1,2, \cdots, \numobs$.  Then there
exists a universal constant $c > 0$ such that for any $\delta \in (0,
1)$ and any bounded symmetric convex set $\UpsSet \subseteq \dualBall$,
we have
\begin{align}
\label{EqnBanachIID}  
\tfrac{1}{\numobs} \sup_{u \in \UpsSet} \inprod{u}{\sum_{i =
    1}^\numobs X_i} \leq \tfrac{c}{\sqrt{\numobs}} \Big \{ \Exs
\big[\sup_{u \in \UpsSet} \inprod{u}{W} \big] + \sqrt{ \sup_{u \in
    \UpsSet} \Exs [ \inprod{u}{W}^2] \cdot \log(\tfrac{1}{\delta}) }
\Big \} + \tfrac{c}{\numobs} \Big \{ \log(\tfrac{1}{\delta}) +
\DudExp(\UpsSet, \myrhotil) \Big \},
\end{align}
with probability at least $1 - \delta$.
\end{lemma}
\noindent See Appendix~\ref{SecBanachIID} for the proof of this
claim. \\

\noindent We next state a bound for the martingale case:
\begin{lemma}
\label{lemma:concentration-in-banach-space-martingale}
Let $\{X_t\}_{t=1}^\numobs$ be a martingale in $\vecspace$ adapted to
the filtration $\{\filtration_t\}_{t=1}^\numobs$. Assume that there
exists a deterministic sequence $\{b_t\}_{t=1}^\numobs$ such that $b_t
\geq \tfrac{1}{\numobs}$ and $\norm{X_t} \leq b_t$ almost surely for
each $t = 1,2, \cdots, \numobs$. Then there exists a universal
constant $c > 0$ such that for any $\delta \in (0, 1)$
\begin{align}
\label{EqnBanachMartingale}
\| \sum_{i = 1}^\numobs X_i\| & \leq c \Big( \DudGauss(\dualBall,
\myrhotil) + \sqrt{\log(1 / \delta)} \Big) \cdot \sqrt{\sum_{i =
    1}^\numobs b_i^2},
\end{align}
with probability at least $1 - \delta$.
\end{lemma}
\noindent See Appendix~\ref{SecBanachMartingale} for the proof of this
claim. \\


\subsection{Proof of Lemma~\ref{lemma:bernstein-in-banach-space-iid}}
\label{SecBanachIID}

Our proof is based on a combination of Talagrand's concentration
inequality~\cite{talagrand1996new}, the generic
chaining~\cite{talagrand2006generic} and a functional Bernstein
inequality~\cite{wainwright2019high}.  The left-hand-side of the
desired inequality is the supremum of an empirical process.  Define
the associated Rademacher complexity \mbox{$\localrad_\numobs(\UpsSet)
  \mydefn \tfrac{1}{\numobs} \Exs[\sup \limits _{y \in \UpsSet}
    R_\numobs(y)]$,} where $R_\numobs(y) \mydefn \avg{\rade_i
  \inprod{y}{\X_i}}$ with $\{\rade_i \}_{i=1}^\numobs$ an
i.i.d. sequence of Rademacher random variables. The expectation is
taken over the randomness of both the Rademacher sequence
$\{\rade_i\}_{i=1}^\numobs$ and the random elements
$(\X_i)_{i=1}^\numobs$.

Our first lemma is a type of functional Bernstein inequality; it
bounds the supremum of the empirical process by the Rademacher
complexity and some additional deviation terms:
\begin{lemma}
\label{lemma:empirical-to-rademacher}
Under the assumptions of
Lemma~\ref{lemma:bernstein-in-banach-space-iid}, we have
\begin{align*}
\tfrac{1}{\numobs} \sup_{u \in \UpsSet} \inprod{u}{\sum_{i=1}^\numobs
  \X_i} \leq 3 \cdot \localrad_\numobs(\UpsSet) + 8 \sqrt{\sup_{u \in
    \UpsSet} \inprod{u}{W}^2 \cdot \tfrac{\log(\tfrac{1}{\delta})}{\numobs}} + c
\cdot \tfrac{ \log(\tfrac{1}{\delta})}{\numobs},
    \end{align*}
 with probability at least $1 - \delta$.
\end{lemma}
\noindent See Section~\ref{SecLemEmpRad} for the proof of this
claim. \\

We now use this auxiliary claim to complete the proof of
Lemma~\ref{lemma:bernstein-in-banach-space-iid}. It suffices to upper
bound the Rademacher complexity $\localrad_\numobs(\UpsSet)$. We
define the pseudometrics
\begin{align*}
\rho_*(x, y) \mydefn \sqrt{\Exs \left[ \inprod{x - y}{X_1}^2 \right]}
\quad \mbox{and} \quad \myrhotil(x, y) \mydefn \sup_{e \in \obssubset
  \cap \ball} \inprod{x - y}{e}, \qquad \mbox{for all $ x,y
  \in \vecspace^\star$.}
\end{align*}
Recalling that $R_\numobs(y) = \avg{\rade_i \inprod{y}{\X_i}}$,
applying Bernstein's inequality yields
\begin{align*}
\Prob \big[ \absVal{R_\numobs(y) - R_\numobs(z)} > t \big] \leq 2
\exp\left\{ - \min \big( \tfrac{\numobs \alpha^2}{2 \rho_* (y, z)^2},
\tfrac{\numobs \alpha}{\myrhotil(y, z)} \big) \right\} \qquad
\mbox{for any $\alpha > 0$.}
\end{align*} 
For $q \geq 1$, we let $\gamma_q$ denote the $q^{th}$-order generic
chaining functional of Talagrand.  With this notation, we have
\begin{align*}
\localrad_\numobs(\UpsSet) = \Exs \left[ \sup_{y \in \UpsSet}
  R_\numobs(y) \right] & \stackrel{(i)}{\leq} \tfrac{c}{\sqrt{\numobs}}
\cdot \gamma_2(\UpsSet, \rho_\star) + \tfrac{1}{\numobs}
\gamma_1(\UpsSet, \myrhotil) \\
    &\stackrel{(ii)}{\leq} \tfrac{c}{\sqrt{\numobs}} \cdot \Exs \left[
  \sup_{u \in \UpsSet} \inprod{u}{W} \right] + \tfrac{1}{\numobs}
\DudExp(\UpsSet, \myrhotil).
\end{align*}
Here step (i) follows from the generic chaining theorem~(see Theorem
1.2.7~\cite{talagrand2006generic}). In step (ii), we bound the first
term using the generic chaining lower bound~(see Theorem
2.1.1~\cite{talagrand2006generic}) and bound the second term using the
fact that $\gamma_1$ functional is upper bounded by the Dudley entropy
integral of order $1$. This completes the proof of
Lemma~\ref{lemma:bernstein-in-banach-space-iid}. It remains to prove
Lemma~\ref{lemma:empirical-to-rademacher}.


\subsubsection{Proof of Lemma~\ref{lemma:empirical-to-rademacher}}
\label{SecLemEmpRad}

The proof of this lemma is based on Talagrand's concentration
inequality for the suprema of empirical
process~\cite{talagrand1996new} and a symmetrization argument.  Define
the random variance $\widehat{\sigma}^2 \mydefn \tfrac{1}{\numobs}
\sup_{y \in\UpsSet} \sum_{i = 1}^\numobs \inprod{y}{\X_i}^2$.  Since
the random variables $X_i$ are bounded and $\UpsSet \subseteq
\dualBall$, we have $\absVal{\sup_{y \in \UpsSet} \; \inprod{y}{\X_i}}
\leq 1$.  Invoking Talagrand's concentration
inequality~\cite{talagrand1996new} yields the tail bound
\begin{align*}
  \Prob \Big[ \sup_{u \in \UpsSet} \inprod{u}{\Xbar_\numobs} \geq \Exs
    \big[ \sup_{u \in \UpsSet} \inprod{u}{\Xbar_\numobs} \big] +
    \alpha \Big] & \leq \exp \Big \{ \tfrac{- \numobs \alpha^2}{56
    \Exs [\widehat{\sigma}^2] + 4 \alpha} \Big \}, \quad \mbox{valid
    for all $\alpha > 0$.}
\end{align*}
Consequently, for any $\delta \in (0,1)$, we have
\begin{align*}
  \sup_{u \in \UpsSet} \inprod{u}{\Xbar_\numobs} \leq \Exs
  \big[\sup_{u \in \UpsSet} \inprod{u}{\Xbar_\numobs} \big] + 8
  \sqrt{\tfrac{\log(\tfrac{1}{\delta})}{\numobs} \Exs[\widehat{\sigma}^2]} + 4
  \cdot\tfrac{ \log(\tfrac{1}{\delta})}{\numobs}
\end{align*}
with probability at least $1 - \delta$.

It remains to upper bound the expected supremum $\Exs \big[\sup_{u \in
    \UpsSet} \inprod{u}{\Xbar_\numobs}\big]$ and the variance term
$\widehat{\sigma}^2$. By a standard symmetrization argument, we have
\begin{align*}
\Exs \left[\sup_{u \in \UpsSet} \inprod{u}{\Xbar_\numobs} \right] \leq
\tfrac{2}{\numobs} \Exs \left[\sup_{u \in \UpsSet} \sum_{i =
    1}^{\numobs} \rade_i \inprod{u}{\X_i} \right] = 2
\localrad_\numobs(\UpsSet),
\end{align*}

Moving onto the bound on $\widehat{\sigma}^2$, we have
\begin{align*}
    \widehat{\sigma}^2 \leq & \tfrac{1}{\numobs} \sup_{y \in \UpsSet}
    \sum_{i = 1}^\numobs \Big \{ \inprod{y}{\X_i}^2 - \Exs\left[
      \inprod{y}{\X_i}^2 \right] \Big \} + \tfrac{1}{\numobs} \sup_{y
      \in \UpsSet} \sum_{i = 1}^\numobs \Exs\left[ \inprod{y}{\X_i}^2
      \right] \; = \; Z_\numobs + \sup_{y \in \UpsSet} \Exs\big[
      \inprod{y}{\X_i}^2 \big],
\end{align*}
where $Z_\numobs \defn \tfrac{1}{\numobs} \sup_{y \in \UpsSet} \sum_{i
  = 1}^\numobs \left(\inprod{y}{\X_i}^2 - \Exs\left[
  \inprod{y}{\X_i}^2 \right] \right)$.  Note that each term
$|\inprod{y}{X_i}|$ is almost surely bounded by $1$, and the map $a
\mapsto a^2$ is $2$-Lipschitz over the interval $[-1, 1]$.
Consequently, letting $\{\rade_i\}_{i=1}^\numobs$ denote an
i.i.d. sequence of Rademacher variables, we have
\begin{align*}
\Exs[Z_\numobs] \; \stackrel{(i)}{\leq} \; \tfrac{2}{\numobs} \Exs
\left[ \sup_{y \in \UpsSet} \sum_{i = 1}^\numobs \rade_i
  \inprod{y}{\X_i}^2 \right] & \stackrel{(ii)}{\leq}
\tfrac{4}{\numobs} \cdot \Exs \left[ \sup_{y \in \UpsSet} \sum_{i =
    1}^\numobs \rade_i \inprod{y}{\X_i} \right]\\
& \stackrel{(iii)}{\leq} \tfrac{\numobs}{64 \log(\tfrac{1}{\delta})}
\localrad_\numobs^2(\UpsSet) + \tfrac{128}{\numobs} \log(\tfrac{1}{\delta}),
\end{align*}
where step (i) follows from a symmetrization argument; step (ii)
follows from the Ledoux--Talagrand contraction; and step (iii) follows
from the Cauchy--Schwarz inequality.  Overall, we have
\begin{align*}
  8 \sqrt{\tfrac{\log(\tfrac{1}{\delta})}{\numobs}\Exs[\widehat{\sigma}^2]} \leq
  \localrad_\numobs + 8 \sqrt{ \sup_{u \in \UpsSet} \Exs \left[
      \inprod{u}{X_1}^2 \right] \cdot \tfrac{\log(\tfrac{1}{\delta})}{\numobs} }
  + \tfrac{c \log(\tfrac{1}{\delta})}{\numobs}.
\end{align*}
Putting together the pieces yields the bound of
Lemma~\ref{lemma:empirical-to-rademacher}.


\subsection{Proof of Lemma~\ref{lemma:concentration-in-banach-space-martingale}}
\label{SecBanachMartingale}

For each vector $u \in \vecspace_*$, we define the random variable
$\Martin_\numobs(u) \mydefn \tfrac{1}{\numobs} \sum_{i = 1}^{\numobs}
\inprod{X_i}{u}$.  Clearly, the sequence $\{\Martin_t \}_{t \geq 1}$ is a
scalar martingale adapted to the filtration $(\filtration_t)_{t \geq
  0}$. Since $b_t \geq \tfrac{1}{\numobs}$, we have
\begin{align*}
\abss{\inprod{X_t}{u}} \leq b_t \cdot \sup_{x \in \ball \cap
  b_t^{-1}\obssubset} \inprod{x}{u} \leq b_t \myrhotil(u, 0)
\end{align*}
almost surely for each $t = 1, 2, \ldots$, where $\myrhotil(\cdot, \cdot)$ is a pseudo-metric on the dual space $\vecspace_*$ defined in~\eqref{eq:defn-pseudo-metric-dual}. For any $u_1, u_2
\in \vecspace^*$, the Azuma-Hoeffding inequality implies that
\begin{align*}
\Prob \Big[ \abss{\Martin_\numobs (u_1) - \Martin_\numobs (u_2)} \geq \alpha \Big]
& \leq \exp \big \{ - \tfrac{\numobs \alpha^2}{\myrhotil(u_1, u_2)^2
  \sum_{i = 1}^\numobs b_i^2} \big \} \qquad \mbox{for each $\alpha >
  0$.}
\end{align*}
Applying the Dudley chaining tail bound (see
e.g.~\cite{van2014probability}, Theorem 5.29) to the sub-Gaussian
process $\{\Martin_\numobs (u)\}_{u \in \dualBall}$, there exist
universal constants $c, c_1 > 0$ such that
\begin{align*}
\Prob \left[ \sup_{u \in \dualBall} \Martin_\numobs (u) \geq c
  \sqrt{\sum_{i = 1}^\numobs b_i^2} \cdot \Big( \int_0^1 \sqrt{\log
    N(s; \dualBall, \myrhotil)} ds + t \Big) \right] \leq c e^{- c_1
  t^2} \quad \mbox{for each $t > 0$.}
\end{align*}
Setting $t = \sqrt{c_1^{-1}\log(1 / \delta)}$ yields the claim.


\section{Proofs of Auxiliary Lemmas}
\label{sec:proof-of-auxiliary-lemmas}
In this section we prove various auxiliary Lemmas that we use
throughout the main proof Section~\ref{sec:Main-proofs}.

\subsection{Proof of Lemma~\ref{lemma:coarse-bound-gronwall}} 
\label{subsubsec:proof-lemma-coarse-bound-gronwall}

From the recursive relation, we have the upper bounds
\mbox{$\norm{\theta_t - \thetastar} \leq \norm{\theta_{t - 1} -
    \thetastar} + \stepsize \norm{\myV_{t}}$,} as well as
\begin{align*}
\norm{\myV_t} & \leq \tfrac{t - 1}{t} \norm{\myV_{t - 1}} + \tfrac{1}{t}
\Big \{ \norm{\theta_{t - 1} - \theta_{t - 2}} + \norm{\Hstoch_t
  (\theta_{t - 1}) - \Hstoch_t (\theta_{t - 2})} \Big \} +
\tfrac{1}{t} \norm{\Hstoch_t (\theta_{t - 1}) - \theta_{t - 1}} \\
& \leq \Big \{ 1 + \tfrac{\stepsize (\Lip + 1)}{t} \Big \} \norm{\myV_{t
    - 1}} + \tfrac{2\Lip}{t} \norm{\theta_{t - 1} - \thetastar} +
\tfrac{1}{t} \boundstar.
\end{align*}
Putting these two inequalities together yields the vector-based
recursion
\begin{align*}
    \begin{bmatrix}
     \norm{\theta_t - \thetastar}\\ \norm{\myV_t}
    \end{bmatrix} \leq \begin{bmatrix}
     1 & \stepsize\\ 0 &1
    \end{bmatrix} \cdot \left(  \begin{bmatrix}
     1 & 0\\ \tfrac{2 \Lip}{t} & 1 + \tfrac{\stepsize (\Lip + 1)}{t}
    \end{bmatrix} \begin{bmatrix}
     \norm{\theta_{t - 1} - \thetastar}\\ \norm{\myV_{t - 1}}
    \end{bmatrix} + \begin{bmatrix}
     0\\ \tfrac{1}{t} \boundstar
    \end{bmatrix} \right),
\end{align*}
where the inequality is taken elementwise.  Solving this vector
recursion yields
\begin{align*}
\norm{\theta_t - \thetastar} + \norm{\myV_t} \leq e^{1 + \stepsize
  \Lip t} \left( \boundstar + \norm{\theta_0 - \thetastar} \right),
\end{align*}
valid for any $t \geq \burnin \geq \tfrac{1}{\stepsize}$.


\subsection{Proof of Lemma~\ref{lemma:burnin-period-bound}}
\label{sec:proof-of-lemma:burnin-period-bound}
By definition, we have
\begin{align*}
\myV_{\burnin} = \tfrac{1}{\burnin} \sum_{t = 1}^{\burnin} \big(
\Hstoch_t (\thetainit) - \thetainit \big) = \big( \hpop(\thetainit) -
\thetainit \big) + \tfrac{1}{\burnin} \sum_{t = 1}^{\burnin} \noise_t
(\thetastar) + \tfrac{1}{\burnin} \sum_{t = 1}^{\burnin} \big(
\noise_t(\thetainit) - \noise_t(\thetastar) \big).
\end{align*}
Lemma~\ref{lemma:bernstein-in-banach-space-iid} guarantees that
\begin{align*}
\|\sum_{t = 1}^{\burnin} \noise_t(\thetastar)\| \leq
      {c}{\sqrt{\burnin}} \Big \{ \noisegauss + \noisevar \sqrt{\log
        (\tfrac{1}{\delta})} \Big \} + c \big \{ \DudExp(\dualBall, \myrhotil) +
      \log(\tfrac{1}{\delta}) \big \},
\end{align*}
with probability $1 - \delta$.  Moreover, for each integer $t \in
[\burnin]$, we have:
\begin{align*}
\norm{\noise_t(\thetainit) - \noise_t(\thetastar)} \leq \Lip
\norm{\thetainit - \thetastar} \leq \tfrac{\Lip}{1 - \contraction}
\norm{\hpop(\thetainit) - \thetainit}.
\end{align*}
Lemma~\ref{lemma:concentration-in-banach-space-martingale} implies
that
\begin{align*}
\|\sum_{t = 1}^{\burnin} \big( \noise_t(\thetainit) - \noise_t
(\thetastar) \big) \| & \leq \tfrac{c \Lip}{1 - \contraction}
\norm{\hpop(\thetainit) - \thetainit} \sqrt{\burnin} \: \Big \{
\constMGtype + \sqrt{\log(\tfrac{1}{\delta})} \Big \}
\end{align*}
with probability at least $1 - \delta$.

By combining these bounds, we find that
\begin{multline*}
\norm{\myV_{\burnin}} \leq \norm{\hpop(\thetainit) - \thetainit} \Big
\{ 1 + \tfrac{c \Lip}{(1 - \contraction) \sqrt{\burnin}} \big[
  \constMGtype + \sqrt{\log(\tfrac{1}{\delta})} \big] \Big \} \\
+ \tfrac{c}{\sqrt{\burnin}} \big( \noisegauss + \noisevar \sqrt{\log
  (\tfrac{1}{\delta})} \big) + \tfrac{c}{\burnin} \Big \{
\DudExp(\dualBall, \myrhotil) + \log(\tfrac{1}{\delta}) \Big \}
\end{multline*}
with probability at least $1 - \delta$.  Substituting the burn-in time
bound~\eqref{eq:burn-in-requirement-in-bootstrap-phase-i} yields the
final claim.

\section{Comments on Theorem~\ref{new--thm:main-nonlinear-general}}

In Section~\ref{sec:localLinQProof}, we prove that the Bellman
optimality operator associated with the optimal $Q$-function estimation
problem satisfies the local linearity
condition~\ref{new--assume:linearization}. Using a similar argument,
in Section~\ref{sec:localLinSSPProof} we show that the Bellman fixed-point operator 
for the stochastic shortest path problem 
satisfies the local linearity condition. 


\subsection{Verifying local linearity for Bellman optimality 
operator}
\label{sec:localLinQProof}

In this section, we verify that the local linearity
assumption~\ref{new--assume:linearization} holds for the Bellman
optimality operator for
$Q$-learning~\cite{watkins1992q,szepesvari1998asymptotic,wainwright2019variance}.
Consider a tabular MDP \mbox{$\MDP = (\rewardQ, \TransMat,
  \contraction)$} with state space $\statesQ$ and action space
$\actionsQ$. For any state-action pair $(\state, \action) \in \statesQ
\times \actionsQ$, the scalar $\rewardQ(\state, \action)$ denotes the
reward when the action $\action$ is taken at state $\state$, and the
scalar $\TransMat_{\action}(\state' \mid \state)$ denotes the
probability of transitioning to state $\state'$ when the action
$\action$ is chosen at state $\state$.

One way to estimate an optimal policy is to calculate the optimal
$Q$-function. Associated with a (deterministic) policy $\pi$ is its
\emph{$Q$-function}
\begin{align*}
\theta^{\pi}(\state, \action) & \defn \Exs \Big[\sum_{k=0}^{\infty}
  \reward(\state_k, \action_k) \mid \state_0 = \state, \action_0 =
  \action \Big], \quad \text { where } \action_k = \policy(\state_k)
\quad \mbox{for all $k = 1, 2, \ldots$.}
\end{align*}
The optimal $Q$-function is given by $\theta^\star(\state, \action)
\defn \sup_{\pi \in \Pi}\; \theta^{\pi}(\state, \action)$, and an
optimal policy can be obtained as $\pi_{\star}(\state) = \arg
\max_{\action} \Qstar(\state, \action)$.

The Bellman optimality operator $\QOp$ acts on the space of
$Q$-functions; more precisely, its action on a given $Q$-function $\Q$
is given by
\begin{align}
\label{eqn:Bellman-Q}
\QOp(\Q)(\state, \action) = \rewardQ(\state, \action) + \contraction
\sum_{\state'} \TransMat_{\action}(\state' \mid \state) \cdot
\max_{\action'} \Q(\state', \action') \quad \text{for all } (\state,
\action) \in \statesQ \times \actionsQ.
\end{align}
By standard results~\cite{bertsekas2012approximate}, the operator
$\QOp$ is $\contraction$-contractive in the $\ell_\infty$-norm, and
the optimal state-action value function $\Qstar$ is its unique fixed
point.

For a given $Q$-function $\Q$, the associated greedy policy $\pi_{\Q}$
is given by
\begin{align}
  \pi_{\Q}(\state) = \arg\max_{\action} \Q(\state,
\action),
\end{align}
 where we break any ties by taking the smallest action (in the enumeration order) that achieves the maximum.  Using this greedy
policy, we can define the right linear operator
\begin{align*}
  \TransMat^{\pi_{\Q}} \Q(\state, \action) = \sum_{\state'}
  \TransMat_{\action}(\state' \mid \state) \Q(\state',
  \pi_{\Q}(\state')).
\end{align*}

Let $\ball(\thetastar, s) \defn \{ \theta \mid \|\theta -
\thetastar\|_\infty \leq s \}$ denote the $\ell_{\infty}$-ball of
radius $s$ around $\thetastar$, and define the set
\begin{align}
\mathcal{A}_{\delrad} = \{ \contraction \cdot \TransMat^{\pi_\Q} \;
\mid \; \pi_\Q \;\; \text{is a greedy policy of } \Q \text{ with } \Q
\in \ball(\Qstar, \delrad)\}
\end{align}
of linear operators.  We use $\pi_\star$ to denote the greedy policy
associated with the optimal $Q$-function $\Qstar$. By definition, the
$Q$-functions $\theta$ and $\thetastar$ satisfy the fixed-point
relations
\begin{align*}
  \QOp(\theta) = \rewardQ + \AQ \Q \qquad \text{and} \qquad \Qstar =
  \rewardQ + \Astar \Qstar.
\end{align*}
Rearranging the last two equations yields
\begin{subequations}
\label{eqn:BellmanDefectQ}    
\begin{align}
\QOp(\Q) - \Q & = \rewardQ + \AQ \Q - \Q = ( \Id - \Astar)(\Qstar -
\Q) + (\AQ - \Astar)\Q \\
\QOp(\Q) - \Q &= \rewardQ + \AQ \Q - \Q = (
\Id - \AQ)(\Qstar - \Q) + (\AQ - \Astar)\Qstar.
\end{align}
\end{subequations}
Next we claim that
\begin{align}
\label{eqn:elementwiseBound}
    (\Id - \Astar)^{-1} (\AQ - \Astar)\Q \stackrel{(a)}{\succcurlyeq}
0 \qquad \text{and}
    \qquad (\Id - \AQ)^{-1} (\AQ - \Astar) \Qstar
    \stackrel{(b)}{\preccurlyeq} 0.
\end{align}
Indeed, since the policy $\pi_\theta$ is greedy for $\Q$, we have the
element-wise inequality $(\AQ - \Astar)\Q \succcurlyeq 0$.  The matrix
$(\Id - \Astar)^{-1}$ has non-negative entries, so that element-wise
inequality (a) holds.  A similar argument, using the fact that
$\pi_\star$ is greedy for $\Qstar$, yields the element-wise (b).

With the last observation in hand, combining the element-wise
inequalities~\eqref{eqn:elementwiseBound} with the two expressions of
the Bellman defect~\eqref{eqn:BellmanDefectQ} yields
\begin{align*}
 |\Qstar - \Q| \preccurlyeq \max\{ |(\Id - \AQ)^{-1}(\QOp(\Q) - \Q)|,
 \; |(\Id - \Astar)^{-1}(\QOp(\Q) - \Q)|\}.
\end{align*}
Finally, note that 
the operator $\AQ \in \mathcal{A}_\delrad$ for
any $\Q \in \ball(\Qstar, \delrad)$. Putting together the pieces 
we conclude that for all 
$\theta \in \ball(\thetastar, \delrad)$
\begin{align*}
  \|\theta - \thetastar\|_\infty
  \leq \sup_{A \in \mathcal{A}_{\delrad}}
  \| (\Id -  A)^{-1} (\hpop(\theta) - \theta) \| 
\end{align*}
Thus, we deduce that the local linearity
condition~\ref{new--assume:linearization} is satisfied for the Bellman
optimality operator $\hpop$ from equation~\eqref{eqn:Bellman-Q} with
$\anorm{\cdot} = \|\cdot\|_\infty$.


\subsection{Verifying local linearity  for the SSP operator}
\label{sec:localLinSSPProof}

Recall from Section~\ref{new--sec:SSP} the definition of a stochastic
shortest path (SSP) problem $(\rewardQ, \TransMat)$ with optimal-$Q$
value $\Qstar$.  For a given $Q$-function $\Q$, consider the greedy
policy $\Pi_{\Q}(\state) = \arg\min_{\action} \Q(\state, \action)$.
We can use it to define the right linear operator
$\TransMat^{\pi_{\Q}} \Q(\state, \action) = \sum_{\state'}
\TransMat_{\action}(\state' \mid \state) \Q(\state',
\pi_{\Q}(\state'))$.  Letting $\ball(\thetastar, s) \defn \{ \theta
\mid \|\theta - \thetastar\| \leq s \}$ denote the
$\ell_{\infty}$-ball of radius $s$ around $\thetastar$, we define the
set
\begin{align}
\label{eqn:A_delta-SSP}
\mathcal{A}_{\delrad} = \{ \TransMat^{\pi_\Q} \; \mid \; \pi_\Q \;\;
\text{is a greedy policy of } \Q \text{ with } \Q \in \ball(\Qstar,
\delrad)\}
\end{align}
of linear operators.  We use $\pi_\star$ to denote the greedy policy
associated with the optimal $Q$-function $\Qstar$.

With this setup in hand, following the same argument as
Section~\ref{sec:localLinQProof}, the local linearity
assumption~\ref{new--assume:linearization} for the Bellman
operator~\eqref{eqn:SSP-Bellman} can be verified with the set
$\mathcal{A}_{\delrad}$ of local linear operators defined in
equation~\eqref{eqn:A_delta-SSP}, and with $\anorm{\cdot} = \|
\cdot\|_\infty$.

\end{document}

%% file: final_macros.tex
\newcommand{\SAstep}{\alpha}

\newcommand{\abs}[1]{|#1|}

\newcommand{\X}{X}
\newcommand{\Xbar}{\widebar{X}}
\newcommand{\avg}[1]{\frac{1}{\numobs}\sum_{i = 1}^\numobs #1}
\newcommand{\absVal}[1]{|#1|}

\newcommand{\dir}{v}

\newcommand{\anorm}[1]{\big\| #1 \big\|_{\mbox{\tiny{\MySet}}}}

\newcommand{\aUB}{\mathcal{D}}

\newcommand{\termOne}{T_1}
\newcommand{\termTwo}{T_2}
\newcommand{\cost}{c}

\newcommand{\nonAbsStates}{\states_{- 1}}

\newcommand{\SSPOp}{\hpop}
\newcommand{\wNorm}[1]{\| #1 \|_{w}}
\newcommand{\NormWeight}{\mathbf{w}}
\newcommand{\SSPdim}{D}
\newcommand{\TranSample}{\mathbf{Z}}
\newcommand{\costRnd}{C}
\newcommand{\SSPstoch}{\Hstoch}

\newcommand{\actionsOne}{\actions_1}
\newcommand{\actionsTwo}{\actions_2}
\newcommand{\TranMG}{\mathbf{P}}

\newcommand{\randReward}{R}

\newcommand{\ProbDistr}[1]{\mathcal{P}(#1)}
\newcommand{\valueMG}{V}
\newcommand{\policySet}{\Pi}
\newcommand{\QfunMG}{\theta^\star}
\newcommand{\MGOp}{\hpop}

\newcommand{\Event}{\mathscr{E}}
\newcommand{\bigthetaEvent}{\Event^{(\theta)}_n (\ratetheta)}
\newcommand{\bigvEvent}{\Event^{(v)}_n (\ratev)}

\newcommand{\numRestarts}{R}

\newcommand{\dualBall}{\Gamma}

\newcommand{\MDP}{M}

\newcommand{\states}{\mathcal{X}}
\newcommand{\state}{x}
\newcommand{\actions}{\mathcal{U}}
\newcommand{\action}{u}
\newcommand{\TransMat}{\mathbf{P}}

\newcommand{\policy}{\pi}

\newcommand{\filtration}{\mathcal{F}}

\newcommand{\dudley}{\mathcal{J}}

\newcommand{\succeqOrthant}{\succeq_{\mathrm{orth}}}
\newcommand{\preceqOrthant}{\preceq_{\mathrm{orth}}}

\usepackage{thmtools, thm-restate}

\newcommand{\stepsize}{\alpha}

\newcommand{\localvar}{\nu}
\newcommand{\localrad}{\mathscr{R}}
\newcommand{\localgauss}{\mathcal{G}}
\newcommand{\noisegauss}{\widebar{\mathscr{W}}}
\newcommand{\noisevar}{\nu}

\newcommand{\ratev}{r_v}
\newcommand{\ratetheta}{r_\theta}
\newcommand{\highorder}{\mathcal{H}}

\newcommand{\numobs}{n}

\newcommand{\bvec}{\ensuremath{b}}

\newcommand{\Amat}{\ensuremath{A}}

\newcommand{\reward}{\ensuremath{r}}

\newcommand{\usedim}{\ensuremath{d}}

\newcommand{\Sphere}[1]{\ensuremath{\mathbb{S}}}

\def\beq{\begin{equation} }
\def\eeq{\end{equation} }

\newcommand{\burnin}{{B_0}}
\newcommand{\constMGtype}{{\DudGauss(\ball^*, \myrhotil)}}

\newcommand{\Lip}{L}

\newcommand{\hpop}{\mathbf{h}}
\newcommand{\Hstoch}{\mathbf{H}}

\newcommand{\linearizationSet}{\mathcal{A}}

\newcommand{\contraction}{\gamma}

\newcommand{\norm}[1]{\vecnorm{#1}{}}
\newcommand{\dualnorm}[1]{\vecnorm{#1}{*}}

\newcommand{\vecspace}{\mathbb{V}}
\newcommand{\rade}{\zeta}
\newcommand{\obssubset}{\Omega}
\newcommand{\boundstar}{b_*}

\newcommand{\DualSpace}{\vecspace^\star}
\newcommand{\LipCon}{L}

\newcommand{\bracketMed}[1]{\left\lbrace #1 \right\rbrace}

\newcommand{\rootSA}{\texttt{ROOT-SA}$\:$}
\newcommand{\rootSANoSpace}{\texttt{ROOT-SA}}
\newcommand{\Seq}[1]{\bracketMed{#1}_{t \geq 1}}
\newcommand{\thetainit}{\theta_0}

\newcommand{\compostep}{m}

\newcommand{\GaussNoise}{W}
\newcommand{\Id}{\mathcal{I}}

\newcommand{\optPolicy}{\policy^\star}

\newcommand{\dims}{\ensuremath{d}}

\newcommand{\real}{\ensuremath{\mathbb{R}}}

\newcommand{\order}[1]{\ensuremath{\mathcal{O}\parenth{#1}}}

\newcommand{\thetastar}{\ensuremath{{\theta^*}}}

\newcommand{\parenth}[1]{\left( #1 \right)}

\newcommand{\abss}[1]{\left| #1 \right |}

\newcommand{\ball}{\ensuremath{\mathbb{B}}}

\newcommand{\mydefn}{\ensuremath{:=}}

\newcommand{\noise}{\ensuremath{\varepsilon}}

\newcommand{\defn}{:=}

\newcommand{\matsnorm}[2]{|\!|\!| #1 | \! | \!|_{{#2}}}
\newcommand{\vecnorm}[2]{\| #1 \|_{#2}}

\newcommand{\opnorm}[1]{\ensuremath{\matsnorm{#1}{\tiny{\mbox{op}}}}}

\newcommand{\inprod}[2]{\ensuremath{\langle #1 , \, #2 \rangle}}

\newcommand{\kull}[2]{\ensuremath{D_{\text{KL}}(#1\; \| \; #2)}}

\newcommand{\Exs}{\ensuremath{{\mathbb{E}}}}
\newcommand{\Prob}{\ensuremath{{\mathbb{P}}}}

\newtheoremstyle{named}{}{}{\itshape}{}{\bfseries}{.}{.5em}{\thmnote{#3's }#1}
\theoremstyle{named}

\theoremstyle{plain}

\newtheorem{theorem}{Theorem}

\newtheorem{lemma}{Lemma}

\newtheorem{corollary}{Corollary}

\newlength{\widebarargwidth}
\newlength{\widebarargheight}
\newlength{\widebarargdepth}
\DeclareRobustCommand{\widebar}[1]{%
  \settowidth{\widebarargwidth}{\ensuremath{#1}}%
  \settoheight{\widebarargheight}{\ensuremath{#1}}%
  \settodepth{\widebarargdepth}{\ensuremath{#1}}%
  \addtolength{\widebarargwidth}{-0.3\widebarargheight}%
  \addtolength{\widebarargwidth}{-0.3\widebarargdepth}%
  \makebox[0pt][l]{\hspace{0.3\widebarargheight}%
    \hspace{0.3\widebarargdepth}%
    \addtolength{\widebarargheight}{0.3ex}%
    \rule[\widebarargheight]{0.95\widebarargwidth}{0.1ex}}%
  {#1}}

\makeatletter
\long\def\@makecaption#1#2{
        \vskip 0.8ex
        \setbox\@tempboxa\hbox{\small {\bf #1:} #2}
        \parindent 1.5em  
        \dimen0=\hsize
        \advance\dimen0 by -3em
        \ifdim \wd\@tempboxa >\dimen0
                \hbox to \hsize{
                        \parindent 0em
                        \hfil
                        \parbox{\dimen0}{\def\baselinestretch{0.96}\small
                                {\bf #1.} #2
                                }
                        \hfil}
        \else \hbox to \hsize{\hfil \box\@tempboxa \hfil}
        \fi
        }
\makeatother

\long\def\comment#1{}
\definecolor{battleshipgrey}{rgb}{0.52, 0.52, 0.51}
\definecolor{darkgray}{rgb}{0.66, 0.66, 0.66}
\definecolor{darkgreen}{rgb}{0.0, 0.2, 0.13}
\definecolor{darkspringgreen}{rgb}{0.09, 0.45, 0.27}
\definecolor{dukeblue}{rgb}{0.0, 0.0, 0.61}
\definecolor{olivedrab7}{rgb}{0.24, 0.2, 0.12}
\definecolor{darkblue}{rgb}{0.0, 0.0, 0.55}
\definecolor{darkscarlet}{rgb}{0.34, 0.01, 0.1}
\definecolor{candyapplered}{rgb}{1.0, 0.03, 0.0}
\definecolor{ao(english)}{rgb}{0.0, 0.5, 0.0}
\definecolor{applegreen}{rgb}{0.55, 0.71, 0.0}
